\documentclass[a4paper, 11pt]{amsart} 
\usepackage{enumitem, amsmath, amssymb, amsthm, xspace, caption, subcaption, comment, amsaddr, mathtools}
\usepackage{tikz}
\usepackage{fullpage}
\usepackage[hidelinks]{hyperref} 
\usepackage{stmaryrd}
\usepackage{xcolor}
\usepackage[normalem]{ulem} 
\usepackage{blkarray} 

\usetikzlibrary{calc}
\usetikzlibrary{decorations.markings}

\newif\ifdraft
\draftfalse

\ifdraft 
    \usepackage[notcite,notref]{showkeys}

  \else
    \renewcommand{\sout}[1]{} 
\fi 

\newtheorem{thm}{Theorem}[section]
\newtheorem*{thm*}{Theorem}
\newtheorem{lem}[thm]{Lemma}
\newtheorem{prop}[thm]{Proposition} 

\newtheorem{conj}[thm]{Conjecture}
\newtheorem{qu}[thm]{Question}
\newtheorem{coro}[thm]{Corollary}
\newtheorem{claim}[thm]{Claim}
\newtheorem{fact}[thm]{Fact}

\newtheorem{remark}[thm]{Remark}

\newtheorem{defi}[thm]{Definition}
\theoremstyle{definition}
\newtheorem*{claim*}{Claim}
\theoremstyle{plain}

\numberwithin{figure}{section}

\newenvironment{romenumerate}[1][-0.25em]{
\vspace{#1}
\begin{enumerate}
\itemsep0pt \parskip0pt \parsep0pt%
 }
 {\end{enumerate} }

\DeclareMathOperator*{\esssup}{ess\,sup}
\DeclareMathOperator*{\essinf}{ess\,inf}

\makeatletter

\newcommand{\By}[2]{\overset{\mbox{\tiny{#1}}}{#2}}
\newcommand{\ByRef}[2]{   \By{\eqref{#1}}{#2} }

\newcommand{\geBy}[1]{    \By{#1}{\ge} }

\newcommand{\leByRef}[1]{ \ByRef{#1}{\le} }

\newcommand{\justify}[1]{\fbox{\tiny{#1}}\quad}

    \newcommand{\trian}{
    \begin{tikzpicture}[baseline=-0.3ex,scale=0.25]
    \tikzstyle{vertex}=[circle,fill=black, minimum size=2pt,inner sep=1pt]
    \node[vertex] (v1) at (-0.5, 0){};
    \node[vertex] (v2) at (0.5,0){};
    \node[vertex] (v3) at (0,0.8){};
    \draw (v1)--(v2)--(v3)--(v1);
    \end{tikzpicture}}
    
    \newcommand{\cherry}{
    \begin{tikzpicture}[baseline=-0.3ex,scale=0.25]
    \tikzstyle{vertex}=[circle,fill=black, minimum size=2pt,inner sep=1pt]
    \node[vertex] (v1) at (-0.5, 0){};
    \node[vertex] (v2) at (0.5,0){};
    \node[vertex] (v3) at (0,0.8){};
    \draw (v1)--(v3)--(v2);
    \end{tikzpicture}}

    \newcommand{\cocherry}{
    \begin{tikzpicture}[baseline=-0.3ex,scale=0.25]
    \tikzstyle{vertex}=[circle,fill=black, minimum size=2pt,inner sep=1pt]
    \node[vertex] (v1) at (-0.5, 0){};
    \node[vertex] (v2) at (0.5,0){};
    \node[vertex] (v3) at (0,0.8){};
    \draw (v1)--(v2);
    \end{tikzpicture}}
    
    \newcommand{\cocherrya}{
    \begin{tikzpicture}[baseline=-0.3ex,scale=0.25]
    \tikzstyle{vertex}=[circle,fill=black, minimum size=2pt,inner sep=1pt]
    \node[vertex] (v1) at (-0.5, 0){};
    \node[vertex] (v2) at (0.5,0){};
    \node[vertex] (v3) at (0,0.8){};
    \draw (v1)--(v3);
    \end{tikzpicture}}

    \newcommand{\cotrian}{
    \begin{tikzpicture}[baseline=-0.3ex,scale=0.25]
    \tikzstyle{vertex}=[circle,fill=black, minimum size=2pt,inner sep=1pt]
    \node[vertex] (v1) at (-0.5, 0){};
    \node[vertex] (v2) at (0.5,0){};
    \node[vertex] (v3) at (0,0.8){};
    \end{tikzpicture}}

\newcommand{\edge}{
  \begin{tikzpicture}[baseline=-0.6ex,scale=0.25]
  \tikzstyle{vertex}=[circle,fill=black, minimum size=2pt,inner sep=1pt]
  \node[vertex] (v1) at (0, -0.5){};
  \node[vertex] (v2) at (0,0.5){};
  \draw (v1)--(v2);
  \end{tikzpicture}}

\newcommand{\coedge}{
  \begin{tikzpicture}[baseline=-0.6ex,scale=0.25]
  \tikzstyle{vertex}=[circle,fill=black, minimum size=2pt,inner sep=1pt]
  \node[vertex] (v1) at (0, -0.5){};
  \node[vertex] (v2) at (0,0.5){};
  \end{tikzpicture}}

\title{From flip processes to dynamical systems on graphons} 
\author{Frederik Garbe}
\address{Faculty of Informatics, Masaryk University, Brno, Czechia}
\author{Jan Hladk\'y, Matas \v{S}ileikis}
\address{Institute of Computer Science of the Czech Academy of Sciences, Prague, Czechia} 
\author{Fiona Skerman}
\address{Department of Mathematics, Uppsala University, Uppsala, Sweden} 
\email{garbe@fi.muni.cz, hladky@cs.cas.cz, matas@cs.cas.cz, fiona.skerman@math.uu.se}

\thanks{FG: Parts of the work were done while affiliated with the Institute of Mathematics of the Czech Academy of Sciences. Research supported by GA\v{C}R project 18-01472Y and RVO: 67985840. Furthermore, this work has received funding from the European Research Council (ERC)
	under the European Union's Horizon 2020 research and innovation
	programme (grant agreement No 648509). This publication reflects only
	its authors' view; the European Research Council Executive Agency is not
	responsible for any use that may be made of the information it contains.
	This author was also supported by the MUNI Award in Science and
	Humanities of the Grant Agency of Masaryk University. JH: Research supported by Czech Science Foundation Project 21-21762X. M\v{S}: Research supported by Czech Science Foundation Project 20‐27757Y, with Institutional Support RVO:67985807. FS: Research supported by the Wallenberg AI, Autonomous Systems and Software
	Program (WASP), the project AI4Research at Uppsala University, and a Simons-Berkeley Research Fellowship.
}

\begin{document}

\newcommand{\sm}{\setminus}
\newcommand{\co}[1]{\overline{#1}}
\newcommand{\lgr}[1]{\mathcal{H}_{#1}}
\newcommand{\tlgr}[1]{\mathcal{F}_{#1}^{\bullet\bullet}}
\newcommand{\rul}{\mathcal{R}}    
\newcommand{\Z}{\mathbb{Z}}    
\newcommand{\R}{\mathbb{R}}    
\newcommand{\N}{\mathbb{N}}    
\newcommand{\Q}{\mathbb{Q}}    
\newcommand{\cf}{\mathcal{F}}    
\newcommand{\ca}{\mathcal{A}}
\newcommand{\cb}{\mathcal{B}}
\newcommand{\cs}{\mathcal{S}}
\newcommand{\eps}{\varepsilon}
\newcommand{\vphi}{\varphi}
\newcommand{\Prob}{\mathbb{P}}
\newcommand{\prob}[1]{\Prob\left( #1 \right)}
\newcommand{\Pc}[2]{\Prob\left(#1\,|\,#2\right)}
\newcommand{\E}{\mathbb{E}}
\newcommand{\Ec}[2]{\E \left(#1\,|\,#2  \right)}
\newcommand{\tind}{{t_{\mathrm{ind}}}}
\newcommand{\tr}[1]{t^{#1}}
\newcommand{\tindr}[1]{t_{\mathrm{ind}}^{#1}}
\newcommand{\Troot}[2]{T_{#1}^{#2}}
\newcommand{\cutn}[1]{\left\lVert #1\right\rVert_{\square}}
\newcommand{\cutm}{\delta_\square}
\newcommand{\cutnd}{d_\square}
\newcommand{\Lone}[1]{\left\lVert #1\right\rVert_{1}}
\newcommand{\Linf}[1]{\left\lVert #1\right\rVert_{\infty}}
\newcommand{\tocutn}{\xrightarrow{\square}}
\newcommand{\toLinf}{\xrightarrow{\infty}}
\newcommand{\vel}[1][\:]{\mathfrak{V}_{#1}} 
\newcommand{\traj}[2]{\Phi^{#1}#2}
\newcommand{\trajRule}[3]{\Phi_{#1}^{#2}#3}
\newcommand{\vv}{\mathbf{v}}
\newcommand{\x}{\mathbf{x}}
\newcommand{\HH}{\mathbf{H}}
\newcommand{\D}{\,\mathrm{d}}
\newcommand{\indic}{\mathbb{I}}
\newcommand{\ind}[1]{\indic_{\left\{ #1 \right\}}}
\newcommand{\Kernel}{\mathcal{W}}
\newcommand{\Gra}[1][]{{\Kernel_{0}^{#1}}} 
\newcommand{\sign}{\operatorname{sign}}
\newcommand{\G}{\mathbb{G}}
\newcommand{\hh}{\mathbf{H}}
\newcommand{\ff}{\mathbf{F}}
\newcommand{\uu}{\mathbf{U}}
\newcommand{\gb}{\mathbf{G}}
\newcommand{\bik}{\binom{k}{2}}
\newcommand{\flo}[1]{\left\lfloor #1 \right\rfloor}
\newcommand{\cei}[1]{\left\lceil #1 \right\rceil}
\newcommand{\dest}[1][]{ { \mathsf{dest}_{\,#1} } }
\newcommand{\source}[1][]{ { \mathsf{orig}_{\,#1} } }
\newcommand{\sdiff}{\ominus}
\newcommand{\norm}[1]{\left\lVert #1 \right\rVert}
\newcommand{\mesh}{\operatorname{mesh}}
\newcommand{\uppdev}{{\overline D}}
\newcommand{\kruh}{\mathsf{circle}}
\newcommand{\age}{\mathsf{age}}
\newcommand{\mdom}[1]{\mathcal{D}_{#1}}
\newcommand{\life}[1]{I_{#1}}

\begin{abstract}
We introduce a class of random graph processes, which we call \emph{flip processes}. Each such process is given by a \emph{rule} which is a function $\mathcal{R}:\mathcal{H}_k\rightarrow \mathcal{H}_k$ from all labeled $k$-vertex graphs into itself ($k$ is fixed). The process starts with a given $n$-vertex graph $G_0$. In each step, the graph $G_i$ is obtained by sampling $k$ random vertices $v_1,\ldots,v_k$ of $G_{i-1}$ and replacing the induced graph $F:=G_{i-1}[v_1,\ldots,v_k]$ by  $\mathcal{R}(F)$. This class contains several previously studied processes including the Erd\H{o}s--R\'enyi random graph process and the triangle removal process. Actually, our definition of flip processes is more general, in that $\mathcal{R}(F)$ is a probability distribution on $\mathcal{H}_k$, thus allowing randomised replacements.

Given a flip process with a rule $\mathcal{R}$, we construct time-indexed trajectories $\Phi:\mathcal{W}_0\times [0,\infty)\rightarrow\mathcal{W}_0$ in the space of graphons. We prove that for any $T > 0$ starting with a large finite graph $G_0$ which is close to a graphon $W_0$ in the cut norm, with high probability the flip process will stay in a thin sausage around the trajectory $(\Phi(W_0,t))_{t=0}^T$ (after rescaling the time by the square of the order of the graph).

These graphon trajectories are then studied from the perspective of dynamical systems. Among others topics, we study continuity properties of these trajectories with respect to time and initial graphon, existence and stability of fixed points and speed of convergence (whenever the infinite time limit exists). We give an example of a flip process with a periodic trajectory.
\end{abstract}
\maketitle
\linespread{1}
\newpage  \tableofcontents
\linespread{1.1}
\section{Introduction}\label{sec:intro}
\subsection{The Erd\texorpdfstring{\H{o}}{\"o}s--R\'enyi random graph process}\label{ssec:introER}
The uniform Erd\H{o}s--R\'enyi random graph model $\G(n,M)$ is 63 years old and  now the theory of random graphs is an active field at the border between graph theory and probability. While studying this model provides insights in many situations, sometimes one is in an evolving setting, of which $\G(n,M)$ is just one snapshot. More precisely, the \emph{Erd\H{o}s--R\'enyi random graph process}, for a given integer $n$, is a random sequence of graphs $G_0,G_1,\ldots,G_{\binom{n}2}$, all on vertex set $[n]$, in which $G_0$ is edgeless and each $G_{\ell+1}$ is created from $G_{\ell}$ by flipping a uniformly random non-edge into an edge. Lastly, $G_{\binom{n}2}$ is a complete graph and so there is no meaningful way to continue the process. 

It is easy to prove that as $n$ goes to infinity, all the graphs $G_0,G_1,\ldots,G_{\binom{n}2}$ are $o(1)$-quasirandom, asymptotically almost surely. Here and throughout the paper, we will use the concept of quasirandomness which was developed in the 1980s by Chung, Graham and Wilson~\cite{Chung1989} and by Thomason~\cite{ThomasonPseudorandom} and reflected the then emerging theory surrounding the Szemer\'edi regularity lemma~\cite{Sze:ReguLemma}. Namely, we say that an $n$ vertex graph is \emph{$\eps$-quasirandom of density $d$} if it has $d\binom{n}2$ edges and each set of vertices $U$ induces $\binom{|U|}{2} d \pm \eps n^2$ edges.

Let us consider a small modification of this process, which we call the \emph{Erd\H{o}s--R\'enyi flip process}. We again start with the edgeless graph on $[n]$. To create $G_{\ell+1}$ from $G_\ell$, we pick a uniformly random pair $xy$ of vertices and make that pair an edge. In particular, if $xy$ was an edge, then $G_{\ell+1}=G_{\ell}$. Let us make some remarks about the Erd\H{o}s--R\'enyi flip process. Firstly, this process is defined for any number of steps, and eventually stabilises at the complete graph (around time $\tau_\mathrm{ER}=n^2\ln n$). Secondly, the probability that a fixed pair $xy$ forms an edge in a graph $G_\ell$ is equal to
\[
1-\left(1-\frac{1}{\binom{n}{2}}\right)^\ell\approx
1-\exp\left(-\frac{2\ell}{n^2}\right)\;,
\]
and indeed it is easy to prove that throughout the evolution the graphs $G_\ell$ stay quasirandom, of density around 
\begin{equation}\label{eq:solutionErdosRenyi}
	1-\exp(-2t)
\end{equation}
for each $t>0$ and $\ell=tn^2$, asymptotically almost surely.

Note the number of edges in the Erd\H{o}s--R\'enyi flip process after $k$ steps is not concentrated on one number (unlike in the uniform Erd\H{o}s--R\'enyi random graph) and is not binomial either (unlike in the binomial Erd\H{o}s--R\'enyi random graph). However, both these models are good proxies since our view and results will be insensitive to changes of $o(n^2)$ edges, as usual in the theory of dense graph limits. 

\subsection{The triangle removal process}\label{ssec:introTR}
The \emph{triangle removal process} was introduced by Bollob\'as and Erd\H{o}s around 1990. We start with the complete graph $G_0$. Each $G_{\ell+1}$ is created from $G_{\ell}$ by firstly selecting a uniformly random triangle in $G_\ell$ and then deleting all its three edges. Thus, the number of edges in each step of the process decreases by three and the process cannot be continued when we arrive at a triangle-free graph. Observe this stopping time is random. A result of Bohman, Frieze, and Lubetzky~\cite{MR3350225} says that there are $n^{3/2+o(1)}$ edges left in this final graph, asymptotically almost surely. To show this, Bohman, Frieze, and Lubetzky set up complicated machinery which allows them to control quasirandomness of the sequence during the evolution. In particular, throughout the evolution, the graphs stay $o(1)$-quasirandom asymptotically almost surely.

As above, we can modify the triangle removal process into the \emph{triangle removal flip process}. In this modification, in each step we sample a triple of vertices uniformly at random and if that triple induces a triangle, then we remove all the three edges. Note that the triangle removal flip process is exactly the triangle removal process, where between each pair of consecutive steps of the triangle removal process we insert a random number of copies of the first graph of the pair, the random number being geometrically distributed with a parameter corresponding to the density of triangles of this first graph. It can be proven that this process reaches a triangle-free graph in $\tau_\mathrm{TR}=n^{4+o(1)}$ steps asymptotically almost surely, at which point it stabilises (it is not trivial to prove this, however we will only need it for informal comments).

Let us now try to understand the evolution of the density in the triangle removal flip process. Suppose that the density of a graph $G_\ell$ is $d$. Then, because $G_\ell$ is quasirandom, the probability that a random triple of vertices induces a triangle is roughly $d^3$. In other words, if we fix $\eps>0$, then the number of triangle removals between rounds $\ell$ and $\ell+\eps n^2$ is around $\eps d^3 n^2$ and hence the density drops to around $d-6\eps d^3$. Note that one source of imprecision is the assumption that within this window, a triangle is still sampled with probability $d^3$ which is true only at the beginning of the window. This imprecision, however, is negligible if $\eps$ is small. So, if for $t \ge 0$ we denote by $d(t)$ the density of the graph $G_{\lfloor tn^2\rfloor}$ and we treat $d(t)$ as a number (actually, it is a random variable, though highly concentrated), then the heuristics given above suggest a differential equation
\begin{align}
	\begin{split}\label{eq:triangleremovaldensity}
		d(0)&=1\;,\\
		d'(t)&=-6d(t)^3\;,
	\end{split}
\end{align}
which has a unique solution on $[0,\infty)$, 
\begin{equation}\label{eq:solutiontriangleremoval}
	d(t)=(1+12t)^{-1/2}\;.
\end{equation}

\subsection{Flip processes}\label{ssec:introFlip}
We can now generalise the idea of local replacements that defined the Erd\H{o}s--R\'enyi flip process and the triangle removal flip process. 

We first define a subclass of flip processes which we call \emph{definite flip processes of order $k$}. Let $\lgr{k}$ be the family of graphs on vertex set $[k]$. A \emph{rule} is an arbitrary function $\mathcal{R}:\lgr{k}\rightarrow \lgr{k}$. We then consider a random discrete time process $(G_\ell : \ell = 0,1, \dots)$ of graphs on the vertex set $[n]$, where $n\ge k$ is arbitrary. We are free to choose the initial graph $G_0$ (deterministically or randomly). In each step, the graph $G_{\ell+1}$ is obtained by sampling $k$ random vertices $v_1,\ldots,v_k$  of $G_{\ell}$ (conditioned to be distinct) and replacing the induced graph $G_{\ell}[v_1,\ldots,v_k]$ --- which we call the \emph{drawn graph} --- by $\mathcal{R}(G_{\ell}[v_1,\ldots,v_k])$ --- which we call the \emph{replacement graph}.

\begin{figure}
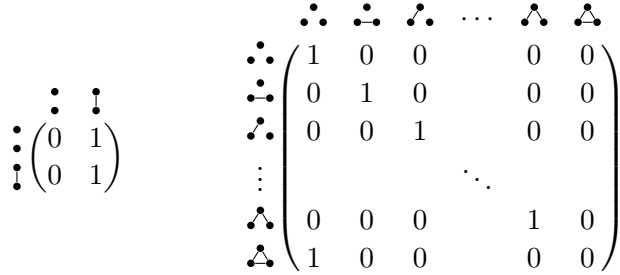

	\centering
	\[
	\begin{blockarray}{ccc}
		& \coedge & \edge  \\
		\begin{block}{c(cc)}
			\coedge   & 0 & 1   \\
			\edge     & 0 & 1   \\
		\end{block}
	\end{blockarray}\:\:\:\:\:\:\:\:\:\:\:\:\:\:\:\:\:\:
	\begin{blockarray}{ccccccc}
		& \cotrian & \cocherry & \cocherrya & \cdots & \cherry & \trian \\
		\begin{block}{c(cccccc)}
			\cotrian  & 1 & 0 &  0 &         & 0 & 0 \\
			\cocherry & 0 & 1 &  0 &        & 0 & 0 \\
			\cocherrya & 0 & 0 & 1 &         & 0 & 0 \\
			\vdots    &   &  & & \ddots    &  & \\
			\cherry   & 0 & 0 &  0 &        & 1 & 0 \\
			\trian    & 1 & 0 &  0 &        & 0 & 0 \\
		\end{block}
	\end{blockarray}
	\]
	\caption{Rules of the Erd\H{o}s--R\'enyi and triangle removal flip processes. 
	}
	\label{fig:rul_matrices}
\end{figure}

Clearly, the Erd\H{o}s--R\'enyi flip process and the triangle removal flip process are definite flip processes of order 2 and 3, respectively, see Figure~\ref{fig:rul_matrices}. Also, note that there is no need to start the Erd\H{o}s--R\'enyi flip process with the edgeless graph and the triangle removal flip process with the complete graph. For example, if we start the triangle removal flip process with a triangle-free graph, then it stabilises immediately.

To generalise these definite flip processes we allow additional randomisation in the replacement procedure. A \emph{rule of a flip process of order $k$} is a matrix $\rul \in [0,1]^{\lgr{k}\times \lgr{k}}$ such that for every $H \in \lgr{k}$ we have 
\begin{equation}\label{eq:stochastic}
	\sum_{J\in \lgr{k}}\rul_{H,J} = 1.
\end{equation}
So, given a rule $\rul$ and integer $n \ge k$, the \emph{flip process with rule $\rul$} is a discrete time process $(G_\ell : \ell = 0, \dots)$ of graphs on the vertex set $[n]$. Again, we start with a given graph $G_0$, and for $\ell \ge 0$ we obtain graph $G_{\ell+1}$ from graph $G_\ell$ as follows. Sample, uniformly at random, an ordered tuple $\mathbf{v} = (v_1,\cdots,v_k)$ of distinct vertices. Let $H \in \lgr{k}$ be the unique graph, called the \emph{drawn graph}, such that $v_iv_j$ is an edge in $G_\ell$ if and only if $ij \in E(H)$. 
Sample a graph $J$ from the \emph{replacement distribution (for the given drawn graph $H$)} $(\rul_{H,J})_J$ and replace $G_{\ell}[\mathbf{v}]$ by $J$, that is, turn $v_iv_j$ into an edge if $ij \in J$ and $v_iv_j$ into a non-edge if $ij \notin J$. We call the resulting graph $G_{\ell + 1}$. Note that the order of the elements in $\mathbf{v}$ matters.

For brevity we will sometimes say ``flip process $\rul$'' meaning ``flip process with rule $\rul$'', having in mind, precisely speaking, the class of processes (given by all possible initial graphs $G_0$). 

Flip processes are defined for all times. However, we shall only be looking at typical behaviours at times $x\cdot n^2$, $x\in [0,+\infty)$, which we call \emph{bounded horizon}. While there are certainly interesting features at later times, as we saw with $\tau_\mathrm{ER}$ and $\tau_\mathrm{TR}$ in Sections~\ref{ssec:introER} and~\ref{ssec:introTR}, these features cannot be captured by our theory.

Given a process of order $k$ with rule $\rul$, we call a graph $H\in \lgr{k}$ \emph{idle} if $\rul_{H,H}=1$. We say that the process is \emph{trivial} if each graph is idle.

Observe that each flip process is a time-homogeneous Markov process.

Let us give a sample questions we address in this paper (of course, with respect to the asymptotics $n\rightarrow\infty$; all statements hold asymptotically almost surely; the number~$99$ below is given only for the sake of concreteness).
\begin{enumerate}[label={(Q\arabic*)}]
	\item\label{q:macrostructure} Is there a flip process which would start with a quasirandom graph of order $n$ and density $0.5$ and in $99n^2$ steps produces a complete balanced bipartite graph?
	\emph{No for two reasons. Firstly, flip processes cannot create additional macroscopic structure in a bounded horizon. Secondly, they cannot achieve densities~0 or~1 (unless the corresponding part of the graph had these densities initially) in a bounded horizon.}
	\item\label{q:wipeout} Is there a flip process which would start with a complete balanced bipartite graph and in $99n^2$ steps produces a quasirandom graph of density $0.5$? 
	\emph{No, flip processes cannot entirely wipe out the macroscopic structure in a bounded horizon.}
	\item\label{q:periodic} Is there a flip process and an initial graph $G_0$ so that the evolution is $99n^2$-periodic? That is, we want that $G_{\ell}$ is similar to $G_{\ell+99n^2}$ for all $\ell$, while it is different from, say, $G_{\ell+40n^2}$. 
	\emph{Yes, we can construct an example of such a periodic process.}
	\item\label{q:backintime} Suppose our partner runs a flip process with a known rule $\rul$ but a hidden initial graph $G_0$ for $99n^2$ steps and shows us the final graph $G_{99n^2}$. Is there a way to approximately reconstruct $G_0$? 
	\emph{Yes.}	
	\item\label{q:conti} Is there a flip process where a small change in the initial graph $G_0$ would result in a large change in $G_{99n^2}$?
	\emph{No.}
\end{enumerate}

\subsection{Passing to graphons}\label{ssec:passingtographons}
In this section, we assume basic familiarity with the theory of graphons. We give more details in Section~\ref{ssec:graphons}. In particular, $\Gra$ is the set of all graphons, that is, symmetric measurable functions defined on the square of a probability space $\Omega^2$ and with values in the interval $[0,1]$. We write $\cutnd(U,W)=\cutn{U-W}$ for the cut norm distance between $U$ and $W$, and regard $\Gra$ as a metric space with metric $\cutnd$.

As we noted earlier, the Erd\H{o}s--R\'enyi flip process, when started from the edgeless graph on a vertex set $[n]$, is typically quasirandom, with edge density around $1-e^{-2t}$ in step $tn^2$. We can thus say that after rescaling the time by $n^2$, the Erd\H{o}s--R\'enyi flip process typically evolves like the function $t \mapsto \traj{t}{}_\mathrm{ER}$, where $\traj{t}{}_\mathrm{ER}\equiv 1-e^{-2t}$ is the constant graphon with density $(1-e^{-2t})$. 

Similarly, we saw that in the triangle removal flip process, when started from the complete graph, the graphs are typically quasirandom with edge density around $(1+12t)^{-1/2}$ at step $tn^2$. We can thus say that after rescaling the time by $n^2$, the triangle removal flip process typically evolves like the function $t \mapsto \traj{t}{}_\mathrm{TR}$, where $\traj{t}{}_\mathrm{TR}\equiv U_{(1+12t)^{-1/2}}$ is the constant graphon with density $(1+12t)^{-1/2}$.

In general, given a fixed rule $\rul$ of order $k$, we shall construct \emph{trajectories}, that is, a we shall define the evolution function $\traj{}{}:\Gra\times[0,+\infty)\rightarrow\Gra$, where we treat the first coordinate as the \emph{initial graphon}, and the second coordinate as the \emph{time}. Instead of writing $\Phi(W,t)$ we put the time in the superscript as $\traj{t}{(W)}$. It will turn out that the trajectories correspond to a dynamical system generalising the differential equation~\eqref{eq:triangleremovaldensity}. However, ~\eqref{eq:triangleremovaldensity} was a differential equation with a single numerical parameter (the edge density). In the general setting, the differential equation will involve a graphon-valued trajectory. So, the key objects here are Banach-space-valued differential equations (the main differential equation is given in~\eqref{eq:PhiDeriv} below). The main statement about the existence of trajectories is Theorem~{thm:flow}.

\begin{figure}
	\centering
	\includegraphics[scale=0.7]{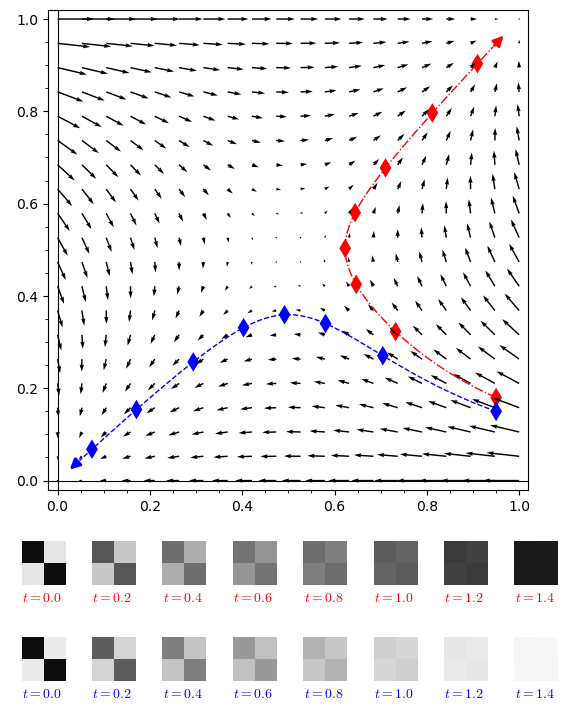}
	\caption{Extremist flip process of order three works as follows:
		if the drawn graph has two or three edges we replace it with the triangle and otherwise we replace it by the edgeless graph. Consider a class of step graphons $\{ U_{x,y} : x, y \in [0,1] \}$ with steps $\Omega_1, \Omega_2$ of measure $1/2$, where $U_{x,y}$ takes value $x$ on $\Omega_1^2 \cup \Omega^2_2$ and value $y$ elsewhere. We represent $U_{x,y}$ by a point $(x,y) \in [0,1]^2$. It turns out that trajectories starting in this class do not leave it. The arrows in the diagram depict the tangent vectors of such trajectories. 
		The red line is the trajectory from initial graphon $U_{0.95, 0.18}$ with diamond markers showing the positions along the trajectory at times $t=0$, $t=0.2$, $\ldots$, $t=1.4$; the graphon represented by each marker is plotted in grayscale underneath the diagram. The blue line is similar, from initial graphon $U_{0.95, 0.15}$. }
	\label{fig:extproc3_sims} 
\end{figure}

As a preview of our general results in Section~\ref{sec:trajectories}, in Figure~\ref{fig:extproc3_sims} we illustrate some trajectories that take values among non-constant graphons.

The \emph{Transference Theorem} (Theorem~\ref{thm:conc_proc}) says that if we start with a large initial graph $G_0$ of order $n$, fix a representation of its vertices on $\Omega$ and hence obtain a graphon representation $W_0$ of $G_0$, then, with high probability, each graph $G_{tn^2}$ (where $t\in[0,\infty)$) when represented as a graphon, is close to $\traj{t}{(W_0)}$ in the cut norm. Obviously, we then must have $\traj{0}{(U)}=U$ for each graphon $U$. Furthermore, the trajectories are continuous with respect to time, and with respect to the initial graphon. The last fact means that $G_{tn^2}$ is also close to $\traj{t}{(U)}$, if $W_0$ is close to $U$ in the cut norm. Let us state an informal version of this key result here.
\begin{thm*}[Informal and simplified consequence of Theorems~\ref{thm:flow}, \ref{thm:genome} and~\ref{thm:conc_proc}]
	Suppose that $\rul$ is a rule.
	Then there exist $\traj{}{}:\Gra\times[0,+\infty)\rightarrow\Gra$ with the following property. Suppose that $T>0$. Suppose that $G_0$ is an $n$-vertex graph and $U$ is a graphon such that the graphon representation $W_0$ of $G_0$ satisfies $\cutnd(W_0,U) =o(1)$. 
	
	Consider the flip process $(G_i)_{i \ge 0}$ with rule $\rul$ started from $G_0$. Write $W_i$ for the graphon representation of $G_i$. Then with probability $1-o(1)$ we have $\max_{i\in\N \cap [0,Tn^2]} \cutnd(W_i,\traj{i/n^2}{U})=o(1)$. (The asymptotics $o(1)$ is for fixed $\rul$, $T$ and as $n\to\infty$.)
\end{thm*}

\begin{figure}
	\includegraphics[scale=0.9]{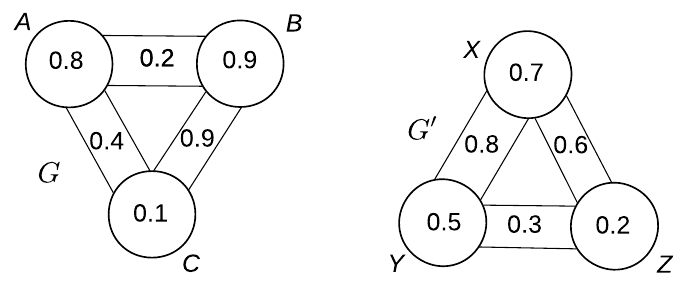}	
	\caption{Two graphs $G$ and $G'$, consisting of three blobs $A, B, C$ and $X, Y, Z$, respectively, each of order $n/3$. Inside the blobs and between the pairs of blobs we have quasirandom and quasirandom bipartite graphs, with densities as indicated. Suppose that we start a flip process with~$G$. Then our results not only say whether after a certain number of steps we typically arrive to a graph with a structure as in $G'$, but they also give, in the positive case, the correspondence of the sets $(A,B,C)$ and of $(X,Y,Z)$.}
	\label{fig:ordermatters}
\end{figure}

Note that we equipped the space of graphons with the cut norm distance $\cutnd$ rather than the cut distance $\cutm$ (which factors out isomorphic graphons). This allows us not only to track how the initial graph evolves up to isomorphism, but also tells us how individual parts of the graph evolve. We give a specific example in Figure~\ref{fig:ordermatters}.

So, Theorem~\ref{thm:conc_proc} transfers problems about random discrete trajectories on graphs to problems about deterministic continuous trajectories on graphons. With this latter view, combinatorial problems such as those we asked in Section~\ref{ssec:introFlip} translate to problems studied in dynamical systems such as the existence and stability of fixed points or periodic trajectories. 

This opens a new perspective and a new set of problems: we answer many of them, yet many more are left unexplored.

\subsection{Related work}
In addition to the Erd\H{o}s--R\'enyi and triangle removal processes there are many other graph processes considered in the literature, the $H$-free process (see~\cite{MR2657427} and references therein), the Achlio\-ptas process (see~\cite{MR2985166} and references therein), or the preferential  attachment process (see \cite{MR2091634}), to name a few. These processes do not fall into our framework, and in particular are typically studied in the sparse regime. An interesting open direction, suggested to us by Mihyun Kang, is to study modifications of flip processes which would incorporate some features of these other standard models.

For example, one could model preferential attachment by introducing `types' of vertices. The initial graph would consist of a single edge between two `active' vertices (or loop at one `active' vertex) and very many isolated `unused' vertices. The rule would then be to pick two vertices at random and if one of them is `active' and one `unused' add an edge between them and change the type of the `unused' vertex to `active'; otherwise do nothing. Note that the rate of growth would depend on the number of `unused' vertices in the construction. Further generalisations with rules of higher order are also possible.

Another related work is that of Keliger~\cite{KeligerMarkovJournal}. In that paper, vertices have a finite number of types, and an update rule at each vertex takes into account the type of that vertex and distribution of the types among the neighborhood. The paper focuses on dense graphs and uses a somewhat similar theory of ODEs.

A different line of research looks at stochastic processes on graphs which satisfy abstract properties and studies their projection into the space of graphons. This direction was started by Crane~\cite{MR3476622}, who studied exchangeable\footnote{A processes is \emph{exchangeable} if it is invariant under arbitrary relabeling of vertices by finite permutations. The concept of exchangeability is connected to graph limits as explained in~\cite{MR2426176,MR2463439}.} stochastic processes on graphs with countably many vertices. Roughly speaking, leaving out various regularity conditions, they proved that these correspond to Markov processes in the space of graphons. In the follow-up work of Athreya, den Hollander, and R\"ollin~\cite{GraphonValuedStochastic}, more complicated processes on graphs are studied, to which the corresponding graphon processes exhibit diffusive behavior.

Last, let us draw a connection between our work and \emph{Wormald's differential equation method} (DEM). This method was introduced in~\cite{wormald1995differential}, unifying many previous ad hoc analyzes, mostly randomized algorithms on random graphs. In such a setting one needs to keep track of random numerical parameters (say, the number of matched vertices in some sequential matching algorithm), and there is a heuristic suggesting that these parameters in expectation evolve according to a certain (deterministic) differential equation. DEM then gives general conditions which guarantee that with high probability the parameters indeed do follow the evolution predicted by the differential equation. In that sense there is a clear similarity to our Transference Theorem (Theorem~\ref{thm:conc_proc}). The main difference is that DEM links typical behavior of numerical parameters, whereas Theorem~\ref{thm:conc_proc} links the typical behavior of the entire structure of the evolving graph (described by the language of graphons) to a graphon-valued differential equation.

\subsection{Organisation of the paper}
We start by giving examples of several interesting classes of flip processes in Section~\ref{sec:examples}. It uses only a moderate amount of references to later parts of the paper and is therefore a good starting point for reading the paper.
Section~\ref{sec:preliminaries} gives the necessary preliminaries, mostly on graphons. In Section~\ref{sec:trajectories} we construct trajectories $\traj{}{}:\Gra\times[0,+\infty)\rightarrow\Gra$ for an arbitrary flip process as advertised in Section~\ref{ssec:passingtographons}. In fact, the construction also extends to negative times, as was already hinted by~\ref{q:backintime} in Section~\ref{ssec:introFlip}.  In Section~\ref{sec:trajectories}, we also state basic properties of trajectories, such as continuity mentioned in~\ref{q:conti}. In Section~\ref{sec:transference} we state and prove the Transference Theorem, which tells us that a typical evolution of a flip process on a finite graph stays close to the graphon trajectory. In Section~\ref{sec:properties} we establish many general features of trajectories, including those that settle questions~\ref{q:macrostructure}--\ref{q:periodic}.

\section{Examples}\label{sec:examples}
We believe that our framework provides countless interesting specific examples. Below we give seven classes of flip processes which we found interesting. See~\cite{Flip2journal} for details. In this section, $k$ is the order of the flip process in question.

\subsection{Ignorant flip processes}\label{ssec:ignorant}
Ignorant flip processes are processes where the replacement graph does not depend on the drawn graph. More  formally, we require that the rule $\rul$ of an ignorant flip process satisfies for each $H,H',J\in \lgr{k}$ that $\rul_{H,J} =\rul_{H',J}$. We then have the \emph{ignorant replacement distribution} $\mathcal{D}=\{\rul_{H,J}\}_J$ (here, $H$ is arbitrary), and its \emph{average density} $d:=\frac{1}{\binom{k}{2}}\cdot \E_{J\in\mathcal{D}}[e(J)]$.

Observe that the Erd\H{o}s--R\'enyi flip process is an ignorant flip process of order~$k=2$ with ignorant replacement distribution which is the Dirac measure on the edge $K_2$.

A crucial (and easy to establish) property is that each ignorant process of order $k\ge 2$ has all trajectories converging to the constant graphon corresponding to the average density of its replacement distribution.

\subsection{Monotone flip processes}
Roughly speaking, a monotone flip process should be one, in which the replacement graph is `monotone' with respect to the drawn graph. Besides an obvious choice whether `monotone' should mean non-decreasing or non-increasing, and whether we require non-strict monotonicity or strict monotonicity,\footnote{Of course, this strictness has to be exempted when the drawn graph is complete (for increasing flip processes) or edgeless (for decreasing flip processes).} there are three further natural definitions: (we now stick to the choice of non-decreasing flip processes) 
\begin{itemize}
	\item \emph{Pointwise non-decreasing flip processes} in which we require that $\rul_{H,J} = 0$ unless $J$ is a supergraph of $H$ (as a labelled graph). 
	\item \emph{Edge non-decreasing flip processes} in which we require that $\rul_{H,J} = 0$ unless $e(J)\ge e(H)$.
	\item \emph{Averaged-edge non-decreasing flip processes} in which we require that  $\sum_{J}\rul_{H,J}\cdot e(J) \ge e(H)$ for every $H$.
\end{itemize}
Obviously, the `averaged-edge' version is weaker than the `edge' version which is weaker than the pointwise version.  Edge non-decreasing flip processes (and thus pointwise non-decreasing flip processes) are always convergent whereas averaged-edge non-decreasing flip processes need not be.

\subsection{Removal flip processes}\label{ssec:removal}
These flip processes generalise the triangle removal flip process. That is, we have a graph $H\in \lgr{k}$. When the drawn graph $H$ is a supergraph of $F$, the replacement graph is $H\setminus F$. Other drawn graphs are idle. 

Obviously,  each removal flip process is a pointwise non-increasing flip process. The set of destinations of all trajectories is the set $\Kernel_{F}$ of $F$-free graphons, that is, the graphons $W$ with $t(F,W)=0$.

The question of speed of convergence in removal flip processes is a nontrivial one with a connection to the removal lemma.

\subsection{The complementing flip process}
\label{ssec:complementing}
The complementing flip process is a process in which the drawn graph is replaced by its complement.

Trajectories of this flip process converge very quickly to the constant-$\frac{1}{2}$ graphon.

\subsection{The component completion flip process}
The component completion flip process is a process in which the drawn graph is replaced by its \emph{component closure}, that is the replacement graph is a disjoint union of cliques which has the same component structure.

For $k\ge 3$, the trajectories of a component completion process converge to a component closure of the original graphon.

\subsection{The stirring flip process}
The stirring process is defined in two variants, \emph{firm} or \emph{loose}, as follows. Suppose that the drawn graph is $H\in \lgr{k}$. Set $M:=e(H)$ and $p := M/\bik$. The replacement graph is the uniform Erd\H{o}s--R\'enyi random graph $\G(k,M)$ (the firm variant) or the binomial Erd\H{o}s--R\'enyi random graph $\G(k,p)$ (the loose variant). 

So, in the firm stirring process we always preserve the number of edges of the underlying graph, and in the loose stirring process we preserve it in expectation.

Each trajectory converges to a constant graphon of the same density as the initial graph. For initial graphons that are regular (i.e., all the degrees are equal to a constant $d$) these trajectories are the same as for an ignorant process of replacement distribution with average density $d$. However, for non-regular initial graphons, the trajectories are different than trajectories of ignorant processes.

\subsection{The extremist process}\label{ssec:extremist}
The extremist process is defined as follows. Suppose $H\in \lgr{k}$ is the sampled graph. Set $\ell:=e(H)$. If $\ell>\binom{k}2/2$, then the replacement graph is the complete graph. If $\ell<\binom{k}2/2$, then the replacement graph is the edgeless graph. If $\ell=\binom{k}2/2$, then $\rul$ is idle.

If $W$ is a constant-$d$ graphon, then the trajectory converges to constant-$0$ for $d<\frac{1}{2}$ and to constant-$1$ for $d>\frac{1}{2}$. Trajectories starting at non-constant graphons are more subtle. For $k\ge 5$ there exists a \emph{density threshold} $d_k>0$ such that the trajectory of each graphon of edge density less than $d_k$ converges to constant-0, and symmetrically, the trajectory of each graphon of edge density more than $1-d_k$ converges to constant-1. Further, $d_k\to\frac12$ as $k\to\infty$. No positive density thresholds exist for $k=3,4$.

\section{Preliminaries}\label{sec:preliminaries}
\subsection{Basic notation and simple facts}
\subsubsection{Basic notation}
We write $(n)_k$ for the \emph{falling factorial}, $(n)_k := n\cdot(n-1)\cdot\ldots\cdot (n-k+1)$. Given a set $A$, we write $\indic_A$ for the indicator function of $A$. The symmetric difference is denoted by $\sdiff$.

\subsubsection{Graphs}
All graphs in this paper are simple, loopless and undirected. Given a graph $G = (V, E) = (V(G), E(G))$, we write $v(G) = |V|$ and $e(G) = |E|$. Given two sets $A,B\subset V$, we write $e_G(A,B)$ for the number of pairs $(a,b)\in A\times B$ such that $ab\in E$. In particular, edges in $A\cap B$ are counted twice.

Let $k\in \N$. We denote the family of all labelled graphs on the vertex set $[k]$ by $\lgr{k}$. 
Given a labeled graph $F = (V(F), E(F))$ and an \emph{ordered} subset $R \subseteq V(F)$ of vertices, we call $(R, V(F), E(F))$ a \emph{rooted graph}. Whenever $R = (a,b)$ is an ordered pair, we simply denote $(R, V(F), E(F))$ by $F^{a,b}$.

\subsection{Graphons}\label{ssec:graphons}
Our notation concerning graphons mostly follows \cite{Lovasz2012}, which is recommended for further reading.

All our graphons will be defined on an atomless measure space $(\Omega, \pi)$. The sigma-algebra is given implicitly and is assumed to be \emph{separable}, that is, there is a countable set $\mathcal C$ of measurable sets such that for every measurable set $S$ and every $\delta > 0$ there is $T \in \mathcal C$ such that $\pi (S \sdiff T) < \delta$. We will often also make the measure implicit and, in particular, for $k\in \N$ we write $\Omega^k$ for the product probability space. By $\Kernel$ we denote the space of measurable functions $W: \Omega^2 \to \R$ which are symmetric (i.e. $W(x,y) = W(y,x)$ almost everywhere\footnote{Here and below we say \emph{almost everywhere} or \emph{a.e.} when we mean \emph{$\mu$-almost everywhere}, whenever the measure $\mu$ is clear from the context; say, in this particular instance $\mu = \pi^2$.}) and bounded (i.e., $\esssup |W| < \infty$). We identify two functions that are equal a.e.\ and tacitly assume $\Kernel$ consists of equivalence classes of such functions. Hence $\Kernel$ is a closed linear subspace of the space $L^\infty(\Omega^2)$ and thus Banach when endowed with the norm $\Linf{\cdot}$. We, however, will also consider other norms on~$\Kernel$. We call elements of $\Kernel$ \emph{kernels} and elements of $\Gra := \left\{ W \in \Kernel : W \in [0,1] \right\}$ \emph{graphons}. Here and elsewhere, we write $W \in S$ to denote that up to a null set the $W$ takes values in $S$, that is, the set $\{(x,y)\in \Omega^2:W(x,y)\notin S\}$ is null. We will also write, for two measurable functions defined on the same measure space, $f = g$ and $f \le g$ meaning that equality or inequality holds a.e. With this notation, we recall that the \emph{essential supremum} $\esssup f$ is the least constant $c$ such that $f\le c$ and the \emph{essential infimum} $\essinf f$ is defined analogously.

The set $\Gra$ is a closed $\Linf{\cdot}$-ball of radius $1/2$ centered at the constant-$\frac12$ graphon. For $\eps\ge -1/2$, we define a ball with the same center but radius $\frac12+\eps$, 
\begin{equation}\label{eq:almostgraphon}
	\Gra[\eps] := \{ W \in \Kernel : W \in [-\eps,1+\eps]\}.
\end{equation}

\subsubsection{Graphon representation}
\label{sss:representation}
Suppose that $G=(V,E)$ is a graph. Suppose that we have a partition $(\Omega_v:v\in V)$ of $\Omega$ into sets of measure $1/|V|$ each. Then we define the \emph{graphon representation of $G$}, which we will denote by $W_G$, as the graphon $W:\Omega^2\rightarrow\{0,1\}$ defined to be~$1$ on each rectangle $\Omega_u\times \Omega_v$ for which $uv\in E$, and~0 otherwise. Note that this graphon representation may depend on the choice of the partition $(\Omega_v:v\in V)$. However, we will frequently keep the partition implicit. The important thing is that we will assume that graphs on the same vertex set are represented using the same partition of $\Omega$.

\subsubsection{Cut norm, cut distance and densities}
The \emph{cut norm} on $\Kernel$ is defined to be
\begin{equation}\label{eq:defcutn}
	\cutn{W} \coloneqq \sup_{A,B\subseteq\Omega}\left|\int_{A\times B}W \D\pi^2\right|,
\end{equation} 
where the supremum is taken over measurable sets.

We will occasionally use the obvious inequalities
\begin{equation}
	\label{eq:cutn_Lone_Linf}
	\cutn{W} \le \Lone{W} \le \Linf{W}, \quad W \in \Kernel.
\end{equation}

The metric induced by the cut norm is called the \emph{cut norm distance},
\[
\cutnd(U,W) \coloneqq  \cutn{U- W}\;.
\]
The \emph{cut distance} is defined for $U, W \in \Kernel$ as
\begin{equation}\label{eq:defcutdist}
	\cutm(U,W) \coloneqq \inf_{\vphi \in \mathbb{S}_{\Omega}} \cutnd\left(U, W^{\vphi}\right),
\end{equation}
where $W^{\vphi}(x,y) = W(\vphi(x), \vphi(y))$ and $\mathbb{S}_{\Omega}$ is the set of all measure preserving bijections from $\Omega$ to $\Omega$.
\begin{remark}\label{rem:factorspace}
	It is known that $\cutm(U,W) = 0$ if and only if $U$ and $W$ are \emph{weakly isomorphic}, that is, when for every finite graph $F$, we have $t(F,U) = t(F,W)$ (using the notion of densities defined in Definition~\ref{def:dens} below). Since weak isomorphism is an equivalence relation, we define $\widetilde {\Gra}$ to be $\Gra$ with equivalent graphons identified, which, together with $\cutm$ turns $\widetilde {\Gra}$ into a metric space.
\end{remark}

The following fact will also be useful.
\begin{lem}\label{lem:complete}
	The metric space $(\Gra,\cutnd)$ is complete.
\end{lem}
A proof of a substantially more general result can be found in~\cite[Proposition~1]{MR3425986}. We give a self-contained proof of Lemma~\ref{lem:complete} in the Appendix.

For a function $f : \Omega^{[k]} \to \R$ and a subset $S \subset [k]$ we write
\begin{equation*}
	\int f \D\pi^S=\int_{\Omega^S} f(x_i : i \in [k]) \prod_{i \in S} \pi(\D x_i).
\end{equation*}
We use this notation to introduce the key concept of density of a finite graph in a kernel.
\begin{defi}[Densities in kernels]\label{def:dens}
	Given a kernel $W \in \Kernel$ and a graph $F = (V,E)$, we define the \emph{density} of $F$ in $W$ as
	\[
	t(F,W) \coloneqq \int \prod_{ij \in E} W(x_i,x_j) \D\pi^{V}.
	\]
	and the \emph{induced density} as
	\[
	\tind(F,W) \coloneqq \int \prod_{ij \in E} W(x_i,x_j) \prod_{ij \notin E} \left( 1 - W(x_i,x_j) \right)\D\pi^{V}.
	\]
\end{defi}

Furthermore, we need to work with a notion of a rooted density. 

\begin{defi}[Rooted densities in kernels]
	Given a kernel $W \in \Kernel$, a rooted graph $F = (R, V, E)$, and ${\x} \in \Omega^R$ we define
	\begin{align}
		\tr{\x}(F,W) &\coloneqq \int \prod_{ij \in E} W(x_i,x_j) \D \pi^{V \sm R}\;,\mbox{and}\\
		\label{eq:tindr}
		\tindr{\x}(F,W) &\coloneqq \int \prod_{ij \in E} W(x_i,x_j) \prod_{ij \notin E} \left[ 1 - W(x_i,x_j) \right] \D \pi^{V \sm R}\;.
	\end{align}
	Note that $\x \mapsto \tr{\x}(F,W)$ and $\x \mapsto \tindr{\x}(F,W)$ are functions in $L^\infty(\Omega^R)$ and are not necessarily symmetric.
	
	For convenience we define the operator $\Troot{F}{a,b}:\Kernel \rightarrow L^\infty(\Omega^2)$, by setting 
	\begin{equation}
		\label{eq:ind_dens_op}
		(\Troot{F}{a,b} W) (x,y) = \tindr{(x,y)}(F^{a,b}, W)\;.
	\end{equation}
\end{defi}

For future reference we note that for any rooted graph $F$ and any $k\geq v(F)$,
\begin{equation}
	\label{eq:tind_t}
	\tr{\x}(F,W) = \sum_{H \in \lgr{k} : E(F) \subseteq E(H)} \tindr{\x}(H,W) \;,
\end{equation}
where in the summation $F$ is treated as a fixed rooted labeled graph, say, with root $R = [|R|]$ and vertex set $V(F) = [v(F)]$, and each $H$ is treated as the rooted graph $(R, [k], E(H))$.
Moreover, for any two graphs $F, G$ with $v(F) \le v(G)$ and for a graphon representation $W_G$ of $G$ we have (see \cite[Exercise 7.7]{Lovasz2012}) that
\begin{equation}
	\label{eq:discrete_error}
	\left| \tind(F,G) - \tind(F, W_G) \right| \le \frac{\binom{v(F)}{2}}{v(G)}\;.
\end{equation}
Similarly, if $W_G$ is a graphon representation of $G$ with respect to a partition $(\Omega_v : v \in V(G))$, then for every $u, v \in V(G)$ and $x \in \Omega_u$, $y \in \Omega_v$, we have
\begin{equation}
	\label{eq:discrete_error_rooted}
	\left| \tindr{(u,v)}(F,G) - \tindr{(x,y)}(F, W_G) \right| \le  \frac{\binom{v(F)}{2}}{v(G)}\;.
\end{equation}
(The proof is also an easy exercise.)

Next, we need a bound on expansion of the operator $\Troot{F}{a,b}$ with respect to the $L^\infty$-metric.
\begin{lem}\label{lem:contdensLinfty}
	Let $F \in \lgr{k}$ and $a,b \in [k]$ be two distinct vertices. Let $\eps \ge 0$ be a constant. If $U, W \in \Gra[\eps]$, then 
	\[
	\Linf{\Troot{F}{a,b} U - \Troot{F}{a,b} W}
	\le (1 + \eps)^{\bik - 1} \bik  \Linf{U - W}\;.
	\]
\end{lem}
We postpone the proof of the lemma to Appendix~\ref{app:graphons}.

\subsubsection{Decorated homomorphisms}
We will also need to consider homomorphism densities of graphs in $m$-tuples of graphons which can be viewed as a density of edge-coloured subgraphs, where each graphon describes a colour class. 

\begin{defi}[Densities in vectors of kernels]
	Given a graph $F = (V,E)$ and a tuple $\mathbf{W} = (W_e : e \in E)$ of kernels, we set
	\[
	t(F,\mathbf{W}) \coloneqq \int \prod_{ij \in E} W_{ij}(x_i,x_j) \D \pi^{V}\;.\]
	We can define similarly the induced density. This time, suppose that we have a graph $F$, a complete graph $K$ on $V(F)$, and a tuple $\mathbf{W} = (W_e : e \in E(K))$ of kernels. Set
	\[
	\tind(F,\mathbf{W}) \coloneqq \int \prod_{ij \in E(K)\setminus E(F)}(1-W_{ij}(x_i,x_j))\prod_{ij \in E(F)} W_{ij}(x_i,x_j) \D \pi^{V(F)}\;.\]
\end{defi}

We will make use of a counting lemma similar to \cite[Lemma~10.24]{Lovasz2012}. It is convenient to extend this to kernels. A proof can be found in Appendix~\ref{app:graphons}.

\begin{lem}\label{lem:counting_kernels}
	Given $\eps\ge 0$, a simple graph $F = (V,E)$, and two tuples $\mathbf{U} = (U_e\in\Gra[\eps] : e \in E)$ and $\mathbf{W} = (W_e\in\Gra[\eps] : e \in E)$ we have
	\[
	|t(F,\mathbf{U}) - t(F,\mathbf{W})| \le \left( 1 + 2\eps \right)^{|E| - 1}\sum_{e \in E} \cutn{U_{e} - W_{e}}.
	\]
\end{lem}

\subsubsection{Step functions and twins}
\label{ssec:twins}
We say that $f \in L^\infty(\Omega^2)$ is a \emph{step function} if there is a partition (finite or countable) $S_1 \cup S_2 \cup \dots$ of $\Omega$ into sets (\emph{steps}) of positive measure so that $f$ is constant almost everywhere on each set $S_i \times S_j$. 
We say that a partition into steps is \emph{minimal} for $f$ if $f$ is not a step function with respect to any coarser partition of $\Omega$. A \emph{step graphon} is a graphon which is a step function.

We say that $A\subset\Omega$ is a \emph{twin-set} of a kernel $W$ if there exists a conull set $A'\subset A$ such that for all $x,x'\in A'$, the (single variable) functions $W(x,\cdot)$ and $W(x',\cdot)$ are equal a.e. A variant of this concept is studied in detail in \cite[Section~13.1.1]{Lovasz2012}. Later we shall need the following easy result, a proof of which can be found in Appendix~\ref{app:graphons}.
\begin{lem}\label{lem:twinsconverge}
	Suppose that $U, U_1,U_2,\ldots$ are kernels and $\lim_{n \to \infty} \cutn{U_n - U} = 0$. If a set $A\subset \Omega$ has positive measure and is a twin-set for each $U_n$, then $A$ is a twin-set for $U$.
\end{lem}

\subsubsection{Sampling from a graphon}
\label{sss:sampling}
Each graphon $W$ defines a probability distribution on $\lgr{k}$ as follows. Sample $k$ independent $\pi$-random elements $\uu_v, v \in [k]$ of $\Omega$ and then, for each pair of distinct vertices $u, v \in [k]$, make $uv$ an edge independently with probability $W(\uu_u, \uu_v)$. For a random graph $G$ obtained this way we write $G \sim \G(k,W)$ or even use $\G(k,W)$ to denote the random graph itself.

Sometimes it will be convenient to write the rooted density defined in~\eqref{eq:tindr} in a probabilistic form as 
\begin{equation}
	\label{eq:tindr_prob}
	\tindr{\x}(F,W) = \Pc{\G(k,W) = F}{\uu_{v} = x_v : v \in R}.
\end{equation}

\subsection{Differential equations}\label{ssec:diffeqBanach}

To develop our framework we will need some theory of differential equations, both for real-valued functions and for Banach-space-valued functions. This is because our trajectories on graphons will be implicitly defined as solutions to differential equations with the derivative which is a function of the flip rule and current position/graphon.

\subsubsection{Gr\"onwall's inequality}
Among real-valued functions, the differential equation
\begin{equation}
	\label{eq:real_diff_eq}
	f'(t)=C f(t)
\end{equation}
with the initial condition $f(0)=z$ has a unique solution $f(t) = ze^{Ct}$. We will often be able to assume only a differential \emph{inequality}, which implies that a function satisfying it is bounded by the solution of a corresponding equation. A statement of this kind are known as Gr\"onwall's inequality. There are many versions of it with different assumptions on signs of constants and continuity/differentiability properties of the function in question. We shall need two versions; the first is borrowed from the literature and the second one (which we could not find a reference for) is derived using standard tools.
\begin{lem}[Integral Gr\"onwall's inequality, {\cite[Theorem 1.10]{tao2006nonlinear}}]\label{lem:Groenwall}
	Let $f : [t_0, t_1] \rightarrow [0, \infty)$
	be a continuous non-negative function.
	Suppose that for some constants $A\ge 0$, $C > 0$ we have
	\begin{equation}
		\label{eq:Viggen}
		f(t) \leq A + C \int_{t_0}^t f(x) \D x
	\end{equation}
	for every $t \in [t_0, t_1]$. Then we have 
	\begin{equation}
		f(t) \leq A e^{ C (t-t_0) } \quad \text{ for all } t \in [t_0, t_1].
	\end{equation}
\end{lem}
\begin{lem}[Differential Gr\"onwall's inequality]\label{lem:Groenwall_diff}
	Let $F : [t_0, t_1] \rightarrow \mathbb{R}$
	be an absolutely continuous function, that is, there exists an integrable function $f:[t_0, t_1]\rightarrow\R$ such that
	\begin{equation*}
		F(t) = F(t_0) + \int_{t_0}^{t} f(x) \D x
	\end{equation*}
	for all $t\in [t_0, t_1]$.
	Suppose that for some constant $C \in \R$ 
	\begin{equation}
		\label{eq:stridsberedskap}
		f(t) \leq C F(t)
	\end{equation}
	for almost every $t \in [t_0, t_1]$. Then we have 
	\begin{equation}
		\label{eq:busvaeder}
		F(t) \leq F(t_0) e^{ C (t-t_0) } \quad \text{ for all } t \in [t_0, t_1].
	\end{equation}
	Moreover if \eqref{eq:stridsberedskap} holds with $\ge$ instead of $\le$, then \eqref{eq:busvaeder} holds with $\ge$ instead of $\le$.
\end{lem}
\begin{proof}
	We revisit the proof of~\cite[Theorem 1.12]{tao2006nonlinear}.
	Define a function $G:[t_0, t_1]\rightarrow\R$ by
	\[
	G(t) = F(t)e^{-C(t-t_0)}\;.
	\]
	Since $G(t)$ is a product of two absolutely continuous functions, it is absolutely continuous (see for example~\cite[\S6.4, Exercise 42]{RoydenFitzpatrick2010}). Hence, we can use the product rule (see for example~\cite[\S6.5, Exercise 52]{RoydenFitzpatrick2010}),
	\[
	G(t) = G(t_0) + \int_{t_0}^{t} \left( f(x)e^{-C(x-t_0)} + F(x) \left(-Ce^{-C(x-t_0)}\right) \right) \D x.
	\]
	By~\eqref{eq:stridsberedskap} the integrand is non-positive almost everywhere, whence $G(t) \le G(t_0)$, which is equivalent to~\eqref{eq:busvaeder}. If we assume the opposite inequality in~\eqref{eq:stridsberedskap}, then we get the opposite inequality in~\eqref{eq:busvaeder}.
\end{proof}

\subsubsection{Banach-space-valued differential equations}
A key step in our construction involves a solution to a first-order autonomous differential equation. This solution is however not real-valued but rather function-valued. This falls within the more general setting of differential equations in Banach spaces.
Recall Definition \ref{def:deriv_normed} of the \emph{derivative}. 

The following theorem is a basic statement about local existence and uniqueness of solutions to ordinary differential equations in Banach spaces. For a reference, see, e.g.,~\cite{AbrMar88}, Lemma~4.1.6 (existence and uniqueness) and Lemma~4.1.8 
(Lipschitz continuity on the initial condition\footnote{In \cite{AbrMar88} the claim is actually weaker, but by examining the proof one easily shows \eqref{eq:deti} .
}). 
\begin{thm}
	\label{thm:diffeq}
	Consider a Banach space $\mathcal{E}$ with norm $\|\cdot\|$. Let $\mathcal{U} \subset \mathcal{E}$ be an open set, and $K, M > 0$ be constants. Let $X : \mathcal{U} \rightarrow \mathcal{E}$ be a $K$-Lipschitz function, that is  
	\[
	\norm{X(x)-X(y)} \le K \norm{x-y} \text{ for all } x,y\in \mathcal{U}.
	\]
	Suppose that
	\begin{equation}
		\label{eq:field_bound}
		\sup\{\|X(x)\| : x\in \mathcal{U} \} \le M.
	\end{equation}
	Let $x_0 \in \mathcal{U}$. If $b$ is a positive number such that the closed ball $\mathcal{B}(x_0, b) :=\{x\in \mathcal{E}: \norm{x-x_0} \le b \}$ is contained in $\mathcal{U}$, then for $\alpha := \min \left\{ 1/K, b/M \right\}$
	there is a unique continuously differentiable function $f = f_{x_0}: [-\alpha, \alpha] \to \mathcal{U}$ satisfying the equation
	\begin{equation*}
		f'(t)= X(f(t))
	\end{equation*}
	and the initial value condition $f(0) = x_0$.
	
	Furthermore, if $x_0,y_0 \in \mathcal{U}$ are such that $\mathcal{B}(x_0,b), \mathcal{B}(y_0,b)$ are contained in $\mathcal{U}$ and $f_{x_0}$ and $f_{y_0}$ correspond, respectively, to the initial values $x_0$ and $y_0$, then for any $t\in [-\alpha, \alpha]$ we have
	\begin{equation}\label{eq:deti}
		\|f_{x_0}(t)-f_{y_0}(t)\|\le e^{K|t|} \|x_0-y_0\|\;.
	\end{equation}
\end{thm}

\section{Trajectories}\label{sec:trajectories}
\subsection{Heuristics for the trajectories}
\label{ssec:heuristics}
Let us explain the main idea behind constructing trajectories $\traj{}{}:\Gra\times[0,+\infty)\rightarrow\Gra$ for a given flip process.  As we explained in Section~\ref{ssec:passingtographons}, the main motivation for trajectories is that we want to have a Transference Theorem, which says that if $G_0$ is an $n$-vertex graph which is close in cut norm to a graphon $W_0$, then performing $tn^2$ steps of the flip process we get a graph $G_{tn^2}$, which is, with high probability, close to $\traj{t}{(W_0)}$. 
The Markov property of the flip process can be informally rephrased as ``my nearest future depends only on where I am now''. Vaguely transferring this to the graphon realm, ``the nearest future'' corresponds to the time derivative $\frac{\D}{\D t}\traj{t}{(W_0)}$, and ``depends only on where I am now'' corresponds to the equation
\[
\frac{\D}{\D t}\traj{t}{(W_0)} = \vel (\traj{t}{(W_0)}),
\]
where $\vel : \Kernel \to \Kernel$ is some operator mapping kernels (and thus graphons) to kernels. In this section we will intuitively describe this (so far) hypothetical operator in the case of triangle removal flip process.

%

To interpret the time derivative, we want to see how a typical graph $G_{\delta n^2}$ evolved in $\delta n^2$ steps from an $n$-vertex graph $G_0$ which is close to $U$. Consider any points $x,y\in \Omega$. The key idea behind the cut norm is that the points $x$ and $y$ correspond to certain sets $X\subset V(G_0)$ and $Y\subset V(G_0)$, where $|X|=|Y|=\gamma n$ (and $\gamma$ should be thought of as infinitesimally small\footnote{Instead of points $x,y$ one can also approximate $U$ with a step graphon of step size $\gamma$ and take sets $I_x, I_y$ containing $x$ and $y$ of measure $\gamma$.}), and that 
\begin{equation}\label{eq:scaleme}
	U(x,y)\approx \frac{e_{G_0}(X,Y)}{\gamma^2 n^2}\;.
\end{equation}
Since in the triangle removal process edges are only removed, $U-\traj{\delta}{U}$ at $(x,y)$ corresponds to the number of edges removed from $G_0[X,Y]$ during the first $\delta n^2$ steps. How can we remove edges in one step, say the first one? Firstly, we need one of the three sampled vertices to go in $X$ and one in $Y$, which happens with probability $\approx 6\gamma^2$ (the approximation justified by the smallness of $\gamma$). Conditioning on this, we need all three vertices to form a triangle, in which case we remove \emph{one} edge from $G_0[X,Y]$. This conditional event happens with a probability approximately reflected in the 2-rooted triangle density $t^{(x,y)}(\trian^{1,2}, U)$ (cf. \eqref{eq:tindr_prob}). Hence,
\[
\E\left[e_{G_0}(X,Y)-e_{G_1}(X,Y)\right]\approx 6\gamma^2 t^{(x,y)}(\trian^{1,2}, U)\;.
\]
Since the overall changes of the host graph are fairly minor for the first $\delta n^2$ steps, the calculations are almost unchanged for $e_{G_1}(X,Y)-e_{G_2}(X,Y)$, \ldots, $e_{G_{\delta n^2}}(X,Y)-e_{G_{\delta n^2-1}}(X,Y)$. Hence, 
\[
\E [ e_{G_0}(X,Y)-e_{G_{\delta n^2}}(X,Y) ] \approx 6\gamma^2 t^{(x,y)}(\trian^{1,2}, U)\cdot \delta n^2.
\]
Recalling the scaling~\eqref{eq:scaleme}, we therefore have $(U-\traj{\delta}{U})(x,y) \approx 6\gamma^2 t^{(x,y)}(\trian^{1,2}, U)\cdot \delta$, and hence we expect that the time derivative $\lim_{\delta\rightarrow 0} [ \traj{\delta}{U}(x,y) - U(x,y) ] / \delta$ equals
\begin{equation}
	\label{eq:TR_vel_heur}
	\vel U(x,y) := -6t^{(x,y)}(\trian^{1,2}, U)\;.
\end{equation}
To generalise to more complicated flip processes, we just need to consider edge erasures and edge additions coming from different subgraph replacements of the host graph.

\subsection{The velocity operator and definition of trajectories}
In this section, we first define the velocity operator $\vel[\rul] :\Kernel \to \Kernel$ for a flip process with a rule $\rul$ of order $k$. As we saw in the motivating Section~\ref{ssec:heuristics}, 2-rooted densities of graphs in $\lgr{k}$ will be key parameters in that definition. Having defined the velocity, we then in Theorem~\ref{thm:flow} construct trajectories $\traj{}{}:\Gra\times[0,+\infty)\rightarrow\Gra$ using Theorem~\ref{thm:diffeq}. In fact it will be more convenient to define trajectories starting at arbitrary kernels, justifying the domain of the operator $\vel[\rul]$.

\begin{defi}\label{def:deriv}
	Recall the operator $\Troot{F}{a,b}:\Kernel \rightarrow L^\infty(\Omega^2)$, defined in \eqref{eq:ind_dens_op}, which, when applied to a graphon, gives the rooted induced density. Given a rule $\rul$ of order $k$, we define an operator $\vel[\rul] :\Kernel \to \Kernel$ by
	\begin{equation}
		\label{eq:velocity}
		\vel[\rul] W := \sum_{F,H\in\lgr{k}}\rul_{F,H} \sum_{1 \le a \neq b \le k} \Troot{F}{a,b} W \cdot (\ind{ab \in H \sm F} - \ind{ab \in F \sm H})\;,
	\end{equation}	
	where $ab$ is a shortening of $\left\{ a,b \right\}$. We call $\vel[\rul] W$ the \emph{velocity at $W$ for rule $\rul$}. We will often omit the rule that is clear from the context, denoting $\vel = \vel[\rul]$. 
\end{defi}

Note that $\Troot{F}{a,b} W$ is not necessarily symmetric in $(x,y)$, but $\Troot{F}{a,b} W$ and $\Troot{F}{b,a} W$ appear with the same coefficient in $\vel[\rul] W$, making it a symmetric function.

With this definition, we may formally calculate the velocity at $W$ for the rule $\rul^\mathrm{TR}$ which defines the triangle removal flip process (see Figure~\ref{fig:rul_matrices}, right). We have $k=3$, $\rul^\mathrm{TR}_{\trian, \;\cotrian }=1$ and $\rul^\mathrm{TR}_{F,F}=1$ for all $F\in \lgr{3}\backslash \{\trian\}$. This gives
\[ \vel[\rul^\mathrm{TR}] W = -\sum_{1 \leq a \neq b \leq 3} \Troot{\trian}{a,b} W = - 6 \Troot{\trian}{1,2} W\;,\]
which matches our heuristic in \eqref{eq:TR_vel_heur}. 

Calculation of the velocity will often be simplified by the following reformulation using notation from Section~\ref{sss:sampling}.
\begin{lem}\label{lem:viewvelocity}
	Suppose that a graphon $W$ and a rule $\rul$ of order $k$ are given. Let $(\ff, \hh)$ be a random pair of graphs in $\lgr{k}$ such that $\ff \sim \G(k,W)$ and $\Pc{\hh = H}{\ff = F} = \rul_{F, H}$ for every $F, H \in \lgr{k}$. We have 
	\begin{equation}
		\label{eq:vel_prob}
		(\vel[\rul] W) (x,y) = \sum_{a \neq b} \left[ \Prob( ab \in \hh | \uu_a = x, \uu_b = y ) - W(x,y) \right]\;.
	\end{equation}
\end{lem}
\begin{proof}
	By~\eqref{eq:tindr_prob},
	\begin{align*}
		\vel[\rul] W(x,y) 
		&= \sum_{a \neq b} \left[ \Pc{ ab \in \hh \sm \ff }{ \uu_a = x, \uu_b = y } - \Pc{ab \in \ff \sm \hh }{ \uu_a = x, \uu_b = y } \right] \\
		&= \sum_{a \neq b} \left[ \Pc{ ab \in \hh }{ \uu_a = x, \uu_b = y } - \Pc{ab \in \ff }{ \uu_a = x, \uu_b = y } \right] \\
		&= \sum_{a \neq b} \left[ \Pc{ ab \in \hh }{ \uu_a = x, \uu_b = y } - W(x,y) \right],
	\end{align*}
	where the second equality follows by adding and subtracting, in each term, the number \newline $\Pc{ab \in \ff \cap \hh }{ \uu_a = x, \uu_b = y}$.
\end{proof}

Let us explain how Definition~\ref{def:deriv} relates to a flip process with rule $\rul$ of order $k$. Suppose that $(\gb_i)_{i \geq 0}$ is the flip process on $n$ vertices. Fix distinct $u, v \in [n]$.
The probability that the $a$th sampled vertex is $u$ and the $b$th is $v$ is $1/(n)_2$, whence
\begin{equation*}
	\begin{split}
		&\Ec{\ind{uv \in \gb_{i+1}} - \ind{uv \in \gb_{i}}}{\gb_i}\\
		=&\tfrac{1}{(n)_2}\sum_{1 \le a \neq b \le k}\sum_{F,H\in\lgr{k}}\rul_{F,H}\cdot \tindr{(u,v)}(F^{a,b}, \gb_i) (\ind{ab \in H \sm F} - \ind{ab \in F \sm H})\;.
	\end{split}
\end{equation*}

Furthermore let $W_i$ be the graphon representation of $\gb_i$ with respect to the same partition. 
Since, by~$\eqref{eq:discrete_error_rooted}$,
\[
|\tindr{(x,y)}(F^{a,b}, \gb_i)-\tindr{(x,y)}(F^{a,b}, W_i)| \le \binom{k-2}{2}\frac{1}{n},
\]
we get the following formula linking the flip process with the velocity operator.
\begin{remark}\label{rmk:Exp_Flip_Deriv}
	Consider the flip process $(\gb_i)_{i\ge 0}$ with rule $\rul$ of order $k$. If $W_i$ is the graphon representation of $\gb_i$, then for a.e. $x, y \in \Omega$,
	\[
	\E[W_{i+1}(x,y)] = W_i(x,y) + \frac{1}{(n)_2}\vel[\rul]W_i(x,y) (1 + O(1/n))\;.
	\]
\end{remark}

We can now define the key object of our paper.
\begin{defi}
	\label{def:trajectory}
	Given a rule $\rul$ and a kernel $W \in \Kernel$, a \emph{trajectory starting at $W$} is a differentiable function $\traj{\cdot}{W} : I \to (\Kernel, \Linf{\cdot})$ defined on an open interval $I \subseteq \R$ containing zero, that satisfies the autonomous differential equation
	\begin{equation}
		\label{eq:PhiDeriv}
		\frac{\D}{\D t}\traj{t}{W} = \vel \traj{t}{W}\;
	\end{equation}
	with the initial condition
	\begin{equation}
		\label{eq:Phi_initial}
		\traj{0}{W} = W\;.\\
	\end{equation}
\end{defi}
The ultimate goal of this section is the following theorem (proved in Section \ref{ssec:traj_proof}), which asserts that trajectories exist, are unique, and whenever they start at a graphon (rather than at a general kernel), they stay among graphons for every positive time.
\begin{thm}\label{thm:flow}
	Suppose that $\rul$ is a rule.
	\begin{romenumerate}
		\item 
		\label{en:exist}
		For each kernel $W \in \Kernel$, there is an open interval $\mdom{W} \subseteq \R$ containing $0$ and a trajectory $\traj{\cdot}{W} : \mdom{W} \to \Kernel$ starting at $W$ such that any other trajectory starting at $W$ is a restriction of $\traj{\cdot}W$ to a subinterval of $\mdom{W}$. 
		
		\item
		\label{en:semi}
		For any $u \in \mdom{W}$ we have $\mdom{\traj{u}{W}} = \{ t \in \R : t + u \in \mdom{W}\}$ and for every $t \in \mdom{\traj{u}{W}}$ we have
		\begin{equation}\label{eq:semigroup_new}
			\traj{t}{ \traj{u}{W} } = \traj{t+u}{W}. 
		\end{equation}
		
		\item
		\label{en:life}
		Whenever $W \in \Gra$, the set 
		$\life{W} := \{ t \in \mdom{W} : \traj{t}W \in \Gra\}$ is a closed interval containing $[0, \infty)$.
	\end{romenumerate}
	
\end{thm}
\begin{defi}\label{def:age}
	For $W\in\Gra$ we write $\age(W) := -\inf \life{W}\in [0, \infty]$\index{$\age(W)$}.  
\end{defi}
We will sometimes (as in the following remark) abuse the notation by saying \emph{trajectory}, when we actually mean the orbit (image) of a trajectory, that is the set $\{ \traj{t}W : t \in I\}$.
\begin{remark}
	\label{rem:no_confluence}
	Note that part~\ref{en:semi} also implies that trajectories from two different kernels $W, U$ are disjoint, unless $W = \traj{t} U$ for some $t$. In other words, trajectories partition $\Kernel$.
\end{remark}

\subsection{Bounds on the velocity}

The first claim of the next lemma, \eqref{eq:vel_bound_by_W}, gives bounds on the velocity in terms of the corresponding graphon. Roughly speaking, it says that $\vel[\rul] W$ cannot be too positive for those $(x,y)$ for which $W(x,y)$ is large, and cannot be too negative for those $(x,y)$ for which $W(x,y)$ is small. Combinatorially, this is quite intuitive: if, for example, $W(x,y)$ is close to~0, then this corresponds to a graph with a sparse bipartite spot $(X,Y)$, where the vertex sets $X$ and $Y$ correspond to the points $x$ and $y$. The lower bound in \eqref{eq:vel_bound_by_W} corresponds to the fact that deletions of $(X,Y)$-edges are indeed rare, since they may only happen when we actually sample an $(X,Y)$-edge.

The following lemma also claims two more facts about the velocity operator on the Banach space $(\Kernel, \Linf{\cdot})$ that are crucial in application of Theorem~\ref{thm:diffeq}: \eqref{eq:vel_Lipsch_Linf} says it is locally Lipschitz, while \eqref{eq:vel_norm} gives a bound on the norm on a ball around the constant-$1/2$ graphon. 

\begin{lem}
	Suppose that $\rul$ is a rule of order $k$. Let $\eps \ge 0$ be a constant. Define
	\begin{equation}\label{eq:def_const}
		C_{k} := (k)_2^2 2^{\binom{k}{2} - 1}, \quad C_\infty(k,\eps) := C_k (1 + \eps)^{\bik - 1}.
	\end{equation}
	For every graphon $W \in \Gra$ we have 
	\label{lem:contderLinfty}
	\begin{equation}
		\label{eq:vel_bound_by_W}
		-(k)_2 \cdot W \le \vel[\rul] W \le (k)_2\cdot (1 - W)\;.
	\end{equation}
	For every $U, W \in \Gra[\eps]$ we have
	\begin{equation}
		\label{eq:vel_Lipsch_Linf}
		\Linf{\vel[\rul] U - \vel[\rul] W }
		\le C_\infty(k,\eps) \Linf{U - W}.
	\end{equation}
	For every $W \in \Gra[\eps]$ we have
	\begin{equation}
		\label{eq:vel_norm}
		\Linf{\vel[\rul] W} \le (k)_2 + C_\infty(k,\eps) \eps.
	\end{equation}
\end{lem}
\begin{proof} 
	Using the probabilistic form of the velocity \eqref{eq:vel_prob}, and noting that the conditional probability takes values in $[0,1]$, inequalities \eqref{eq:vel_bound_by_W} follow. 
	
	To obtain~\eqref{eq:vel_Lipsch_Linf}, using Definition \ref{def:deriv} we have
	\begin{align*}
		\Linf{ \vel[\rul] U - \vel[\rul] W  } 
		&= \Linf{ \sum_{a \neq b}\sum_{F \in \lgr{k}} \left( \Troot{F}{a,b}U - \Troot{F}{a,b}W\right) \sum_{H \in \lgr{k}}\rul_{F,H} \left(\ind{ab \in H \sm F} - \ind{ab \in F \sm H}\right) } \\ 
		&\le  \sum_{a \neq b}\sum_{F \in \lgr{k}} \Linf{ \Troot{F}{a,b}U - \Troot{F}{a,b}W } \sum_{H \in \lgr{k}}\rul_{F,H}  \\ 
		\justify{by \eqref{eq:stochastic}}&=  \sum_{a \neq b}\sum_{F \in \lgr{k}} \Linf{ \Troot{F}{a,b}U - \Troot{F}{a,b}W }  \\ 
		\justify{Lemma~\ref{lem:contdensLinfty}}&\le (k)_2 |\lgr{k}| \cdot (1 + \eps)^{\binom{k}{2} - 1} \binom{k}{2} \Linf{U - W} = C_\infty(k,\eps) \Linf{U - W}.
	\end{align*}
	
	Finally, to see~\eqref{eq:vel_norm}, given $W$, pick the nearest $U \in \Gra$ (so that $\Linf{W - U} \le \eps$). Using the triangle inequality,  \eqref{eq:vel_bound_by_W} and~\eqref{eq:vel_Lipsch_Linf}, we infer
	\[
	\Linf{\vel W} \le \Linf{\vel U} + C_\infty(k,\eps) \Linf{W - U} \le (k)_2 + C_\infty(k,\eps) \eps.
	\]
\end{proof}

The continuity property~\eqref{eq:vel_Lipsch_Linf} allows us to rewrite the differential equation~\eqref{eq:PhiDeriv} as an integral equation. This is formulated in the remark below.
\begin{remark}
	\label{rem:integral_form}
	Since $\vel$ is continuous on $(\Kernel, \Linf{\cdot})$ by~\eqref{eq:vel_Lipsch_Linf},
	Theorem~\ref{thm:fund_thm_calc} tells us that the trajectory $\traj{\cdot}W : I \to \Kernel$ satisfies
	\begin{equation*}
		\traj{t} W = W + \int_0^t \vel[\rul] \traj{\tau} W \D \tau, \quad t \in I,
	\end{equation*}
	where the integral is defined with respect to the norm $\Linf{\cdot}$.
\end{remark}

We will also need that the velocity operator is locally Lipschitz with respect to the cut norm.
\begin{lem}\label{lem:square_Lip}
	For any rule of order $k$ and every $\eps \ge 0$ and $U, W \in \Gra[\eps]$ we have 
	\begin{equation*}
		\cutn{\vel[\rul]U-\vel[\rul]W} \le C_\square \cutn{U - W},
	\end{equation*}
	where $C_\square = C_\square(k, \eps) = (k)_2^2 \left( 2 + 4\eps \right)^{\binom{k}{2}}$.
\end{lem}

\begin{proof}
	First, we want to show that for every $a,b\in [k]$, $a \ne b$, and $G\in\lgr{k}$ the following holds.
	
	\begin{claim}\label{clm:densgr}
		If $U, W \in \Gra[\eps]$, then 
		\[
		\cutn{\Troot{G}{a,b} U - \Troot{G}{a,b} W} \le (k)_2(1 + 2\eps)^{\binom{k}{2}} \cutn{U - W}.
		\]
	\end{claim}
	\begin{proof}[Proof of Claim~\ref{clm:densgr}]
		Fix an arbitrary set $A\subseteq \Omega$ and let $\indic_{A \times A}: \Omega^2 \to \{0,1\}$ be the indicator function of the set $A \times A$.
		Given $W$, $G$, $a$ and $b$, we define a tuple $\mathbf W = (W_e : e \in E(K_k))$ by
		\[
		W_e = 
		\begin{cases}
			W, &\text{ if } e \neq ab, e \in E(G) \\
			1 - W, &\text{ if } e \neq ab, e \notin E(G) \\
			W\cdot \indic_{A\times A}, &\text{ if } e = ab \in E(G) \\
			(1 - W)\cdot \indic_{A\times A}, &\text{ if } e = ab \notin E(G) \; .
		\end{cases}
		\]
		and similarly $\mathbf U = (U_e : e \in E(K_k))$, with $W$ replaced by $U$. We claim that 
		\begin{equation}\label{eq:UeWe}
			\cutn{U_{e} - W_{e}} \le \cutn{U-W}
		\end{equation}
		for every $e \in K_k$. 
		When $e \neq ab$, regardless whether $e \in E(G)$ or not, we have $U_e - W_e = U - W$ and therefore $\cutn{U_{e} - W_{e}} = \cutn{U-W}$. 
		
		When $e = ab$, we have
		$  U_e - W_e = (U - W) \indic_{A \times A}$ and therefore
		\begin{align*}
			\cutn{U_e - W_e} = \cutn{(U -W)\indic_{A \times A})} 
			&=
			\sup_{S, T \subseteq \Omega} \left|\int_{S\times T}(U - W)\cdot \indic_{A\times A} \D\pi^2\right|\\
			&=
			\sup_{S, T \subseteq \Omega}\left|\int_{(S\cap A)\times (T\cap A)}(U - W) \D\pi^2\right|
			\le \cutn{U-W}\;.
		\end{align*}
		Note that
		\begin{align*}
			\int_{A\times A} \Troot{G}{a,b}U \D \pi^2 &= \int \indic_{A \times A}(x_a,x_b)\int  
			\prod_{ij \in E(G) } U(x_i,x_j)\prod_{ij \notin E(G)} (1 - U(x_i,x_j)) \D \pi^{[k] \sm \{a,b\}} \D \pi^{\{a,b\}} \\
			&= t(K_k, \mathbf{U}).
		\end{align*}
		And similarly $\int_{A \times A} \Troot{G}{a,b}W \D \pi^2 = t(K_k, \mathbf{W})$. Hence
		\begin{align*}
			\left|\int_{A\times A} (\Troot{G}{a,b}U - \Troot{G}{a,b}W) \D\pi^2\right| 
			&= \left|t(K_k, \mathbf U) - t(K_k, \mathbf W)\right| \\
			\justify{Lemma~\ref{lem:counting_kernels}} &\le (1 + 2\eps)^{\binom{k}{2}} \sum_{e \in K_k} \cutn{U_e - W_e} \\
			\justify{by~\eqref{eq:UeWe}}&\le (1 + 2\eps)^{\binom{k}{2}}\binom{k}{2} \cutn{U - W}\;.
		\end{align*}
		Since this was verified for each set $A$, the result follows from Fact~\ref{fact:cutnorm1}.
	\end{proof}
	
	Now, using the triangle inequality and Claim~\ref{clm:densgr}, we get	
	\begin{align*}
		\cutn{\vel[\rul]U-\vel[\rul]W} 
		&\leq \sum_{1 \le a \ne b \le k} \sum_{(G,F) \in \lgr{k}\times\lgr{k}} \rul_{G,F} \cdot \cutn{\Troot{G}{a,b}U-\Troot{G}{a,b}W} \\
		\justify{by \eqref{eq:stochastic}}  &\le \sum_{1 \le a \ne b \le k} \sum_{G \in \lgr{k}} \cutn{\Troot{G}{a,b}U-\Troot{G}{a,b}W}\\
		&\leq (k)_22^{\binom{k}{2}} \cdot (k)_2(1 + 2\eps)^{\binom{k}{2}} \cutn{U - W} < C_{\square}\cutn{U - W}\;,
	\end{align*}
	which proves the lemma.
\end{proof}

\subsection{Proof of Theorem~\ref{thm:flow}}
\label{ssec:traj_proof}
Fix a rule $\rul$ of order $k$ and further write $\vel = \vel[\rul]$.

We start with existence of a time-local solution for every $W \in \Kernel$.
\begin{claim}
	\label{clm:small_int}
	For every $\rho \ge 0$ there is $T_\rho > 0$, so that for every $W \in \Gra[\rho]$ there is a unique trajectory $\traj{\cdot}W : J_\rho = (-T_\rho, T_\rho) \to \Gra[1 + \rho]$. Moreover, for every $t \in J_\rho$ the function $W \mapsto \traj{t}W$ is continuous on $(\Gra[\rho], \Linf{\cdot})$.
\end{claim}
\begin{proof}
	We apply Theorem~\ref{thm:diffeq} for the Banach space $\mathcal{E} = (\Kernel, \Linf{\cdot})$ with function $ X = \vel$. Using the notation of Theorem~\ref{thm:diffeq}, we set $\mathcal{U}$ to be the interior of $\Gra[1 + \rho]$ (equivalently, an open ball around constant-$\frac12$ of radius $\frac32 + \rho$), that is,
	\[
	\mathcal{U} = \bigcup_{\beta>0}\Gra[1 + \rho-\beta]\;.
	\]
	Let function $C_\infty$ be as in \eqref{eq:def_const}. By~\eqref{eq:vel_Lipsch_Linf}, operator $\vel$ is $K$-Lipschitz on $\mathcal{U}$ with $K = C_\infty(k, 1 + \rho)$. Setting $b := 0.9$, say, note that for every $x_0 = W \in \Gra[\rho]$, we have that the closed ball $\mathcal{B}(x_0, b)$ lies in $\mathcal{U}$. By~\eqref{eq:vel_norm}, condition \eqref{eq:field_bound} is satisfied by $M := (k)_2 + (1 + \rho) C_\infty(k,1 + \rho)$. Theorem~\ref{thm:diffeq} implies the claim with $T_\rho = \alpha = \min\{b/M, 1/K \}$, with the continuity of $W \mapsto \traj{t}W$ following from~\eqref{eq:deti}. 
\end{proof}
For a given rule, consider all trajectories starting at $W$ that are defined on an open interval containing $0$ and let $\mdom{W}$ be the union of their domains.  Using uniqueness given by Claim~\ref{clm:small_int} it is easy to see (cf. \cite[Section 2.6]{Teschl}) that any two such trajectories agree on the intersection of their domains and therefore determine a unique trajectory $\traj{\cdot}W : \mdom{W} \to \Kernel$ with a maximal domain. This proves~\ref{en:exist}.

Uniqueness also implies~\ref{en:semi}, by noting that $f_1(t) = \traj{t+u}{W}$ and $f_2(t) = \traj{t}{\traj{u}{W}}$ are two solutions for the same initial condition $f(0) = \traj{u}{W}$.

It remains to prove~\ref{en:life}. Despite being ``heuristically obvious'', its analytic proof is somewhat tedious. The first step is to show this for graphons that lie in the interior of $\Gra$ with respect to $\Linf{\cdot}$. 
\begin{claim}
	\label{clm:forward_graphon1}
	For every $\eps > 0$ and every $W\in \Gra$ such that $\eps \le W \le 1 - \eps$ and every $t\in \mdom{W} \cap [0, \infty)$ we have $\traj{t}{W}\in \Gra$.
\end{claim}
\begin{proof}
	Let $T := \sup \left\{ t \in \mdom{W} : \traj{u}W \in \Gra \forall u \in [0,t]  \right\}$ and suppose for contradiction that $T < \infty$. 
	
	The continuity\footnote{The continuity with respect to $\Linf{\cdot}$ is implied by the fact that Theorem~\ref{thm:diffeq} outputs a function which is even differentiable with respect to $\Linf{\cdot}$.} of the function $t \mapsto \traj{t}W$ implies the following three facts: $T > 0$ (since $W$ belongs to the interior of $\Gra$ with respect to $\Linf{\cdot}$); $\traj{T}{W} \in \Gra$ (since $\Gra$ is a closed set); by maximality of $T$, $\traj{T}{W}$ lies on the boundary of $\Gra$, that is,
	\begin{equation}\label{eq:boundary}
		\essinf \traj{T}{W}=0 \quad\mbox{or}\quad \esssup \traj{T}{W}=1\;.
	\end{equation}
	We will obtain a contradiction by showing that $\traj{T}{W}$ lies in the interior of~$\Gra$. 
	
	Define continuous kernel-valued functions $f_\cdot, g_\cdot$ with domain $[0, T]$ by $f_t := \vel\traj{t}W$ and $g_t = - (k)_2\traj{t}W$. By~\eqref{eq:vel_bound_by_W} we have that $f_t \ge g_t$ for $t \in [0,T]$ and by Remark~\ref{rem:integral_form} $\traj{t}W = W + \int_0^t f_\tau \D \tau$. In view of Proposition~\ref{prop:int_pointw}, we can assume that for every $(x,y)\in \Omega^2$,
	\begin{itemize}
		\item 
		the function $\tau \mapsto \vel \traj{\tau}W(x,y)$ is continuous on $[0,T]$;
		\item
		$\vel \traj{\tau}W(x,y) \ge - (k)_2 \traj{\tau}W(x,y)$, $\tau \in [0,T]$;
		\item and
		\[
		\traj{t} W (x,y) = W(x,y) + \int_0^t \vel \traj{\tau} W (x,y) \D \tau, \quad t \in [0, T].
		\]
	\end{itemize}
	Let $h(t) = \vel \traj{t} W(x,y)$ so that $H(t) := \traj{t}W(x,y)$ satisfies $H(t) = H(0) + \int_0^t h(x) \D x$ and $h(t) \ge -(k)_2H(t)$. By Lemma~\ref{lem:Groenwall_diff}, $H(T) \ge H(0) e^{-(k)_2 T}$. Since $W \ge \eps$, we conclude that $\traj{T}W \ge \eps e^{-(k)_2 T} =: \eps' > 0$. By an analogous argument applied to graphon $1 - W$, we get that $1 - W \ge \eps'$. This is a contradiction to~\eqref{eq:boundary}.
\end{proof}

We now extend Claim~\ref{clm:forward_graphon1} to all graphons.
\begin{claim}
	\label{clm:forward_graphon}
	For every $W\in \Gra$ and every $t\in \mdom{W} \cap [0, \infty)$ we have $\traj{t}{W}\in \Gra$.
\end{claim}
\begin{proof}
	Let $T_0$ be given by Claim~\ref{clm:small_int} for $\rho = 0$. 
	Observe that it is enough to prove $\traj{t}W \in \Gra$ only for times $t\in \mdom{W} \cap [0, T_0)$. Indeed, with such a weaker statement we can write any time $t$ as $t=t_1+\ldots +t_\ell$ where each $t_s\in (0,T_0)$ and then inductively for $s=1,\ldots,\ell$ deduce that $\traj{\sum_{i=1}^st_i}{W}\in\Gra$. The advantage of working on this smaller time interval is that Claim~\ref{clm:small_int} tells us that for each $\tau\in(-T_0,T_0)$, the map $U\mapsto \traj{\tau}{U}$ is $\Linf{\cdot}$-continuous on $\Gra[0] =\Gra$. 
	
	Let $W \in \Gra$ and $t \in \mdom{W} \cap [0, T_0)$ be arbitrary. Then for every $\eps > 0$ let $U_\eps := \eps + (1-2\eps)W$. Noting that $\eps \le U_\eps \le 1- \eps$ we conclude by Claim~\ref{clm:forward_graphon1} that $\traj{t}U_\eps \in \Gra$. Since $\Linf{W - U_\eps} \to 0$ as $\eps\to 0$, continuity of the map $U \mapsto \traj{t}U$ implies that $\traj{t}W \in \Gra$, as was needed.
\end{proof}
Claim~\ref{clm:forward_graphon} and \ref{en:semi} imply that $\life{W} \subseteq \R$ is an (open or closed) interval containing $[0, \infty)$. Since $\traj{\cdot}W : \mdom{W} \to (\Kernel, \Linf{\cdot})$ is continuous and $\Gra$ is closed, interval $\life{W}$ is closed with respect to $\mdom{W}$. Hence, recalling the notation $\age(W) = - \inf \life{W}$, we have that either $\life{W} = \mdom{W}$ or $\life{W} = [-\age(W), \infty)$. It remains to check that in the former case $\age(W) = \infty$, so that $\life{W} = \mdom{W} = \R$. Suppose for contradiction that we have $\age(W) < \infty$, so that $\mdom{W} = (-\age(W), \infty)$. Let $T_0$ be given by Claim~\ref{clm:small_int} for $\rho = 0$. Let $t = -\age(W) + T_0/4 \in \life{W}$. Since $\traj{t}W \in \Gra$, we can use Claim~\ref{clm:small_int} to extend the trajectory from the graphon $\traj{t}W$ to the interval $(-T_0,T_0)$. In particular, we have extended the trajectory of the graphon $W$ by at least $3T_0/4$ before $-\age(W)$, contradicting the maximality of $\mdom{W}$.

\subsection{Smoothness of trajectories in time}

The following lemma says that for any graphon the trajectory is Lipschitz in time. Inequality \eqref{eq:LinfPhi_eps} will be used in the proof of concentration of flip processes (Theorem~\ref{thm:conc_proc}).
\begin{lem}
	\label{lem:LinftyCont}
	For every rule $\rul$, every graphon $W \in \Gra$ and every $t,  u \in \life{W}$,
	\begin{equation}\label{eq:Phi_t_Lipschitz}
		\Linf{\traj{t}{W} - \traj{u}{W}} \le (k)_2 |t - u|,
	\end{equation}
	and for every $\delta > 0$,
	\begin{equation}\label{eq:LinfPhi_eps}
		\Linf{\frac{\traj{\delta}{W} - W}{\delta} - \vel W} \le C_k(k)_2 \delta,
	\end{equation}
	where $C_k$ is defined in \eqref{eq:def_const}.
\end{lem}
\begin{proof}
	By Remark~\ref{rem:integral_form} we have, for every $t, u \ge - \age(W)$,
	\begin{equation}
		\label{eq:integral_form_tu}
		\traj{t}{W} - \traj{u}{W} =\int_{u}^t \vel\traj{\tau}{W}\D \tau \;.
	\end{equation}
	Combining~\eqref{eq:norm_of_int} and the inequality $\Linf{\vel U} \le (k)_2$ for $U\in\Gra$ (see \eqref{eq:vel_norm}), we obtain
	\[
	\Linf{ \traj{t}{W} - \traj{u}{W} } \le |t-u| \sup_{\tau \in [u, t]} \Linf{ \vel\traj{\tau}{W}} \le (k)_2 |t-u|,
	\]
	thus proving~\eqref{eq:Phi_t_Lipschitz}.
	
	To prove~\eqref{eq:LinfPhi_eps}, setting $t = \delta, u = 0$ in~\eqref{eq:integral_form_tu} becomes
	\begin{equation}
		\label{eq:integral_form}
		\traj{\delta}{W} - W =\int_{0}^\delta \vel\traj{\tau}{W}\D \tau \;.
	\end{equation}
	By \eqref{eq:Phi_t_Lipschitz} we have, for every $\tau \ge 0$,
	\[
	\Linf{ \traj{\tau}{W} - W } = \Linf{ \traj{\tau}{W} - \traj{0}{W}} \le (k)_2\tau\;,
	\]
	and since $\traj{\tau}{W} \in \Gra$, by~\eqref{eq:vel_Lipsch_Linf},
	\begin{equation*}
		\Linf{ \vel\traj{\tau}{W} - \vel W } \le C_\infty(k,0) (k)_2 \tau = C_k (k)_2 \tau \;.
	\end{equation*}
	Hence, using \eqref{eq:integral_form}, and elementary properties of the integral (see Section~\ref{app:Banach_calc}),
	\begin{align*}
		\Linf{\traj{\delta}{W} - W - \delta \vel W} &= \Linf {  \int_0^\delta (\vel \traj{\tau}{W} - \vel W ) \D \tau } \\
		&= (\delta - 0) \sup_{\tau \in [0, \delta]} \Linf{\vel\traj{\tau}{W} - \vel W} \le C_k (k)_2 \delta^2.
	\end{align*}
	Dividing by $\delta$ we obtain \eqref{eq:LinfPhi_eps}.
\end{proof}

\subsection{Dependence on the initial conditions}
The following theorem implies that on any bounded time interval the trajectories depend continuously on the initial conditions. We have already used such a fact in the proof of Claim~\ref{clm:forward_graphon}, where it relied on the well-known machinery of Banach-space-valued differential equations and on~\eqref{eq:deti} in particular. However, there the metric was induced by the $L^\infty$-norm. For many combinatorial problems the cut norm is more relevant, and to deduce Lipschitz continuity we need to be more careful. Even though the velocity is also Lipschitz with respect to the cut norm (see Lemma~\ref{lem:square_Lip}), we are not in the setting of Banach spaces.\footnote{That is, we do not know if $(\Kernel, \cutn{\cdot})$ is complete. Note that the subset of uniformly bounded functions is complete by Lemma~\ref{lem:complete}.}

\begin{thm}
	\label{thm:genome}
	Given $k \in \N$, let $C_\square = C_\square(k, 0)$ as in Lemma~\ref{lem:square_Lip}.
	Let $\rul$ be a rule of order $k$, and $U,W$ be any graphons. If $t \in \life{U} \cap \life{W}$, then 
	\[
	\cutn{\traj{t}{U} - \traj{t}{W}} \le e^{C_{\square}|t|}\cutn{U - W}\;.
	\]
\end{thm}
\begin{proof}
	For $t \in \life{U} \cap \life{W} = [-\min(\age(U),\age(W)),+\infty)$, let $W_t := \traj{t}W$ and $U_t := \traj{t}U$.
	Let us first consider the case $t \ge 0$ and set $F(t):=\cutn{U_t - W_t}$. By Remark~\ref{rem:integral_form}, for any $s$ and $t$ such that $-\min(\age(U),\age(W)) \le s \le t < \infty$ we have
	\begin{equation}
		\label{eq:integral_rep}
		W_t - W_s = \int_s^t \vel W_x \D x \;,
	\end{equation} 
	where the integral is understood in terms of the Banach space $(\Kernel, \Linf{\cdot})$. Since $x \mapsto \vel W_x$ is also continuous with respect to the cut norm, integral can be treated with respect of the completion $\mathcal{E}$ of the normed space $(\Kernel, \cutn{\cdot})$. Note that $\Gra$ is a closed subset of $\mathcal{E}$ by Lemma~\ref{lem:complete}. Hence, using the triangle inequality, property of the integral \eqref{eq:norm_of_int}, and Lemma~\ref{lem:square_Lip}, we obtain, for $t \ge s \ge 0$,
	\begin{align*}
		&\justify{triangle inequality}&    F(t) - F(s) &\le \cutn{U_t - U_s - (W_t - W_s)} \\
		&\justify{by \eqref{eq:integral_rep}}&    &\le \cutn{ \int_s^t \left( \vel U_x - \vel W_x \right) \D x } \\
		&\justify{by \eqref{eq:norm_of_int}}& &\le \int_s^t \cutn{\vel U_x - \vel W_x} \D x  \\
		&\justify{Lemma~\ref{lem:square_Lip}}&    &\le C_\square \int_s^t \cutn{U_x - W_x} = C_\square \int_s^t F(x) \D x \;.
	\end{align*}
	In particular this implies that $F$ is $C_\square$-Lipschitz and so continuous. Moreover, it implies that $F$ satisfies \eqref{eq:Viggen} with $t_0 = 0$, $A = F(0)$ and $C = C_\square$. Hence Lemma~\ref{lem:Groenwall} implies Theorem for $t \ge 0$. 
	
	Now considering $t < 0$, we define 
	\[
	F(x) = \cutn{U_{-x} - W_{-x}}, \quad x \in [0, \min\left\{ \age(U), \age(W) \right\}]
	\]
	and similarly as above obtain that, for $0 \le s \le t \le \min \left\{ \age(U), \age(W) \right\}$
	\[
	F(t) - F(s) \le C_\square \int_{s}^{t} \cutn{U_{-x} - W_{-x}} \D x  = C_\square \int_s^t F(x) \D x \;,
	\]
	whence by Lemma~\ref{lem:Groenwall},
	\[
	F(t) \le F(0)e^{C_\square t}, \quad t \in [0, \min\left\{ \age(U), \age(W) \right\}].
	\]
\end{proof}

\section{The Transference Theorem}\label{sec:transference}
In this section we prove that for a quadratic number of steps the flip process starting with some large graph $G$, viewed as a graphon-valued process, is likely to stay close (with respect to the cut norm) to the deterministic trajectory $t \mapsto \traj{t}{W_{G}}$.

Let $W_G$ be the graphon representation of a graph $G$ with respect to some partition $(\Omega_v : v \in [n])$ of $\Omega$ (see Section~\ref{sss:representation}).

\begin{thm}\label{thm:conc_proc}
	For every $k \in \N$ there is a constant $C > 0$ so that the following holds.
	Given a rule $\rul$ of order $k$ and a graph $G$ on the vertex set $[n]$, let $(G_i)_{i \ge 0}$ be the flip process starting with $G_0 = G$. 
	
	For any $T > 0$ and $\eps > 0$ have
	\begin{equation}
		\label{eq:conc_proc}
		\max_{i \in (0, T n^2] \cap \Z} \cutn{W_{G_{i}} - \traj{i/n^2}{W_G}} < \eps\ \;
	\end{equation}
	with probability at least
	\begin{equation}\label{eq:probboundt}
		1 - \frac{CTe^{CT}}{\eps} \exp \left( (2 \ln 2) n - \frac{C \eps^3n^2}{e^{CT}} \right).
	\end{equation}
\end{thm}
\begin{remark}\label{rem:fromclean}
	A typical application of Theorem~\ref{thm:conc_proc} will also invoke Theorem~\ref{thm:genome}. That is, instead of relating $W_{G_{i}}$ to a trajectory started at $W_G$ (at time $i/n^2$), it seems more useful to relate it to a cleaned version of that trajectory started at some (simple) graphon $W$ which is close to $W_G$. For example, if $G$ is a quasirandom graph of density $\frac12$, then taking $W\equiv\frac12$ seems to be a good choice.
\end{remark}
\begin{remark}\label{rem:longtimes}
	While throughout the paper we are motivated by the idea that $T$ is an arbitrarily large constant (corresponding to `bounded horizons'), the probability bound~\eqref{eq:probboundt} is useful even for $T=O(\log n)$. We are not sure if this bound can be pushed further.
	
	This exact limit of applicability does not seem that important because in practice we will also be limited by Theorem~\ref{thm:genome} as we explained in Remark~\ref{rem:fromclean}. For example, suppose that we have a rule of order $3$ in which the edgeless graph stays edgeless and other graphs are replaced by a triangle. Let $G$ be a graph on $n$ vertices and $o(n^2)$ edges. This is a sparse graph and so $W\equiv 0$ seems like a good choice of a starting point for the trajectory. Obviously, we have $\traj{t}{W}\equiv 0$ but since $\cutn{W_G - 0} \ge 1/n^2$, we cannot conclude anything useful about $G_i$ for $i \gg  n^2\log n$ due to the Lipschitz constant in Theorem~\ref{thm:genome}.
\end{remark}

For convenience, we also include an easy corollary of Theorem~\ref{thm:conc_proc} phrased in terms of the cut distance.
\begin{coro}\label{cor:transferenceunlabvelled}
	For every $k \in \N$ the following holds with the same constant $C$ as in Theorem~\ref{thm:conc_proc}.
	Given a rule $\rul$ of order $k$ and an $n$-vertex graph $G$, let $(G_i)_{i \ge 0}$ be the flip process starting with $G_0 = G$. Let $W_0,W_1,\ldots$ be arbitrary graphon representations of the graphs $G_0,G_1,\ldots$.
	For any $T > 0$ and $\eps > 0$ have
	$\max \{\cutm(W_i, \traj{i/n^2}{W_0}) : i \in (0, T n^2] \cap \Z\}< \eps$ with probability at least~\eqref{eq:probboundt}.
\end{coro}

\subsection{Proof of Theorem~\ref{thm:conc_proc}}

The idea is to subdivide the time domain into intervals of size $\delta n^2$ for some small enough $\delta$ (depending on $k$, $\eps$ and $T$) and approximate, in each of the intervals, the expected change of the graphon $W_{G_i}$ linearly in terms of the velocity at the start of the interval. The following lemma quantifies the error of this approximation.

\begin{lem}\label{lem:expectedflip}
	Let $k, t, n \in \N$ be such that $n\geq 2$ and $\delta:= t/n^2 \ge 1/n$. 
	Furthermore, let $\rul$ be a rule of order $k$ and $G$ be a graph with vertex set $[n]$. Let $(G_i)_{i \ge 0}$ be the flip process with initial graph $G_0=G$.
	With probability at least
	\begin{equation}
		\label{eq:conc_prob}
		1 - 2^{2n + 1} \exp \left( -\frac{\delta^3n^2}{2(k)_2^2} \right)
	\end{equation}
	we have
	\begin{equation}
		\label{eq:linflip}
		\cutn{W_{G_{t}} - (W_{G}+\delta \vel[\rul] W_{G})}\leq C'''\delta^2\;,
	\end{equation}
	with a constant 
	\begin{equation}
		\label{eq:constantly_numb}
		C''' := \tbinom{k}{2}\left(2 \tbinom{k}{2} |\lgr{k}| + 6 + C_\square\right) + 1,
	\end{equation}
	where $C_\square = C_\square(k, 0)= (k)_2^2 2^{\binom{k}{2}}$, as in Lemma~\ref{lem:square_Lip}. 
\end{lem}

\begin{proof}
	In addition to $C_\square$ and $C'''$ we need to define two further constants, 
	\[
	C' := 2 \left( \binom{k}{2}^2|\lgr{k}| + (k)_2 \right) \quad\text{ and }\quad C'':=C_\square \bik.
	\]
	
	For simplicity write $W_i = W_{G_i}$.
	For $S \subseteq [n]$ write $\Omega_S = \cup \{\Omega_u :u \in S\}$.
	Since all graphons and kernels involved in the proof are constant on each of the sets $\Omega_u \times \Omega_v$, in view of Fact~\ref{fact:cutnorm2} we will bound the cut norm only on sets of the form $\Omega_S\times \Omega_T$.
	
	We plan to use the Hoeffding--Azuma inequality (see Section~\ref{ssec:MCDiarmid}) to show concentration of
	\[
	e_{G_t}(S,T) - e_{G_0}(S,T)
	\]
	for fixed $S, T \subseteq [n]$ and then apply the union bound over the $2^{2n}$ choices of $S, T$.
	For $i \in [t]$, let $\vv_i$ be the $k$-tuple of vertices chosen at step $i$ from graph $G_{i-1}$ and $\HH_i \in \lgr{k}$ be the graph by which we replace $G_{i-1}[\vv_i]$. We consider a filtration $\cf_0 \subseteq \cf_1 \subseteq \dots \subseteq \cf_t$, where $\cf_i$ is generated by $(\vv_j, \HH_j)_{1 \le j \le i}$.
	
	For each $i\in[t]$, define a random variable $X_i := e_{G_i}(S,T) - e_{G_{i-1}}(S,T)$. We have $e_{G_t}(S,T) - e_{G_0}(S,T) = \sum_{i = 1}^t X_i$.
	Note that $X_i$ is determined by $(\vv_j, \HH_j)_{j \in [i]}$, and so is $\cf_i$-measurable. Defining $Y_i := X_i - \Ec{X_i}{\cf_{i-1}}$, 
	we have that $Y_1, \dots, Y_t$ is a sequence of martingale differences.
	
	We will show that, with probability one,
	\begin{equation}
		\label{eq:exp_changes_sum}
		\left| \sum_{i=1}^t \Ec{X_i}{\cf_{i-1}} - t \int_{\Omega_S \times \Omega_T} \vel[\rul]W_{G_0} \right| \le (C''' - 1) \delta^2 n^2,
	\end{equation}
	and
	\begin{equation}
		\label{eq:martintail}
		\prob{\left| \sum_{i=1}^t Y_i \right| \ge \delta^2n^2  } \le 2\exp \left( - \frac{\delta^4n^4}{2t (k)_2^2 } \right) = 2\exp \left( -\frac{\delta^3n^2}{2(k)_2^2} \right).
	\end{equation}
	The triangle inequality and inequalities \eqref{eq:exp_changes_sum} and \eqref{eq:martintail} imply that, with probability at least $1 - 2\exp \left( -\delta^3n^2/\left( 2(k)_2^2 \right) \right)$,
	\begin{equation*}
		\left| e_{G_t}(S,T) - e_{G_0}(S,T) - t \int_{\Omega_S \times \Omega_T} \vel[\rul] W_{G_0} \right| \le  C'''\delta^2 n^2.
	\end{equation*}
	Taking the union bound over $2^{2n}$ choices of $S,T$, we obtain that  
	\[
	\cutn{W_{G_t} - W_{G_0} - \delta \vel[\rul]W_{G_0}} \le C'''\delta^2
	\]
	holds with probability \eqref{eq:conc_prob}.
	
	To prove \eqref{eq:exp_changes_sum}, we fix $i \in [t]$. Then for any two distinct vertices $u$ and $v$ we define random variables 
	\[
	Y_{uv} = Y_{vu}:= \Ec{\indic_{uv \in G_i} - \indic_{uv \in G_{i-1}}}{\cf_{i-1}} \;, 
	\]
	and write 
	\[
	\Ec{X_i}{\cf_{i-1}} = \sum_{u \in S, v \in T : u \neq v} Y_{uv} \;.
	\]
	Note that $Y_{uv}$ is a function of $G_{i-1}$. Hence we can fix (``condition on'') $G_{i-1}$ and calculate the expectation of $\indic_{uv \in G_i} - \indic_{uv \in G_{i-1}}$ with respect to the random pair $(\vv_i, \HH_i) \in [n]^k \times \lgr{k}$. 
	
	Writing $\vv_{i} = (v_1, \dots, v_k)$, note that events $E_{a,b} := \{ v_a = u, v_b = v \}$, $1 \le a < b \le k$, are mutually exclusive and equally likely, each of probability $1/(n)_2$. If we condition on $E_{a,b}$, then for any $F, H \in \lgr{k}$ the probability that the function $j \mapsto v_j, j \in [k]$ defines an isomorphism between $F$ and $G_{i-1}[\vv_i]$ is exactly $\tindr{(u,v)}(F^{a,b}, G_{i-1})$. The probability that then $G_{i-1}[\vv_i]$ is replaced by $H$ is $\rul_{F,H}$. Hence 
	\begin{equation}
		\label{eq:stika_zere_mouchy}
		Y_{uv} =\sum_{1 \le a \neq b \le k}\frac{1}{(n)_2}\sum_{F,H\in\lgr{k}} \tindr{(u,v)}(F^{a,b}, G_{i-1}) \cdot \rul_{F,H} \cdot (\ind{ab \in H \sm F} - \ind{ab \in F \sm H})\;.
	\end{equation}
	Now we want to relate \eqref{eq:stika_zere_mouchy} to the function $\vel[\rul]W_{i-1}$, which is constant on the set $\Omega_u \times \Omega_v$.
	By~\eqref{eq:discrete_error_rooted}, for $x \in \Omega_u, y \in \Omega_v$ we have 
	\[
	\left| \tindr{(u,v)}(F^{a,b}, G_{i-1}) - \tindr{(x, y)}(F^{a,b}, W_{i-1}) \right| \le \frac{\bik}{n},
	\]
	which allows us to rewrite~\eqref{eq:stika_zere_mouchy} as
	\begin{equation*}
		Y_{uv} =\sum_{1 \le a \neq b \le k}\frac{1}{(n)_2}\sum_{F,H\in\lgr{k}} \left(\tindr{(x, y)}(F^{a,b}, W_{i-1})\pm\frac{\bik}{n}\right) \cdot \rul_{F,H} \cdot (\ind{ab \in H \sm F} - \ind{ab \in F \sm H})\;.
	\end{equation*}
	To bound the aggregate contribution of the error term $\pm\frac{\binom{k}{2}}{n}$, we need to include a factor of $\bik$ because of the summation over $1 \le a \neq b \le k$, a factor of $|\lgr{k}|$ because of the summation of $H\in \lgr{k}$, and a factor of~$1$ because of the summation $F\in \lgr{k}$ (here, we make use of the term $\rul_{F,H}$ and~\eqref{eq:stochastic}). Hence, 
	\begin{equation*}
		Y_{uv} =\sum_{1 \le a \neq b \le k}\frac{1}{(n)_2}\sum_{F,H\in\lgr{k}} \tindr{(x, y)}(F^{a,b}, W_{i-1}) \cdot \rul_{F,H} \cdot (\ind{ab \in H \sm F} - \ind{ab \in F \sm H}) \;\pm\frac{\bik^2|\lgr{k}|}{n}\;,
	\end{equation*}
	and we recover~\eqref{eq:velocity}. Thus, for any $x \in \Omega_u, y \in \Omega_v$,
	\begin{align}
		\label{eq:Y_vel_diff}
		\left|Y_{uv} - \tfrac{1}{(n)_2} \vel[\rul]W_{i-1}(x,y)\right| \le \frac{\binom{k}{2}^2|\lgr{k}|}{n^2(n-1)}\;.
	\end{align}
	Since $\vel[\rul]W_{i-1}(x,y) = n^2\int_{\Omega_u \times \Omega_v} \vel[\rul]W_{i-1}$ for any $x \in \Omega_u, y \in \Omega_v$, and~\eqref{eq:vel_norm} implies 
	\[
	\left|\int_{\Omega_u \times \Omega_v} \vel[\rul]W_{i-1}\right| \le \frac{(k)_2}{n^2},
	\]
	we obtain 
	\begin{equation}
		\label{eq:Y_int_diff}
		\left|Y_{uv} - \int_{\Omega_u \times \Omega_v} \vel[\rul]W_{i-1}\right| \le \frac{ \binom{k}{2}^2|\lgr{k}| + (k)_2}{n^2(n-1)} \le \frac{C'}{n^3},
	\end{equation}
	where the extra term $(k)_2$ is due to replacement of $(n)_2$ by $n^2$ in \eqref{eq:Y_vel_diff}. Again using the bound $\Linf{\vel W_{i-1}} \le (k)_2$ from \eqref{eq:vel_norm},
	\begin{equation}
		\label{eq:vel_diag}
		\begin{split}
			&\left|\int_{\Omega_S \times \Omega_T} \vel[\rul]W_{i-1} - \sum_{u \in S, v \in T: u \neq v} \int_{\Omega_u \times \Omega_v} \vel[\rul]W_{i-1}\right| \\
			&= \left|\sum_{u \in S \cap T} \int_{\Omega_u \times \Omega_v} \vel[\rul]W_{i-1}\right| \le \frac{(k)_2|S \cap T|}{n^2} \;.
		\end{split}
	\end{equation}
	Summing \eqref{eq:Y_int_diff} over $u \neq v$ and using \eqref{eq:vel_diag}, we obtain
	\begin{equation}
		\label{eq:splouchy}
		\left| \Ec{X_i}{\cf_{i-1}} - \int_{\Omega_S \times \Omega_T} \vel[\rul]W_{i-1} \right| \le \frac{C'|S||T|}{n^3} + \frac{(k)_2|S \cap T|}{n^2} \le \frac{C' + (k)_2}{n}.
	\end{equation}
	Now Lemma~\ref{lem:square_Lip} implies that, for $i = 1, \dots, t$,
	\begin{equation}
		\label{eq:zblunky}
		\begin{split}
			&\Big| \int_{\Omega_S \times \Omega_T} \vel[\rul]W_i - \int_{\Omega_S \times \Omega_T} \vel[\rul]W_0 \Big| \le C_\square \cutn{W_i - W_0} \\
			&\le  C_\square \int_{\Omega^2} |W_i - W_0| = \frac{C_\square \left| E(G_i) \sdiff E(G_0) \right|}{n^2}  \\
			&\le C_\square\binom{k}{2}\frac{i}{n^2} \le C_\square \bik \delta = C'' \delta\;.
		\end{split}
	\end{equation}
	Recalling the assumption $\delta \ge 1/n$, from \eqref{eq:splouchy} and \eqref{eq:zblunky} we infer that the left-hand side of \eqref{eq:exp_changes_sum} is at most
	\[
	(C' + (k)_2) \delta n + C'' \delta^2n^2 \le \left( C' + (k)_2 + C'' \right) \delta^2n^2 = (C''' - 1) \delta^2 n^2.
	\]
	For every $i$ we have $X_i \in [-\bik, \bik]$ and hence $\Ec{X_i}{\mathcal{F}_{i-1}} \in [-\bik, \bik]$. Consequently $|Y_i| \le (k)_2$.
	Hence, \eqref{eq:martintail} follows by applying Lemma~\ref{lem:MBoundDiff} to both $\sum_{i=1}^t Y_i$ and $\sum_{i=1}^t (-Y_i)$.
\end{proof}

We can now prove Theorem~\ref{thm:conc_proc}.
\begin{proof}[Proof of Theorem~\ref{thm:conc_proc}]
	For simplicity we write $W_i=W_{G_i}$ and $\vel W_i=\vel[\rul] W_{G_i}$. 
	
	We first show that it is enough to fix a small number $\delta > 0$ such that $t := \delta n^2$ is an integer and to prove \eqref{eq:conc_proc} only for $i$ being multiples of $t$ that fall in the interval $(0, Tn^2]$. For any integer $i \in [Tn^2]$ let $i'$ be the largest multiple of $t$ not larger than $i$, namely $i' = t\flo{i/t}$. By the triangle inequality,
	\begin{align*}
		\left|\:\cutn{W_i - \traj{i/n^2}{W_0}} -  \cutn{W_{i'} - \traj{i'/n^2}{W_0} } \right|\le \cutn{W_i - W_{i'}} + \cutn{\traj{i'/n^2}{W_0} - \traj{i/n^2}{W_0}}.
	\end{align*}
	We have $\cutn{W_i - W_{i'}} \leByRef{eq:cutn_Lone_Linf} \Lone{W_i - W_{i'}} = 2|E(G_i) \sdiff E(G_i')|/n^2 \le 2(i - i')\bik/n^2$. Moreover, by~\eqref{eq:cutn_Lone_Linf} and~\eqref{eq:Phi_t_Lipschitz}, $\cutn{\traj{i'/n^2}{W_0} - \traj{i/n^2}{W_0}} \le (k)_2 (i - i')/n^2$. Hence
	\begin{align*}
		&\left|\:\cutn{W_i - \traj{i/n^2}{W_0}} -  \cutn{W_{i'} - \traj{i'/n^2}{W_0} } \right| \\
		\justify{triangle inequality} &\le \cutn{W_i - W_{i'}} + \cutn{\traj{i'/n^2}{W_0} - \traj{i/n^2}{W_0}} \\
		&\le \frac{2(i - i')}{n^2}\binom{k}{2} + (k)_2 \frac{i - i'}{n^2} \le   4 \bik \delta.
	\end{align*}
	Hence if we choose $\delta$ small enough to make sure that, say, $4\bik\delta \le \eps/2$, it suffices to prove that with high probability
	\begin{equation}
		\label{eq:chobotnicka}
		\max_{i \in (0, Tn^2] \cap t\Z} \cutn{W_i - \traj{i/n^2}{W_0}} \le \frac{\eps}{2}.
	\end{equation}
	With foresight we define a constant (depending on $k$)
	\begin{equation*}
		C^* := \frac{C_\square}{2 (C''' + C_k(k)_2 + 1) },
	\end{equation*}
	where $C_\square = C_\square(k,0)$ is the constant defined in Lemma~\ref{lem:square_Lip}, $C_k$ is defined in \eqref{eq:def_const}, and $C'''$ is defined in~\eqref{eq:constantly_numb}.
	Now we set
	\begin{equation}
		\label{eq:delta_upper}
		\delta = \frac{C^* \eps}{e^{C_\square T}}\;.
	\end{equation}
	We further set 
	\[
	r := |\left\{ i \in (0, Tn^2] \cap t\Z \right\}| = \flo{\frac{Tn^2}{t}} = \flo{\frac{T}{\delta}}\;.
	\]
	
	By Lemma~\ref{lem:expectedflip} and union bound 
	we have
	\begin{equation}\label{eq:short_deriv}
		\max_{i \in [r]} \cutn{  W_{it} - W_{(i-1)t} - \delta \vel W_{(i-1)t}  } \le \left( C''' + 1 \right)\delta^2\;
	\end{equation}
	with probability at least 
	\begin{align*}
		1 &- r 2^{2n+1} \exp\left(-\frac{\delta^3n^2}{2(k)_2^2} \right) \\
		&\ge 1 - \frac{Te^{C_\square T}}{C^*\eps} \exp \left( (2 \ln 2)n - \frac{C^*\eps^3n^2}{2(k)_2^2e^{3 C_\square T}} \right)\;,
	\end{align*}
	For the following fix an outcome of the flip process such that~\eqref{eq:short_deriv} holds. For every $i\in[0,r]$ we have that
	\begin{equation}\label{eq:cutnPhi_eps}
		\cutn{\Phi^{\delta}W_i - W_i - \delta \vel W_i} \leByRef{eq:cutn_Lone_Linf}
		\Linf{\Phi^{\delta}W_i - W_i - \delta \vel W_i} \leByRef{eq:LinfPhi_eps}
		C_k(k)_2 \delta^2\;.
	\end{equation}
	From \eqref{eq:short_deriv} and \eqref{eq:cutnPhi_eps} we get that for every $i \in [r]$ that
	\begin{align}
		\cutn{W_{it} - \traj{\delta}{W_{(i-1)t}} } \nonumber &\le \cutn{W_{it} - W_{(i-1)t} - \delta \vel W_{(i-1)t}} \\
		&+ \cutn{W_{(i-1)t} + \delta\vel W_{(i-1)t} - \traj{\delta}{W_{(i-1)t}}}\nonumber \\
		&\le (C''' + 1 + C_k(k)_2)  \delta^2=: \eps_2\;. \label{eq:short}
	\end{align}
	Theorem~\ref{thm:genome} implies that for every $i \in [r]$
	\begin{equation}
		\label{eq:Spargel}
		\cutn{\traj{\delta}{W_{(i-1)t}} - \traj{\delta}{\traj{(i-1)\delta}{W_0}}} \le e^{C_\square \delta} \cutn{W_{(i-1)t} - \traj{(i-1)\delta}{W_0}}\;.
	\end{equation}
	We further prove by induction that for $i = 0, 1, \dots, r$ we have
	\begin{equation}\label{eq:long}
		\cutn{W_{it}  - \traj{i\delta}{W_0}} \le \frac{\left( e^{C_\square \delta i} - 1 \right) \eps_2}{e^{C_\square \delta} - 1}\;.
	\end{equation}
	
	The base $i = 0$ is trivial and for the induction step assuming $i \ge 1$ we have
	\begin{align*} 
		\cutn{W_{it}  - \traj{i\delta}{W_0}} &\leByRef{eq:semigroup_new} \cutn{W_{it} - \traj{\delta}{W_{(i-1)t}}} + \cutn{\traj{\delta}{W_{(i-1)t}} - \traj{\delta}{\traj{(i-1)\delta}{W_0}}} \\
		\justify{by \eqref{eq:short}, \eqref{eq:Spargel}} &\le \eps_2 + e^{C \delta}\cutn{W_{(i-1)t} - \traj{(i-1)\delta}{W_0}} \\
		\justify{induction hypothesis} &\le \eps_2 +  \frac{e^{C_\square \delta}\left( e^{C_\square \delta (i-1)} - 1 \right) \eps_2}{e^{C_\square \delta} - 1} = \frac{\left( e^{C_\square \delta i} - 1 \right) \eps_2}{e^{C_\square \delta} - 1}\;,
	\end{align*}
	as was needed.
	
	This puts us in a position to prove~\eqref{eq:chobotnicka}, which was needed to finish the proof. Indeed, for every $i \in [r]$ we have
	\begin{align*}
		\cutn{W_{it} - \traj{i\delta}{W_0}} &\leByRef{eq:long}\frac{e^{C_\square \delta r} - 1}{e^{C_\square \delta} - 1} \eps_2 \\
		\justify{$r \le T/\delta$, $e^x \ge 1 + x$} &\le e^{C_\square T} \eps_2/C_\square\delta  \\
		\justify{by \eqref{eq:short}}  &= e^{C_\square T} \cdot (C''' + C_k(k)_2 + 1) \delta/C_\square  \\
		\justify{by \eqref{eq:delta_upper}}&= \eps/2\;.
	\end{align*}
\end{proof}

\subsection{Updates along a Poisson point process}\label{ssec:Poisson}
In the model of flip processes we discussed so far, the flips occur at discrete times. The following modification may be more suitable for applications. Starting with a graph $G=G^*_0$ we consider a continuous time Markov process $(G^*_t)_{t\ge 0}$ in which we update the graph along a Poisson point process with intensity~1 according to the same procedure as in the discrete case. We then have a counterpart of Theorem~\ref{thm:conc_proc}. Indeed, at any given time $Tn^2$, by the law of large numbers, the number of arrivals is $(1 + o(1))Tn^2$, with exponentially high probability.

\section{Properties of trajectories}\label{sec:properties}
In this part we collect various useful properties of trajectories with respect to a fixed flip process. Some of these properties are interesting per se, other are more of an auxiliary nature and used later on. Let us quickly summarise these results. 
In Section~\ref{ssec:backintime} we study backwards evolution of trajectories. Recall that in Definition~\ref{def:age}, we defined the age of each graphon $W$ as the maximum time by which the backward trajectory stays in the space of graphons. We define the \emph{origin of $W$} as $\traj{-\age(W)}{W}$. We also prove that $\age$ is upper semicontinuous.
For many natural flip processes and many initial graphons $W$, the trajectory $\traj{t}{W}$ converges as $t\to\infty$. In Section~\ref{ssec:destinations} we define these limits as \emph{destinations}.
Section~\ref{ssec:block} treats graphons which have a step structure, and in particular proves that this step structure is preserved throughout the evolution. In Section~\ref{ssec:val01} we prove that values~0 and~1 cannot emerge in a graphon throughout its evolution. Section~\ref{ssec:zeroonasection} proves the plausible fact that a zero velocity on a part $\Omega_1\times\Omega_2$ of $\Omega^2$ is a sufficient condition for a graphon to stay constant on $\Omega_1\times\Omega_2$ throughout its evolution. Section~\ref{ssec:constantgraphons} treats evolution of constant graphons. In particular, we get that each flip process has at least one \emph{fixed point}, that is, a graphon where the trajectory stays constant. In Section~\ref{ssec:notfixed} we prove that any nontrivial rule has at least one nontrivial trajectory. In Section~\ref{ssec:stablesinks} we introduce stable destinations as counterparts to asymptotically stable fixed points in the theory of dynamical systems.

The most interesting results about trajectories are given in Sections~\ref{ssec:oscillatory} and~\ref{ssec:complicatedtrajectories}. We try to understand the situation in which a given trajectory does not converge, thus complementing Section~\ref{ssec:destinations}.
In Section~\ref{ssec:oscillatory} we prove the existence (in a non-constructive way) of a periodic trajectory, using the Poincar\'e--Bendixson theorem as the main tool. In Section~\ref{ssec:complicatedtrajectories} we study how complicated a single trajectory can be. For example, we prove that a single trajectory cannot be dense in the space of graphons.

Section~\ref{ssec:speed} studies the speed of convergence, that is, how fast convergent trajectories approach their destination. The last Section~\ref{ssec:nonuniqueness} is short, but contains an important uniqueness question: what can we say about rules which yield the same trajectories?

\subsection{Going back in time}\label{ssec:backintime}
In this section, we look at origins of each graphon with respect to our trajectories. Recall that, given a graphon $W$, $\age(W)$ is the time for which the trajectory can be extended from $W$ backwards within the space $\Gra$ of graphons. For many flip processes, the age of many graphons is finite. For example, in the Erd\H{o}s--R\'enyi flip process, when going backwards, the values at each point of the graphon decrease faster and faster, eventually becoming negative.

If $\age(W) < \infty$, we define \emph{the origin of $W$}
\[
\index{$\source$} \source[\rul](W) := \traj{-\age(W)}W.
\]
which, by Theorem~\ref{thm:flow}\ref{en:life}, is well-defined and is indeed a graphon.
\begin{remark}
	Given a rule $\rul$, let $W\in\Gra$ be a graphon such that $\age(W)<\infty$. Writing $U:= \source(W)$, we have $\essinf U = 0$ or $\esssup U = 1$, since otherwise $U$ belongs to the interior of $\Gra$ with respect to the norm $\Linf{\cdot}$, and the trajectory starting at $U$ stays, by continuity (cf. Lemma~\ref{lem:LinftyCont}), inside $\Gra$ for some negative time contradicting the fact that $t = -\age(W)$ was the earliest time $\traj{t}W$ was a graphon.
\end{remark}

Let us now look at continuity properties of the map $W\mapsto \age(W)$ in the space of graphons. As usual, the topology of interest is that of the cut norm. Unfortunately, the map $\age$ need not be continuous. Consider the  triangle removal flip process. The age of the constant-0 graphon is~$\infty$. However, for an arbitrarily small $\alpha>0$ we can take a set $A\subset \Omega$ of measure $\alpha$ and a graphon that is zero everywhere except on a set $A\times A$ where it is~1. The age of that graphon is~$0$. To be more exact, this example shows that the map $\age$ need not be lower semicontinuous. However, upper semicontinuity is guaranteed by the following theorem.

\begin{thm}\label{thm:agesemicont}
	For any rule $\rul$, if $(W_n)_{n=1}^\infty$ is a sequence of graphons that converges to a graphon $W$ in the cut norm, then $\limsup_n \age(W_n) \le \age(W)$.
\end{thm}
\begin{proof}
	The statement is vacuous when $\limsup_n \age(W_n)=0$, so assume the contrary. We need to prove that for any (finite) $\gamma\in(0,\limsup_n \age(W_n))$ the trajectory of $W$ is well-defined in backward time until $-\gamma$, that is, there exists a graphon $U$ so that $\traj{\gamma}{U}=W$. By passing to a subsequence, we can assume that for every $n$ we have $\age(W_n) > \gamma$, and hence graphons $U_n:=\traj{-\gamma}{W_n}$ exist. We apply Theorem~\ref{thm:genome} with $t=-\gamma$ and graphons $(W_n)_n$, and see that the sequence $(U_n)_n$ is Cauchy with respect to the cut norm distance. Hence, by Lemma~\ref{lem:complete} the sequence $(U_n)_n$ converges to some graphon $U$ in the cut norm. Applying Theorem~\ref{thm:genome}, this time with $t=\gamma$ and to graphons $U, U_1, U_2, \dots$, we see that $(\traj{\gamma}{U_n}=W_n)_n$ converges to $\traj{\gamma}{U}$, implying $\traj{\gamma}U = W$.
\end{proof}

\subsubsection{Going back in time even more}
\begin{itemize}
	\item If $\age(W) < \infty$, by Theorem~\ref{thm:flow} the trajectory $\traj{t}{W}$ is defined on an \emph{open} interval $\mdom{W}$ containing $[-\age(W), \infty)$. In particular it is defined for some $t < - \age(W)$. We may have $\mdom{W} \neq (-\infty, \infty)$. For example, we saw that the trajectory of the triangle removal flip process along constant graphons given by~\eqref{eq:solutiontriangleremoval} satisfies $\lim_{t\searrow -\frac1{12}}\traj{t}{(1)}=+\infty$, so $\mdom{1} = (-\frac{1}{12}, \infty)$.
	
	\item Theorem~\ref{thm:flow}\ref{en:life} implies that if we follow the trajectory from a graphon backwards, once we leave the graphon space, we never come back.
	\item We do not have a combinatorially meaningful interpretation of the trajectories outside the space of graphons. For example the velocity of the constant-$(-1)$ kernel with respect to the triangle removal flip process is positive, contrary to the intuition that the triangle removal flip process decreases values.
\end{itemize}

\subsection{Destinations}\label{ssec:destinations}
Suppose that we have a flip process with a rule $\rul$. Given a graphon $U$, we say that \emph{the trajectory of $U$ is convergent} if there exists a \emph{destination of $U$}, denoted $\dest_\rul(U)$ such that $\traj{t}{U} \tocutn \dest_\rul(U)$ as $t \to \infty$. We say that a flip process with rule $\rul$ is \emph{convergent} if all trajectories converge.

Trivially each fixed point $W$ is a destination (of the trajectory $\traj{\cdot}W \equiv W$). We observe that a destination is necessarily a fixed point.
\begin{prop}\label{prop:destinations}
	Suppose that $\rul$ is a flip process. If $U$ is a graphon with a destination $W$, then $\vel[\rul] W=0$.
\end{prop}
\begin{proof}
	For every $t \ge 0$, by the cut norm continuity of $W \mapsto \traj{t} W$ (see Lemma~\ref{lem:square_Lip})  we have that (will all limits with respect to the cut norm)
	\[
	W = \lim_{u \to \infty} \traj{t + u}U = \lim_{u \to \infty} \traj{t}\traj{u}U = \traj{t} (\lim_{u \to \infty} \traj{u} U) = \traj{t}W
	\]
	whence $\vel W = \lim_{t \searrow 0} (\traj{t} W - W)/t \equiv 0$.
\end{proof}

The next conjecture asserts that destinations are not only cut norm limits but also $L^1$-norm limits. It could be that they are even $L^\infty$-norm limits.
\begin{conj}\label{conj:sinkconvergeL1}
	Suppose that $\rul$ is a flip process. If $U$ is a graphon with a destination $W$, then $(\traj{t}{U})_{t\rightarrow \infty}$ converges to $W$ also in $L^1$.
\end{conj}

\subsection{Permuting the ground space}\label{ssec:repermuting}
Recall that $\mathbb{S}_{\Omega}$ is the set of all measure preserving bijections from $\Omega$ to $\Omega$. Given a rule $\rul$ and $\vphi\in \mathbb{S}_{\Omega}$, it is easy to check that for every graphon $W$ we have $\vel (W^\vphi) =(\vel W)^\vphi$. Since the operator $W \mapsto W^\vphi$ is continuous on $L^\infty(\Omega^2)$, it easily follows from the definition that for every graphon $W$ the function $t \mapsto (\traj{t}W)^{\vphi}$ is a trajectory starting at $U := W^\vphi$. By uniqueness, $\life{U}=\life{W}$ and $\traj{\cdot}U = (\traj{\cdot}W)^\vphi$. Continuity also implies that $\dest{U}=(\dest{W})^\vphi$, should at least one of the destinations exist.

\subsection{Step structure}\label{ssec:block}
One of the structurally simplest classes of graphons are step graphons. 
Related random graph models called \emph{stochastic block models} were actually studied prior to the notion of graphons (recall the terminology from Section~\ref{ssec:twins}). One can check that the velocity of a step graphon $W$ (for any fixed flip process) is a kernel with the same step structure, which in turn yields that $\traj{t}(W)$ still has the same step structure (cf.~Remark~\ref{rem:integral_form}). We derive this claim from the following slightly stronger proposition.

\begin{prop}\label{prop:twinsstay} Consider an arbitrary flip process and a graphon $W$. If $A\subset\Omega$ is a twin-set of $W$, then for any $t \in \life{W}$, the set $A$ is also a twin-set for $\traj{t}W$. 
\end{prop}
\begin{proof} 
	The set $\Kernel_A := \left\{ W \in \Kernel : A \text{ is a set of twins in } W \right\}$ is clearly a linear subspace of the space of kernels $\Kernel$. Lemma~\ref{lem:twinsconverge} tells us that $\Kernel_A$ is a closed subspace with respect to the $\Linf{\cdot}$-norm. 
	Crucially, observe that whenever $U \in \Kernel_A$, then it follows from Definition~\ref{def:deriv} that $\vel U \in \Kernel_A$. 
	Applying Theorem~\ref{thm:diffeq}, along the same lines as in Theorem~\ref{thm:flow}, to the Banach space $(\Kernel_A, \Linf{\cdot})$ (rather than to the whole $\Kernel$) it follows that $\traj{t}{W} \in \Kernel_A$.
\end{proof}

The following corollary says that the step structure of a graphon stays the same along the trajectory. A special case of the corollary is that if $W$ is not a constant graphon, then $\traj{t}(W)$ is not a constant graphon for every $t\in\life{W}$. 
\begin{coro}\label{coro:steps_stay}
	If $W$ is a step graphon with a minimal step partition $S_1\cup S_2\cup \dots$, then for every $t\in\life{W}$ the graphon $\traj{t}W$ is a step graphon with a minimal step partition $S_1\cup S_2\cup \dots$ .
\end{coro}
\begin{proof}
	If follows immediately from Proposition~\ref{prop:twinsstay} that $S_1\cup S_2\cup \dots$ is a step partition for any $\traj{t}W$. It remains to argue that it stays minimal.
	
	Suppose $\traj{t}{W} =: U$ is a step graphon with respect to a coarser step partition $T_1\cup T_2\cup \dots$. By Proposition~\ref{prop:twinsstay}, applied to graphon $U$, graphon $W = \traj{-t}{U} \in \Gra$ is a step graphon with partition $T_1\cup T_2\cup \dots$, contradicting the minimality of partition $S_1\cup S_2\cup\dots$ for $W$.
\end{proof}

\subsection{Values 0 and 1}\label{ssec:val01}
Here, we give a quick proof of the fact that we cannot arrive to values~0 and 1. This result will be useful in Section~\ref{ssec:complicatedtrajectories}.
\begin{prop}\label{prop:valueszeroone}
	For any flip process, any graphon $W$, 
	and every $t > 0$ the sets
	\[
	A_t := \left\{ (x,y) \in \Omega^2 : W(x,y) > 0 , \traj{t}W(x,y) = 0 \right\}, 
	\]
	\[
	B_t := \left\{ (x,y) \in \Omega^2 : W(x,y) < 1 , \traj{t}W(x,y) = 1 \right\}
	\]
	have measure zero.
\end{prop}
\begin{proof}
	The proof is similar to the proof of Claim~\ref{clm:forward_graphon1}, but things are now simpler since we know $\traj{t}W \in \Gra$ for every $t \ge 0$. 
	By~\eqref{eq:vel_bound_by_W}, we have $\vel \traj{t} W \ge -(k)_2 \traj{t} W$, so using Remark~\ref{rem:integral_form} and Proposition~\ref{prop:int_pointw} for justification, we can fix $x, y \in \Omega$ and assume that $F(t) = \traj{t}W (x,y)$ and $f(t) = \vel \traj{t}W(x,y)$ satisfy $F(t) - F(0) = \int_0^t f(\tau) \D \tau$ and $f(t) \ge -(k)_2 F(t)$ for $t \ge 0$. Lemma~\ref{lem:Groenwall_diff} implies that
	\[
	\traj{t}W(x,y) = F(t) \ge F(0)e^{-(k)_2t} = e^{-(k)_2 t} W(x,y). 
	\]
	It follows that $\traj{t}W > 0$ wherever $W > 0$. Since, when applying Proposition~\ref{prop:int_pointw}, we modified $W$ and $\traj{t}W$ only on sets of measure zero, the proof is complete for the set $A_t$. The proof for $B_t$ is analogous.
\end{proof}

\subsection{Zero velocity on a section}\label{ssec:zeroonasection}
Consider a subspace of graphons which agree on a fixed subset $S$. Our next proposition states an intuitive and sometimes ``obvious'' thing: if the velocity of every graphon in the subspace equals zero on $S$, then any trajectory started in the subspace stays inside it.
\begin{prop}\label{prop:zeroonasection}
	Consider a flip process of order $k$. Given two sets $\Omega_1,\Omega_2\subset \Omega$, let $S:=\Omega_1\times\Omega_2\cup \Omega_2\times\Omega_1$ and let $\alpha:S\rightarrow[0,1]$ be a symmetric function. Let $\mathcal{A}=\{W\in\Gra:W_{\restriction S}=\alpha\}$. Suppose that for each $W\in \mathcal{A}$ we have $(\vel W)_{\restriction S}=0$. 
	
	If $U\in \mathcal{A}$ and $t\ge 0$, then $\traj{t}{U}\in \mathcal{A}$\;.
\end{prop}
\begin{proof}
	We denote the $L^\infty$-norm on the space of functions with domain $S$ as $\|\cdot\|_*$. That is to make a distinction from the $L^\infty$-space on $\Omega^2$ for which we use the symbol $\|\cdot\|_\infty$. Given an arbitrary graphon $W\in \Gra$, let its \emph{$\alpha$-surgery} be a graphon $W^{[\alpha]}$ which is equal to $\alpha$ on $S$ and to $W$ on $\Omega^2\setminus S$. Observe that we always have
	\begin{equation}\label{eq:oko}
		\|W-W^{[\alpha]}\|_\infty =\|W_{\restriction S}-\alpha\|_*\;,
	\end{equation}
	and also
	\begin{equation}\label{eq:okooko}
		\|W-U\|_\infty \ge\|W_{\restriction S}-U_{\restriction S}\|_*\;.
	\end{equation}
	
	Suppose that $U\in \mathcal{A}$ is given. For each $t\ge 0$ we define $U_t:=\traj{t}{U}$, and write $f(t):=\| (U_t)_{\restriction S}-\alpha\|_*$. In particular, $f(0)=0$.
	Using that $U_{\restriction S} = \alpha$ and $\left( \vel \left( U_\tau^{[\alpha]} \right) \right)_{\restriction S} \equiv 0$, we can write
	\[
	\alpha = \left( U + \int_{\tau = 0}^t \vel \left( (U_\tau)^{[\alpha]} \right) \right)_{\restriction S}\;.
	\]
	We have 
	\begin{align*}
		f(t) & =\left\|
		\left(U+\int_{\tau=0}^t \vel U_\tau\right)_{\restriction S}
		-
		\left(U+\int_{\tau=0}^t \vel \left((U_\tau)^{[\alpha]}\right)\right)_{\restriction S}
		\right\|_*
		\\
		\justify{by~\eqref{eq:okooko}}& \le \left\|
		\left(U+\int_{\tau=0}^t \vel U_\tau\right)
		-
		\left(U+\int_{\tau=0}^t \vel \left((U_\tau)^{[\alpha]}\right)\right)
		\right\|_\infty
		\\
		\justify{triangle inequality}
		&\le\int_{\tau=0}^t
		\left\|
		\vel U_\tau
		-\vel \left((U_\tau)^{[\alpha]}\right)
		\right\|_\infty
		\\
		\justify{by~\eqref{eq:vel_Lipsch_Linf}}&\le C_k \int_{\tau=0}^t
		\left\|
		U_\tau -\left((U_\tau)^{[\alpha]}\right)
		\right\|_\infty
		\\
		\justify{by~\eqref{eq:oko}}&= C_k \int_{\tau=0}^t
		\left\|
		(U_\tau)_{\restriction S} - \alpha
		\right\|_*
		=
		C_k \int_{\tau=0}^t f(\tau)\;.
	\end{align*}
	Applying Lemma~\ref{lem:Groenwall} (with $A = 0$), we obtain, for every $t \ge 0$, that $f(t)\le A \exp(C_k t) = 0$, which implies that $f(t)=0$ and therefore $\traj{t}U = U_t \in \mathcal{A}$, as was needed.
\end{proof}

\subsection{Constant graphons}\label{ssec:constantgraphons}
Suppose that $\rul$ is a flip process of order $k$. We shall study the trajectories when restricted to the space of constant graphons (recall that by Corollary~\ref{coro:steps_stay} trajectories started from constants stay constant and that no trajectory started at a non-constant reaches a constant). Moreover, the velocity of a constant graphon is also constant.
Given a constant-$d$ graphon $W$, let $\ff$ be the drawn graph and $\hh$ be the replacement graph. By integrating \eqref{eq:vel_prob} over $\Omega^2$, it follows that
\begin{equation}
	\label{eq:vel_prob_int}
	\vel d = \int_{\Omega^2} \vel d = 2 \E \left( e(\hh) - e(\ff) \right).
\end{equation}

For each $\ell\in \{0,\ldots,\bik\}$, write 
\begin{equation}\label{eq:defineDeltas}
	\Delta_\ell := \Ec{e(\hh) - e(\ff)}{e(\ff) = \ell}
	=\binom{\bik}{\ell}^{-1}\sum_{H \in \lgr{k}:e(H)=\ell}\sum_{H' \in \lgr{k}}\rul_{H,H'}\cdot (e(H')-\ell)
	\;.
\end{equation}
That is, $\Delta_\ell$ is the expected change of the number of edges in the replacement graph in one step of the flip process, if the drawn graph is a uniformly random graph with $\ell$ edges. Note that 
\begin{equation}\label{eq:extremeDeltas}
	\Delta_0\ge 0 \quad\mbox{and}\quad \Delta_{\binom{k}{2}}\le 0\;.
\end{equation}

Since $\ff$ is a binomial random graph $\G(k, d)$, for each constant $d\in [0,1]$ the velocity satisfies
\begin{equation}\label{eq:velconst}
	\vel d = 2\sum_{\ell=0}^{\bik} \binom{\bik}{\ell} d^\ell (1-d)^{\bik-\ell}\Delta_\ell\;.
\end{equation}
Viewed as a function in $d$, this is a polynomial of degree at most $\bik$. 
Firstly, let us characterise the case when this polynomial is trivial.
\begin{fact}
	$\vel d\equiv 0$ if and only if $\Delta_\ell=0$ for all $\ell\in \{0,\ldots,\bik\}$.
\end{fact}
\begin{proof}
	If $\Delta_\ell=0$ for all $\ell\in \{0,\ldots,\bik\}$, then obviously $\vel d\equiv 0$. Suppose, on the other hand, that there exists $\Delta_h\neq 0$, and let $h$ be minimum such. Then looking at~\eqref{eq:velconst}, we have $\lim_{d\to 0} \frac{\vel d}{d^{h}}=\binom{\bik}{h}\Delta_h$ which implies that $\vel d$ is indeed non-zero.
\end{proof}
So, let us proceed under the assumption that $\vel d$ is nontrivial. The roots of $\vel d$ on the interval $[0,1]$ are those constant graphons which are fixed points of $\rul$. By~\eqref{eq:extremeDeltas} and the intermediate value theorem, there is at least one such fixed point. Hence we get the following.
\begin{prop}\label{prop:atleastonefixed}
	Any flip process has at least one fixed point.
\end{prop}
A trajectory started at a value $d\in (r,r')$ between two consecutive roots $r$ and $r'$ converges either to $r$ or $r'$ depending on the signum of $\vel$ on that interval. Further, polynomial speed of convergence (see Section~\ref{ssec:speed} for general treatment of non-constant graphons) can be easily deduced. 

So, we have reduced the problem of studying trajectories on constants to questions about roots of polynomials of the form 
\begin{equation}\label{eq:polynomialform}
	f(x)=\sum_{\ell=0}^{\binom{k}{2}}c_\ell x^\ell (1-x)^{\binom{k}{2}-\ell} \;,
\end{equation}
where $c_\ell \in 2\binom{\bik}{\ell} \cdot [-\ell,\bik-\ell]$. (It is clear from the second form of $\Delta_\ell$ in \eqref{eq:defineDeltas} that one can choose a rule $\rul$ so that $\vel[\rul]d = f(d)$.) This coefficient restriction can be relaxed by multiplying by a small constant. After this the only restrictions that remain are $c_0\ge 0$ and $c_{\bik}\le 0$, which translates as
\begin{equation}\label{eq:twoconditions}
	\mbox{$f(0)\ge 0$ and $f(1)\le 0$}\;.
\end{equation}
The space of polynomials of the form~\eqref{eq:polynomialform} (with arbitrary real coefficients) is a linear space of dimension~$\bik+1$. This can be seen from the fact that each of the $\bik+1$ coefficients $c_\ell$ can be chosen, and the only constant-0 polynomial is the one with all coefficients zero. Hence this space is the space of all polynomials of degree at most $\bik$.

Recall that $f$ was the velocity. From~\eqref{eq:twoconditions} we conclude that, for a given $k$ and arbitrary points $0\le x_1\le\ldots\le x_{m} \le 1$, where $m\le \bik$, there exists a rule of order $\rul$, which, when restricted to constant graphons, has fixed points $\{x_1,\ldots,x_m\}$ and non-fixed trajectories are monotone and take values in the intervals between the roots.

\subsection{Examples of graphons not fixed for any (nontrivial) rule}\label{ssec:notfixed}
In this section, we shall find graphons $W$ for which we cannot have $\vel[\rul] W \equiv 0$ for any nontrivial rule $\rul$.

\begin{figure}[t]
	\centering
	\includegraphics[scale = 0.7]{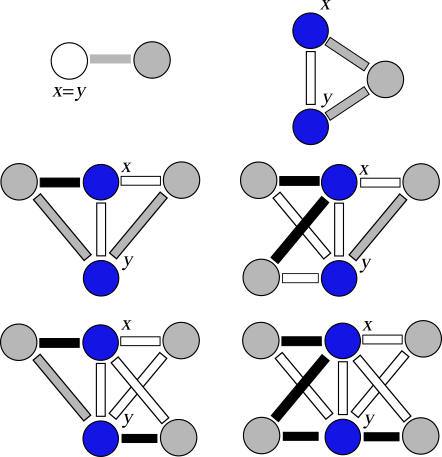}
	\caption{Examples of configurations in a step graphon preventing $\vel[\rul] W \equiv 0$ for any nontrivial rule $\rul$. White has density $0$, black has density $1$. Different grey objects allow different densities (but necessarily in $(0,1)$). Blue clusters have unrestricted densities. Changing density of $W(x,y)$ from $0$ to $1$ gives six more configurations.}
	\label{fig_no_fixed_point}
\end{figure}

We shall work with step graphons with a small number of blocks. Hence, throughout this section $\Omega$ is finite and $\pi$ is a measure such that $\pi(\{\omega\}) > 0$ for every $\omega \in \Omega$. Suppose there are steps $x, y, z \in \Omega$ (possibly $x = y$) such that $W(x, y) = 0$ and $W(x,z), W(y,z), W(z,z) \in (0,1)$, see Figure~\ref{fig_no_fixed_point}, top two examples. The value $\vel[\rul] W(x,y)$ is a linear combination of numbers $\rul_{G,H}$, where the coefficient of each $\rul_{G,H}$ is proportional to
\begin{equation*}
	\sum_{ab \in H \sm G} \tindr{(x,y)}(G^{ab}, W) - \sum_{ab \in G \sm H} \tindr{(x,y)}(G^{ab}, W)=\sum_{ab \in H \sm G} \tindr{(x,y)}(G^{ab}, W)
\end{equation*}
because of $W(x,y) = 0$. Note that the first sum is positive whenever $H \sm G \neq \emptyset$, since if we fix $ab \in H \sm G$, then by mapping $a$ to $x$, $b$ to $y$ and the remaining vertices of $G$ to $z$, we see that $\tindr{x,y}(G^{ab}, W) > 0$.
So in order to have $\vel W (x,y) = 0$ we need $\rul_{G,H} = 0$ whenever $|H \sm G| \ge 1$. In other words, we can have $\rul_{H, G} > 0$ only if $H \subseteq G$. 
Suppose that $\rul_{G,H} > 0$ for some $G,H$ such that $H \subsetneq G$. Since in $\vel W(z,z)$ the coefficient of such $\rul_{G,H}$ is, because of $W(z,z) \in (0,1)$,
\begin{equation*}
	- \sum_{ab \in G \sm H} \tindr{(z,z)}(G^{ab}, W) < 0,
\end{equation*}
implying that $\rul_{G,H} = 0$ for $H \subsetneq G$ as well. This leaves only the possibility $\rul_{G,H} > 0$ when $G = H$, hence the rule $\rul$ is trivial.

By a symmetric argument $\vel[\rul] W \equiv 0$ implies $\rul$ is trivial if value $0$ for $W(x,y)$ is replaced by value $1$. There are a few more configurations which prevent a step graphon from being a fixed point in a nontrivial rule, see Figure~\ref{fig_no_fixed_point}.

The most important corollary is the following.
\begin{coro}\label{cor:allgraphonsfixed}
	If a rule $\rul$ is such that every graphon is a fixed point, then $\rul$ is trivial.
\end{coro}

\subsection{Stable destinations}\label{ssec:stablesinks}
We say that a graphon $S\in\Gra$ is a \emph{stable destination} for a given flip process if there exists $\delta>0$ such that for each $W\in \Gra$ with $\cutnd{(W,S)}<\delta$ we have $\dest(W)=S$. We call the supremum of such $\delta$ the \emph{radius of the stable destination}.
In the dynamical systems terminology a stable destination is called an \emph{asymptotically stable fixed point}.

Many nontrivial flip processes have at least one stable destination. However, this probably is not always the case. In~\cite{Flip2journal} it is shown that the ``extremist flip process'' (see Section~\ref{ssec:extremist}) for $k=3,4$ does not have stable destinations at its most likely locations (though it does not rule out the existence of a stable destination at some other locations).

It is plausible that Theorem~\ref{thm:conc_proc} (see also Remarks~\ref{rem:fromclean} and~\ref{rem:longtimes}) can be strengthened for trajectories converging to a stable destination. That is, it is likely that a flip process started with a finite $n$-vertex graph close to such a trajectory typically stays in the neighbourhood of that destination for $\exp(\Theta(n^2))$ steps, until a really rare event kicks in.

\subsection{A periodic trajectory}\label{ssec:oscillatory}
In this section we prove existence of a flip process which has a periodic trajectory in $\Gra$ (which is not a fixed point). 

Let us first present the main idea, introducing necessary notation on the way.
Fix a partition $\Omega = \Omega_1 \cup \Omega_2$ with $\pi(\Omega_1) = \pi(\Omega_2) = 1/2$. We define a map $\Gamma:\R^2\to\Kernel$ which associates to each $(p_1,p_2) \in \R^2$ a step function 
\begin{equation}\label{eq:defGamma}
	\Gamma(p_1,p_2) = p_1 \indic_{\Omega_1^2} + p_2 \indic_{\Omega_2^2},
\end{equation}
that is, $\Gamma(p_1,p_2)$ is constant $p_i$ on $\Omega_i^2$ for $i = 1, 2$, and zero elsewhere. Define a class of graphons $\mathcal{U}_0:=\Gamma([0,1]^2)$ and a class of kernels $\mathcal{U}:=\Gamma(\R^2)$. We will show that for some $\rul$ there is a periodic trajectory taking values in $\mathcal{U}_0$. We will sometimes omit the symbol $\Gamma$ and make no distinction between points in $\R^2$ and elements of $\mathcal{U}$. In particular, we write $\vel[\rul] (p_1, p_2) := \vel[\rul] (\Gamma(p_1,p_2))$. It is immediate from the definition that $\vel[\rul] (p_1, p_2)$ is a step function with respect to the partition $\Omega = \Omega_1 \cup \Omega_2$. Furthermore, the trajectory will be confined to $\mathcal{U}_0$, which is equivalent (with the help of Proposition~\ref{prop:zeroonasection}) to
\begin{equation}\label{eq:zerooffdiag}
	\vel[\rul] (p_1, p_2) \text{ is zero on } \Omega_1\times \Omega_2.
\end{equation}
Hence, we can treat $\vel[\rul] (p_1, p_2)$ as a point in $\R^2$ via $\Gamma^{-1}$. 

For $r>0$ and $p\in\mathbb R^2$, we denote by $B(r,p)$ the open ball of radius $r$ around a point $p$, by $\overline{B(r,p)}$ the corresponding closed ball and we set $\kruh(r,p)=\overline{B(r,p)}\setminus B(r,p)$.

\begin{figure}[t]
	\centering
	\includegraphics[scale = 0.25]{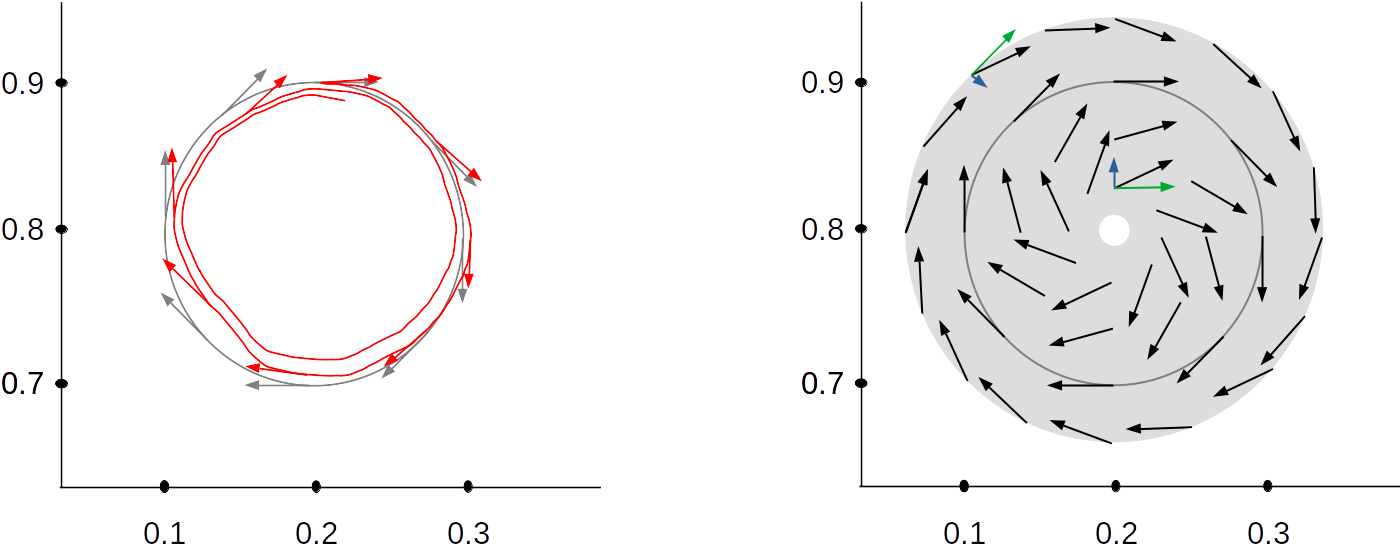}
	\caption{\emph{Left:} Ideal periodic trajectory given by the field $g$ shown in grey arrows. Red is the actual trajectory (if the rule comes only from $g$) which is perturbed by small inaccuracies and rare events. \emph{Right:} Vector field $g+h$ in black defines $\rul$ on an annulus containing $\kruh(r,q)$. Vectors of $g$ and $h$ appear in green and blue, respectively.}
	\label{fig:tangent}
\end{figure}(The following description uses tacitly the map $\Gamma$, as discussed above.) Our first attempt is to have the trajectory constrained to $\kruh(r,q)$ with $r=0.1$ and $q=(a,b):=(0.2,0.8)$. If we want to achieve that, then for each point $(x,y)\in \kruh(r,q)$, the velocity $\vel (x, y)$ has to be a constant multiple of the tangent field $g:\R^2\setminus\{q\}\rightarrow\R^2$,
\begin{equation}\label{eq:defG}
	g(x,y)=(-y+b,x-a)/\sqrt{(-y+b)^2+(x-a)^2}\;.
\end{equation}
See Figure~\ref{fig:tangent}(left).
So let us aim at $\vel[\rul] (x, y)=g(x,y)$. Assume for a while an idealised setting where we also have that
\begin{enumerate}[label={(\Alph*)}]
	\item\label{en:AA} in the drawn graph $F$, each vertex carries the information whether it was sampled from $\Omega_1$ or from $\Omega_2$, and also that it carries the information about the density $x$ on $\Omega_1\times \Omega_1$ and the density $y$ on $\Omega_2\times\Omega_2$, and
	\item\label{en:BB} the replacement distribution of the rule $\rul$ is allowed to depend on these additional input parameters.
\end{enumerate}
In such a setting, we could easily add or remove edges on the vertices sampled from $\Omega_1$ according to the first coordinate of $g(x,y)$ (even cases where, for example, the first coordinate of $g(x,y)$ is irrational could be handled by randomised outputs), and similarly for the vertices sampled from $\Omega_2$. We shall leave the pairs of vertices corresponding to $\Omega_1\times \Omega_2$ intact because of~\eqref{eq:zerooffdiag}.

In reality we do not have~\ref{en:AA}. But when the order $k$ of the flip process is large, then we can, with high probability, get a fairly precise estimate of all this data. To this end we use basic large deviation tools and further standard calculations to see that with probability $1-o(1)$ for the drawn graph $F$ for the rule, $F\sim \G(n,\Gamma(x,y))$ satisfies the following properties: 
\begin{enumerate}
	\item[($\clubsuit$)] $F$ consists of two components, each of size $(1 + o(1))\frac{k}2$ (henceforth until the end of section asymptotics are for $k \to \infty$). Further, one component forms a graph of density $d_1=x + o(1)$ and all its vertices were sampled from $\Omega_1$. The other component forms a graph of density $d_2 = y + o(1)$ and all its vertices were sampled from $\Omega_2$.
\end{enumerate}

Recall that for each $(x,y)\in \kruh(r,q)$ we have $x\le 0.3$ and $y\ge 0.7$. This means that with high probability the vertices in the sparser component were sampled from $\Omega_1$ and the vertices from the denser component from $\Omega_2$. So, to circumvent the unavailability of~\ref{en:BB} we might define the rule $\rul$ to be nontrivial on graphs of the form ($\clubsuit$) and for each such drawn graph the rule would change the density in the first and the second component according to the first and the second coordinate of $g(d_1,d_2)$, respectively.

Unfortunately, such a rule alone cannot guarantee that $\kruh(r,q)$ will be the trajectory, since we do not have any control over the $o(1)$-errors. To give a specific example, even when the current graphon is $(x,y)$, with a small but positive probability it may happen that the drawn graph $F$ will consist of two components of respective densities $d_1=\hat{x}\pm o(1)$ and $d_2=\hat{y}\pm o(1)$ for some very different $(\hat{x},\hat{y})\in \kruh(r,q)$. Hence, unwanted changes to the densities of the components of $F$ will be applied. 

To amend the situation, rather than focusing solely on the circle, we work in an annulus $C$ containing $\kruh(r,q)$  in its interior, and $\rul$ will involve all possible graphs $F$ as in $(\clubsuit)$ where $(d_1,d_2)\in C$. On $\kruh(r,q)$, the rule $\rul$ is as before. On the rest of the annulus, the rule is given by the tangent vector field $g$ composed with a vector field $h:\R^2\setminus\{q\}\rightarrow\R^2$ oriented towards $\kruh(r,q)$, 
\begin{equation}\label{eq:defH}
	h(x,y)=\left(\tfrac{r}{\sqrt{(x-a)^2+(y-b)^2}}-1\right)(x-a,y-b)\;.
\end{equation}
See Figure~\ref{fig:tangent}(right). For drawn graphs $F$ which do not consist of exactly two components whose joint densities lie in $C$, we define the replacement graph~$H$ arbitrarily, subject to the condition that
\begin{equation}\label{eq:noconnection}
	\mbox{no edges are introduced between any two components of $F$.}
\end{equation}

The point is that even though we still know that $\vel[\rul](x,y)$ is only approximately equal to $(h+g)(x,y)$ on $C$, we can prove that the undesired perturbation from rare events cannot cause that the actual velocity vector escapes from $C$. The classical Poincar\'e--Bendixson theorem then implies that there must exist a periodic trajectory. Note that~\eqref{eq:noconnection} is needed to preserve~\eqref{eq:zerooffdiag}; here even rare violations of~\eqref{eq:noconnection} would kill~\eqref{eq:zerooffdiag}.

Lemma~\ref{lem:2dimapproxrule} below allows us to produce a rule whose velocities approximate a given vector field; in fact such a rule is produced in the way we described above. To motivate some of the conditions of the lemma, we recall that we relied on the fact that the two components of the graph $F$ are distinguishable by their densities (see~\ref{cond:densdiff}), and also that in the planned change, either component cannot have density less than~0 or more than~1 (see~\ref{cond:0dens1}).

\begin{lem}\label{lem:2dimapproxrule}
	Let $\Gamma$ be defined by~\eqref{eq:defGamma}.
	For every $\eps,\delta > 0$ exists $k\in\mathbb N$ such that the following holds. If $S_1,S_2\subseteq [0,1]^2$ and $f : S_1 \to \R^2$ are such that
	\begin{romenumerate}
		\item $B(\delta,p)\subseteq S_1$ for every $p\in S_2$,\label{cond:dball}
		\item $0<\delta\leq p_1\leq p_2-\delta\leq 1-2\delta$ for every $(p_1,p_2)\in S_1$,\label{cond:densdiff}
		\item $f$ is uniformly continuous on $S_1$, and\label{cond:unif_cont}
		\item $0\leq p_i + (k)_2^{-1}f_i(p_1,p_2) \leq 1$ for every $i=1,2$ and $(p_1,p_2)\in S_1$\label{cond:0dens1},
	\end{romenumerate}
	then there exists a rule $\rul$ of order $k$ such that \eqref{eq:zerooffdiag} holds for every $p = (p_1, p_2) \in S_2$ and such that
	\begin{equation}
		\label{eq:f_appr}
		\sup_{p \in S_2} |f(p) - \Gamma^{-1}(\vel[\rul]\Gamma(p))| \le \eps.
	\end{equation}
\end{lem}

\begin{proof}
	Within the proof the asymptotic notation is with respect to $k \to \infty$ and whenever we say ``for $k$ large enough'', we mean $k$ larger than some constant depending on $\eps, \delta$.
	We fix some sequences $\alpha_1 = \alpha_1(k), \alpha_2 = \alpha_2(k)$ such that
	\[
	\tfrac{1}{\sqrt{k}} \ll  \alpha_1 \ll \alpha_2 \ll 1. 
	\]
	
	We define a rule $\rul$ as follows. 
	Assume that $F \in \lgr{k}$ is the drawn graph. We define the probabilities $\rul_{F, \cdot}$ by describing a random replacement graph $\hh$ . If $F$ has an isolated vertex or the number of connected components in $F$ is not two, we do nothing, that is $\hh := F$. If $F$ has two components and each of them consists of at least~2 vertices, let $F_1$ be the component of $F$ with smaller edge density and let $F_2$ be the other one (the ties can be broken arbitrarily). For $i = 1, 2$, we write $t_i := e(F_i)/\binom{v(F_i)}{2}$ for the edge density of $F_i$. If $t_i+ f_i(t_1, t_2)<0$ or $t_i+ f_i(t_1, t_2)>1$ for some $i\in[2]$, then we again set $\hh := F$. Otherwise we replace $F_i$ by a random binomial graph with edge probability $t_i + (k)_2^{-1}f_i(t_1, t_2)$. Hence $\hh$ is the vertex-disjoint union of these two (independently drawn) random graphs.

	Write $p_{1,1}:=p_1$, $p_{2,2}:=p_2$ and $p_{1,2}=p_{2,1}:=0$. Let $\ff \sim \G(k, \Gamma(p))$ be a random $k$-vertex graph sampled from $\Gamma(p_1, p_2)$. We get the same random graph by sampling independent uniform elements $\uu_v \in [2]$ for each $v \in [k]$ and then connecting each pair $u, v$ independently with probability $p_{\uu_u, \uu_v}$. By Lemma~\ref{lem:viewvelocity},
	\begin{equation}
		\label{eq:vel_p1p2}
		(\vel[\rul]\Gamma(p))_{\restriction \Omega_i\times\Omega_j} = \sum_{a\neq b, a,b\in[k]}(\Pc{ab\in \hh}{ \uu_a=i,\uu_b=j} - p_{i,j}), \quad i,j\in[2]\;.
	\end{equation}
	Whenever $i \ne j$, the condition $\uu_a=i,\uu_b=j$ implies that pair $ab$ is a non-edge in $\hh$ and thus each term in \eqref{eq:vel_p1p2} is zero. This implies that~\eqref{eq:zerooffdiag} holds for every $p \in [0,1]^2$.
	
	We further assume that $p = (p_1,p_2)\in S_2$ and prove \eqref{eq:f_appr}.
	Let $D_i$ be the edge density of the subgraph $\ff[i]$ of $\ff$ induced by $\{v : \uu_v = i\}$. A standard Chernoff bound implies that with probability at least $1-2\exp(-\alpha_1^2k/6)$ we have $v(\ff[i]) = (1 \pm \alpha_1)k/2$ for $i\in[2]$. We denote this event by $\ca$. Furthermore, $\Ec{e(\ff[i])}{ \ca} = \binom{(1\pm\alpha_1)k/2}{2}p_i \sim k^2p_i/8$ and, assuming that $\ca$ holds, again a standard Chernoff bound implies that with probability at least $1-2\exp(-\alpha_1^2(1-\alpha_1)^2k^2p_i/24)$ we have $e(\ff[i]) = (1\pm\alpha_1)\binom{(1\pm\alpha_1)k/2}{2}p_i$ and therefore $D_i=(1\pm\alpha_2)p_i$. Lastly, assuming $\ca$ holds, with probability at least 
	\[
	1 - \binom{v(\ff[i])}{2} \left( 1 - p_i^2 \right)^{v(\ff[i]) - 2} \ge 1 - \binom{k}{2} (1-p_i^2)^{(1 - \alpha_1)k/2 - 2} \;.
	\]
	we have that in $\ff[i]$ every two vertices have a common neighbour and hence $\ff[i]$ is connected. We denote the event such that all of the three events above hold by $\mathcal{C}$ and note that $\Prob(\mathcal{C}) = 1-\exp(-\Omega(\alpha_1^2 k) ) = 1 - o(1)$. In particular $\mathcal{C}$ together with \ref{cond:densdiff} implies that $D_2\geq p_2-\alpha_2> p_2-\delta/2\geq p_1+\delta/2 > p_1+\alpha_2 \geq D_1$ and also that $(D_1,D_2)\in B(\delta,(p_1,p_2))$ and therefore by \ref{cond:dball} and \ref{cond:0dens1} we have that $0\leq D_i + (k)_2^{-1}f_i(D_1,D_2)\leq 1$. 
	In view of \ref{cond:dball}, for large enough $k$ we have $p(1 \pm \alpha_2) \in S_1$, so that condition \ref{cond:unif_cont} implies $|f_i(1 \pm \alpha_2)p) - f_i(p)| = o(1)$, uniformly over all choices of $p \in S_2$.
	
	Hence for $i\in[2]$ and $a,b\in[k]$ (note that by symmetry of the rule the following does not depend on the actual choice of $ab$ as long as $\uu_a=\uu_b=i$), 
	\begin{align*}
		\Prob(ab\in \hh\mid \uu_a=\uu_b=i)&=\Pc{ab \in \hh }{\uu_a = \uu_b = i,\mathcal{C}}\cdot \Prob(\mathcal{C}) \pm \Prob(\overline{\mathcal{C}})\\
		&=\left( (1\pm\alpha_2)p_i + (k)_2^{-1}f_i\left((1\pm\alpha_2)p_1,(1\pm\alpha_2)p_2\right ) \right)\cdot \Prob(\mathcal{C}) \pm \Prob(\overline{\mathcal{C}})\\
		\justify{by \ref{cond:unif_cont}} &=((1\pm\alpha_2)p_i+(k)_2^{-1}f_i(p_1,p_2) + o(1))\cdot(1- o(1)) + o(1)\\
		\justify{for large enough $k$} &=p_i+(k)_2^{-1}f_i(p_1,p_2)\pm\eps/(2(k)_2)\;.
	\end{align*}
	
	Therefore, recalling \eqref{eq:vel_p1p2}, we finally get that
	\begin{equation*}
		(\vel[\rul]\Gamma(p))_{\restriction \Omega_i\times\Omega_i} = f_i(p_1,p_2) \pm \eps/2, \quad i \in [2]\;,
	\end{equation*}
	and in particular
	\[
	\sup_{p \in S_2} |f(p) - \Gamma^{-1}(\vel[\rul]\Gamma(p))| \le \eps.
	\]

\end{proof}

We now construct the set where we will find a periodic trajectory. Fix values $r=0.1$, $r_s=0.09$, $r_w=0.11$, $\theta < 0.0001$, $a=0.2$, $b=0.8$ and set $q=(a,b)$. Recall that functions $g,h:B(0.15,q)\setminus{q}\rightarrow \mathbb R^2$ were defined by~\eqref{eq:defG} and~\eqref{eq:defH}. Define a function $f:B(0.15,q)\setminus{q}\rightarrow \mathbb R^2$ by 
\[
f(x,y) = \theta \cdot (g(x,y) + h(x,y)).
\]

Note that the dynamical system over $B(0.15,q)\setminus{q}$ defined by $f$ has a periodic trajectory $\{(x,y)\in\R^2\mid (x-a)^2+(y-b)^2=r^2\}$. To finish our construction we apply Lemma~\ref{lem:2dimapproxrule} to the function $f$, with the constants
\[
\eps<\min\{|r-r_s|,|r-r_w|, \theta \} \quad\text{and} \quad \delta<0.005,
\]
and the sets $S_1=B(0.15,q)\setminus\overline{B(0.005,q)}$ and $S_2=B(0.145,q)\setminus\overline{B(0.01,q)}$.
Note that $f$ is continuous on the closure of $S_1$, which is compact, implying that condition~\ref{cond:unif_cont} holds. Verifying the other conditions is straightforward.
Hence by Lemma~\ref{lem:2dimapproxrule} there is a rule $\rul$ of order $k$ so that for the vector field $f_\rul:S_2\rightarrow\R^2$ defined by $f_\rul(p) = \Gamma^{-1}(\vel[\rul] \Gamma(p))$ we have 
\begin{equation}
	\label{eq:f_ftilde}
	|f(p)-f_\rul(p)|\leq\eps \quad \text{for all} \quad p\in S_2. 
\end{equation}

Using the continuity and linearity of $\Gamma$ and $\Gamma^{-1}$ we see that the trajectory $\trajRule{\rul}{\cdot}\Gamma(p)$ (bijectively) corresponds to the solution $x(t) = \Gamma^{-1} \left( \traj{t}{\Gamma(p)} \right)$ of the dynamical system defined by $F = f_\rul$ with initial condition $F(0) = p$. Of course, one is periodic if and only if the other is periodic. 
We will show that the dynamical system defined by $F$ has a periodic trajectory in a closed annulus using the Poincar\'e--Bendixson theorem, see \cite[Theorem 1, p.~245]{Perko2001}.

\begin{thm}\label{thm:PB}[Poincar\'e--Bendixson]
	Suppose $S\subseteq \mathbb R^2$ is an open set and a continuously differentiable function $F:S\rightarrow \R^2$ defines the dynamical system 
	\begin{equation}
		\label{eq:dyn_sys}
		\tag{$*$} \quad x'(t)=F(x(t))\;. 
	\end{equation}
	If $C\subseteq S$ is a compact set and
	\begin{romenumerate}
		\item\label{en:f1} there exists $p\in C$ and a trajectory with $x(0)=p$ such that $x(t)\in C$ for every $t\geq 0$, and
		\item\label{en:f2} there are no points $p\in C$ with $F(p)=0$,
	\end{romenumerate}
	then the system \eqref{eq:dyn_sys} has at least one periodic trajectory inside $C$.
\end{thm}

We apply Theorem~\ref{thm:PB} to the function $F = f_\rul$ and the sets $S = S_2$ and $C=\overline{B(r_w,q)}\setminus B(r_s,q)$ The following claim shows that our choices fulfil the conditions of Theorem~\ref{thm:PB}, implying the existence of a periodic trajectory.

\begin{figure}[t]
	\centering
	\begin{tikzpicture}[
		point/.style = {draw, circle, fill=black, inner sep=0.5pt}, scale=1.4
		]
		
		\def\rad{1cm}
		\def\radw{2.5cm}
		
		\node (P)  at +(55:\rad) [point]{};
		\node[scale=0.8]  at ($(0,0.2)+(55:\rad)$) {$p$};
		
		
		\draw (\rad,0) arc[start angle=0, end angle=90, radius=\rad];
		\draw (\radw,0) arc[start angle=0, end angle=90, radius=\radw];
		\node[scale=0.8] at ($(0,\rad)+(-0.6,0)$) {$\kruh(r_s,q)$};
		\node[scale=0.8] at ($(0,\radw)+(-0.6,0)$) {$\kruh(r_w,q)$};
		\node[scale=0.8] at ($(0,\rad/2+\radw/2)$) {$C$};
		
		\def\TangentAngle{10}
		
		\draw (0.2,0.7)
		to [out=30,in=180+\TangentAngle] (P)
		to [out=\TangentAngle, in=160] ($(P)+(0.8,-0.08)$);
		\node[scale=0.8]  at (0.2,0.5) {$x(t)$};

		\def\ArcAngle{40}
		\def\ArcRadius{1}
		\def\LableExtend{1.15}%
		
		\def\HorizAngle{\TangentAngle-\ArcAngle/2}
		
		\pgfmathsetmacro{\XValueTangent}{\ArcRadius*cos(\TangentAngle)}%
		\pgfmathsetmacro{\YValueTangent}{\ArcRadius*sin(\TangentAngle)}%
		
		\pgfmathsetmacro{\XValueHoriz}{\ArcRadius*cos(\HorizAngle)}%
		\pgfmathsetmacro{\YValueHoriz}{\ArcRadius*sin(\HorizAngle)}%
		
		\pgfmathsetmacro{\XValueArc}{\ArcRadius*cos(\ArcAngle+\HorizAngle)}%
		\pgfmathsetmacro{\YValueArc}{\ArcRadius*sin(\ArcAngle+\HorizAngle)}%
		
		\pgfmathsetmacro{\XValueLabel}{\ArcRadius*cos(\HorizAngle+\ArcAngle/2)}%
		\pgfmathsetmacro{\YValueLabel}{\ArcRadius*sin(\HorizAngle+\ArcAngle/2)}%
		
		\draw [->, thin, dashed] (P) -- ($(P)+(\XValueTangent,\YValueTangent)$);
		\node[scale=0.8] at ($(P)+(\XValueTangent,\YValueTangent)+(0.3,0)$) {$F(p)$};
		
		\draw [->, thin, dotted] (P) -- ($(P)+(\XValueHoriz,\YValueHoriz)$);
		
		\draw [->, thin, dotted] (P) -- ($(P)+(\XValueArc,\YValueArc)$);
		
		\draw [dotted] ($(P)+(\XValueHoriz,\YValueHoriz)$)
		arc (\HorizAngle:(\HorizAngle+\ArcAngle):\ArcRadius);
		
		\node[scale=0.8] at ($(P)+(\XValueArc,\YValueArc)+(0,0.2)$) {$D_\delta$};

	\end{tikzpicture}
	\caption{The trajectory $x(t)$ entering the set $C$ at the point $p$, is approximated by the tangent $F(p)$ and thus stays in the wedge-shaped set $D_\delta$.}
	\label{fig:claiminterior}
\end{figure}
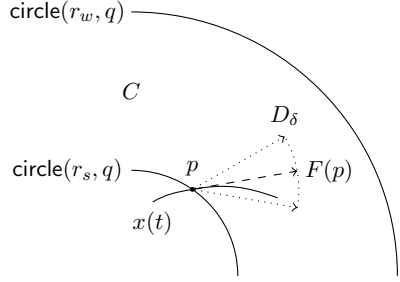

\begin{claim}\label{clm:cond_PB}
	The following is true for $F = f_\rul$.
	\begin{romenumerate}
		\item\label{en:IlB} $F$ \emph{points inwards} at every point $p$ on the boundary of $C$, that is, $F(p)$ is not collinear with the tangent of the boundary at $p$ and $p + \delta F(p) \in C$ for small enough $\delta >0$.
		\item\label{en:IlC} If $p$ lies on the boundary of $C$, that is, $\kruh(r_w, q) \cup \kruh(r_s, q)$, and $x(T) = p$ for some trajectory defined by \eqref{eq:dyn_sys} and some $T > 0$, then for some $\tau > 0$ we have $x([T, T + \tau]) \in C$.
		\item\label{en:IlD} Condition \ref{en:f1} of Theorem~\ref{thm:PB} holds (in fact, for every $p \in C$).
		\item\label{en:IlA} Condition \ref{en:f2} of Theorem~\ref{thm:PB} holds.
	\end{romenumerate}
\end{claim}

\begin{proof}
	
	For~\ref{en:IlB}, assume first that $(x,y)=p\in \kruh(r_s,q)$. The component $g(p)$ of $f(p)$ is tangent to $\kruh(r_s,q)$ and the component $h(p)$ points away from the point $q$. Since $p\in \kruh(r_s,q)$, we have $\sqrt{(x-a)^2+(y-b)^2}=r_s$, therefore $|h(p)|=|r-r_s|>\eps$ and hence, in view of \eqref{eq:f_ftilde}, vector $F(p)$ points inwards.
	An analogous argument shows that $F(p)$ also points inwards for every $p \in \kruh(r_w, q)$.
	
	We now prove \ref{en:IlC}. By \ref{en:IlB}, vector $F(p)$ forms a non-zero angle $\gamma \le \pi/2$ with the tangent at $p$ to a boundary circle and by adding to $p$ a small multiple of $F(p)$ we end up in $C$. Let $D_\delta$ be a wedge-shaped set consisting of vectors of length at most $\delta$ at angle at most $\gamma/2 (\le \pi/4)$ with $F(p)$ (see Figure~\ref{fig:claiminterior}). For sufficiently small $\delta$ we have that $p + y \in C$ for every $y \in D_\delta$. Since the trajectory $x$ has the derivative $x'(t) = F(p)$, we have 
	\[
	|x(t+\tau) - (p + \tau F(p) ) | = o(\tau)\;.
	\]
	Note that for small enough $\tau$ (depending on $\delta$) the vector $\tau F(p)$ is $\Omega(\tau)$-far from the edge of $D_\delta$. This implies that $x(t+\tau) - p \in D_\delta$, for small enough $\tau$, and hence $x(t+\tau) \in C$.
	
	Now we prove \ref{en:IlD}. Pick an arbitrary $p \in C$ and let $x(t)$ be the solution of \eqref{eq:dyn_sys} with the initial condition $x(0) = p$. Writing $T = \sup\{ t \ge 0 : x([0,t])\subset C \}$, we need to show $T = \infty$. Suppose the contrary, that is, $T < \infty$. Clearly, $x([0,T)) \subset C$, and  $x(T)$ lies on the boundary of $C$ (as otherwise the trajectory would stay inside $C$ a bit longer). By~\ref{en:IlC} there exists $\tau>0$ such that $x([0,T+\tau])\subseteq C$ contradicting the maximality of $T$.
	
	Finally, to prove \ref{en:IlA}, fix $p \in C$. Since $C \subset S_2$, $g(p)$ is perpendicular to $h(p)$, and $|g(p)|=1$, we get that 
	\[
	|F(p)| \geBy{\eqref{eq:f_ftilde}} |f(p)| - \eps = \theta |g(p) + h(p)| - \eps \ge  \theta |g(p)| - \eps = \theta - \eps > 0\;,
	\]
	which proves \ref{en:IlA}.
\end{proof}

\begin{remark}\label{rem:convergencetooscilatory}
	We stated the Poincar\'e--Bendixson theorem in an abridged form, sufficient for our purposes. Using the original version of the Poincar\'e--Bendixson theorem (see \cite[Theorem 1, p.~245]{Perko2001}) and Claim~\ref{clm:cond_PB}.\ref{en:IlD}, it actually follows that every trajectory starting in $C$ is periodic or spirals towards a periodic orbit.
\end{remark}

We note that our proof of the existence of a periodic trajectory is not constructive, since the Poincar\'e--Bendixson theorem is not. We leave it as an open problem to find an explicit construction.

\subsection{Complicated trajectories}\label{ssec:complicatedtrajectories}
Given a rule, a natural question about a rule and an initial graphon $W$ is to understand which graphons a trajectory approaches infinitely often, namely to determine the set 
\[
\omega_+(W) := \{U \in \Gra : \cutn{\traj{t_n}W - U} \to 0 \text{ for some sequence } t_n \to \infty\}. \]
Convergent trajectories are simplest in the sense that for them $\omega_+(W)$  contains only $\dest(W)$. Even though the compactness of $\widetilde{\Gra}$ implies that every trajectory approximates some graphon $U$ in the cut \emph{distance} $\cutm$, we still do not have an answer to the following question.
\begin{qu}\label{qu:limit_point}
	Is it true that, for every rule $\rul$, and every graphon $W$, set $\omega_+(W)$ is nonempty?
\end{qu}
We have seen examples where $\omega_+(W)$ is a continuum set, namely a periodic trajectory from Section~\ref{ssec:oscillatory} and a spiraling trajectory in Remark~\ref{rem:convergencetooscilatory}.

The following question has been suggested to us by Joel Spencer.
\begin{qu}\label{qu:spencer}
	Does there exist a flip process and an initial graphon $W$ so that the edge density oscillates between constant graphons~0 and~1, that is $\{0, 1\} \in \omega_+(W)$?
\end{qu}

What we find the most intriguing open question in this project is to find flip processes and graphons with really complicated trajectories. There is a topological way to express this (and it would be interesting to see alternative notions of complexity of trajectories). Recall that the metric $\cutnd(U,W) = \cutn{U - W}$ induces a non-compact topology on $\Gra$. We say that the trajectory starting at $W$ is \emph{complicated} if $\{\traj{t}{W}: t\ge0\}$ is not totally bounded, or equivalently, the closure of its orbit, that is, $\overline{\{\traj{t}{W}: t\ge0\}}$ is not compact with respect to the metric $\cutnd$. Then the main question is this.
\begin{qu}\label{qu:noncompacttrajectory}
	Is there a flip process and a graphon $W$ such that $\overline{\{\traj{t}{W}: t\ge0\}}$ is not compact?
\end{qu}
A negative answer to Question~\ref{qu:limit_point} (which would be rather surprising) would give a positive answer to Question~\ref{qu:noncompacttrajectory}.

The proposition below excludes convergent trajectories from the search.
\begin{prop}\label{prop:convsinkNotCompl}
	If a rule and a graphon $W$ are such that $\dest(W)$ exists, then $\overline{\{\traj{t}{W}: t\ge0\}}$ is compact.
\end{prop}
\begin{proof}
	Let $\rul$ be a rule of order $k$ and $W$ be the initial graphon with $U := \dest(W)$. We need to prove total boundedness, that is, for any given $\eps>0$ there exist finitely many $\cutnd$-balls $B_1,\ldots,B_h$ in $\Gra$ of radius $\eps$ so that $\bigcup B_i\supset \{\traj{t}{W}: t\ge0\}$. Since $U = \dest(W)$, there exists $t_0\ge 0$ such that for each $t\ge t_0$, $\cutnd(\traj{t}{W},U)<\eps$. Now, the appropriate balls have centers $U$ and $\traj{\tau}{W}$ where $\tau$ ranges over multiples of $\eps/k^2$ from~0 until (it exceeds) $t_0$. To prove the statement of the proposition, we need to argue that, for any $t\ge 0$, the graphon $\traj{t}{W}$ is covered by one of the balls above. If $t\ge t_0$, then $\traj{t}{W}$ is covered by the ball around $U$. If $t\in [0,t_0)$, then $\traj{t}{W}$ is covered by the ball around $\traj{\tau}{W}$ for the largest $\tau$ with $\tau\le t$. To see this, recall that by~\eqref{eq:Phi_t_Lipschitz}, $\traj{\tau}{W}$ and $\traj{t}{W}$ at most $(k)_2(t-\tau)\le\eps$ apart in the $L^\infty$-norm (and hence also in $\cutnd$).
\end{proof}

The strongest possible (positive) answer to Question~\ref{qu:noncompacttrajectory} (and also to Question~\ref{qu:spencer}) would be a trajectory whose closure covers the entire graphon space. However, there is no such trajectory.
\begin{prop}\label{prop:nopeanotrajectory}
	There is no flip process and a graphon $W$ such that $\overline{\{\traj{t}{W}: t\ge0\}}=\Gra$ \;.
\end{prop}
\begin{proof}
	Suppose that such a flip process exists, say with a rule $\rul$ of order $k$, and let $W$ be a graphon whose closure of the trajectory covers $\Gra$. 
	
	Graphon $W$ cannot have a set of twins of positive measure. Indeed, if it had, Proposition~\ref{prop:twinsstay} and Lemma~\ref{lem:twinsconverge} would tell us that we can approximate only graphons with that set of twins, making it impossible to approximate graphons without a set of twins of positive measure. In particular,
	\begin{equation}\label{eq:norepres}
		\mbox{$W$ is not a graphon representation of any $k$-vertex graph.}
	\end{equation}
	
	The following claim is the core of the proof.
	\begin{claim*}
		For every $G\in\lgr{k}$ we have $\rul_{G, G}=1$.
	\end{claim*}
	\begin{proof}[Proof of Claim]
		Suppose that this is not the case, that is, there exists a graph $H\in\lgr{k}\setminus\{G\}$ such that $p:=\rul_{G,H}>0$. Fix a pair $ij\in E(G)\sdiff E(H)$. Let us assume that~(*) $ij\notin E(G)$ and $ij\in E(H)$. The argument in the case~(**) $ij\in E(G)$ and $ij\notin E(H)$ is symmetric and we will not include it. However, to clarify the adjustments needed and also for its own sake, let us give an overview of the proof of~(*). We shall consider a graphon representation $W_G$ of $G$ with respect to a fixed partition $(\Omega_v)_{v\in[k]}$. Our focus will be on the square $\Omega_i\times\Omega_j$ where $W_G$ is constant~0. The core argument (expressed in~\eqref{eq:veryimportantvelocity}) is that if the trajectory $U=\traj{t}{W}$ gets somewhat close to $W_G$ then necessarily $\int_{\Omega_i\times\Omega_j}\vel U$ is positive (with an explicit positive lower bound) thus repelling from $W_G$ on $\Omega_i\times\Omega_j$. Then inapproximability of $W_G$ follows. Obviously, if we have (**) instead of (*) then the counterpart of \eqref{eq:veryimportantvelocity} is a negative upper bound on $\int_{\Omega_i\times\Omega_j}\vel U$ showing that an ambition of having a trajectory approaching $W_G$ must fail due to $(W_G)_{\restriction\Omega_i\times\Omega_j}=1$.
		
		Let $W_G$ be the graphon representation of $G$ with respect to an arbitrary fixed partition $(\Omega_v)_{v\in[k]}$. 
		By~\eqref{eq:norepres} we have that $\traj{0}{W}\neq W_G$. Also, as $W_G$ is $\{0,1\}$-valued, we cannot have $\traj{t}{W}=W_G$ for any $t>0$ by Proposition~\ref{prop:valueszeroone}. 
		In particular this means that there are times $t_1' < t_1 < t_2' < t_2 < \ldots$ such that 
		\begin{equation}
			\label{eq:bangos}
			\cutnd{(W_G,\traj{t_\ell'}{W})}>\frac{1}{4}, \quad  \cutnd{(W_G,\traj{t_\ell}{W})}<\frac{1}{\ell}.
		\end{equation}
		To see this, note that there are two $\{0, 1\}$-valued graphons with distance more than $1/2$ between them (say, representation of complete and empty graphs); hence at least one of them, $F$, is such that $\cutnd{(W_F, W_G)} > 1/4$; since $W_F$ also does not lie on the trajectory, the trajectory has to alternate between  small neighbourhoods of $W_G$ and $W_F$.
		
		Let $C_\square=C_\square(k,0)$ be given by Lemma~\ref{lem:square_Lip}. Set $\delta := pk^{-k}/(2C_\square)$. We claim that for every graphon $U$ with $\cutnd(U,W_G) \le \delta$ we have \begin{equation}\label{eq:veryimportantvelocity}
			\int_{\Omega_i\times\Omega_j}\vel U\ge \frac{p k^{-k}}2\;.
		\end{equation}
		To this end, first we compute $\vel W_G$, focusing on the value on the square $\Omega_i\times\Omega_j$ (on which $\vel W_G$ is constant $0$). So, let us fix $x\in \Omega_i$, $y\in \Omega_j$ and look at the summands in the representation~\eqref{eq:vel_prob}. Since $W_G(x,y)=0$, each term in \eqref{eq:vel_prob} is nonnegative. Keeping the term for $ab = ij$, we get (recall that $\uu_v$s are independent $\pi$-random variables)
		\begin{align*}
			\vel W_G(x,y) &\ge \Pc{\G(k,W_G) = G}{\uu_i = x, \uu_j = y}\rul_{G,H} \\
			&\ge \prob{\uu_v \in \Omega_v \forall v \in [k]} \rul_{G,H} = k^{-k}p.
		\end{align*}
		By Lemma~\ref{lem:square_Lip}, $\int_{\Omega_i\times\Omega_j}U\ge \int_{\Omega_i\times\Omega_j}W_G-C_\square\cutnd(U,W_G)\ge p k^{-k}/2$.
		This proves~\eqref{eq:veryimportantvelocity}.
		
		Let us show how~\eqref{eq:veryimportantvelocity} leads to a contradiction, thus proving the Claim. 
		Letting $\tau_\ell$ be the last time (see \eqref{eq:bangos}) smaller than $t_\ell$ such that $\cutnd{(W_G,\traj{\tau}{W})}\ge \delta$, we have that 
		\[
		\cutnd{(\traj{\tau_\ell}{W},\traj{t_\ell}{W})}\ge \delta-\frac{1}{\ell}\;,
		\]
		and hence~\eqref{eq:Phi_t_Lipschitz} (taking into account that $\cutn{\cdot} \le \Linf{\cdot}$) implies that $t_\ell-\tau_\ell\ge (\delta - 1/{\ell})/{(k)_2}$. Starting with the integral form for the trajectory (see Remark~\ref{rem:integral_form}), $\traj{t_\ell}{W}=\traj{\tau_\ell}{W}+\int_{\tau_\ell}^{t_\ell} \vel\traj{t}{W}\D t$ and making use of the fact that for each $t\in [\tau_\ell,t_\ell]$ we can apply~\eqref{eq:veryimportantvelocity} with $U = \traj{t}{W}$, we get
		\begin{align}
			\label{eq:not_too_small}
			\begin{split}
				\int_{\Omega_i\times \Omega_j}\traj{t_\ell}{W} &= \int_{\Omega_i\times \Omega_j}\left(\traj{\tau_\ell}{W}\right) + \int_{\Omega_i\times \Omega_j}\left(\int_{\tau_\ell}^{t_\ell} \vel\traj{t}{W}\D t\right)\\
				&\ge 0 + \frac{\delta - \frac{1}{\ell}}{(k)_2}\cdot \frac{p k^{-k}}2\;.
			\end{split}
		\end{align}
		To justify the last inequality, note that by definition of the integral (see Appendix~\ref{app:Banach_calc}) we can approximate (in $L^\infty$-norm and thus in the cut norm) $\int_{\tau_\ell}^{t_\ell} \vel\traj{t}{W}\D t$ by a linear combination of kernels taking values in $\left\{ \vel\traj{t}W : t \in [\tau_\ell, t_\ell]  \right\}$ and applying \eqref{eq:veryimportantvelocity} to each of them.
		
		By choosing $\ell$ large enough we get that 
		$  \int_{\Omega_i\times \Omega_j}\traj{t_\ell}{W} > \frac{1}{\ell}.$
		On the other hand, $ij\notin E(G)$ implies
		$\int_{\Omega_i\times \Omega_j}W_G = 0$.
		We conclude that 
		\[
		\cutnd(\traj{t_\ell}{W}, W_G)\ge \int_{\Omega_i \times \Omega_j} \left( \traj{t_\ell}W - W_G \right) 
		> 1/\ell,
		\]
		contradicting \eqref{eq:bangos}.
	\end{proof}
	The Claim indeed proves the main statement, since by applying it to all graphs $G\in\lgr{k}$ we deduce that the flip process in question is trivial. This is clearly a contradiction. 
\end{proof}

\begin{remark}\label{rem:strengtheNoPeano}
	A strengthening of Proposition~\ref{prop:nopeanotrajectory} can be proved for the factor space $\widetilde {\Gra}$ (recall Remark~\ref{rem:factorspace}). Namely, there is no trajectory whose closure in the cut distance covers $\widetilde {\Gra}$. This follows by making minor changes to the proof of Proposition~\ref{prop:nopeanotrajectory}: in \eqref{eq:bangos} we allow the graphon representation $W_G$ to depend on $\ell$.
\end{remark}
\begin{remark}\label{rem:strengtheNoPeano2}
	In this section, we looked at the closure of the \emph{forward} orbit $\overline{\{\traj{t}{W}: t\ge0\}}$. As an alternative, we could look at the closure of the portion of the ``whole'' orbit $\overline{\{\traj{t}{W}: t\in\life{W}\}}$. Obviously, the qualitative difference between these two images can be only in the case $\age(W)=\infty$. Question~\ref{qu:noncompacttrajectory} is also open in this alternative setting as is the question about non-emptiness of the set
	\[
	\omega_-(W) := \{U \in \Gra : \cutn{\traj{t_n}W - U} \to 0 \text{ for some sequence } t_n \to -\infty\}. \]
	Nevertheless, a counterpart to Proposition~\ref{prop:convsinkNotCompl} holds if we assume, in addition, that $\traj{t}W$ converges as $t \to -\infty$. Also, Proposition~\ref{prop:nopeanotrajectory} holds even in this stronger setting. Indeed, the repulsion argument used to show the inapproximability of $W_G$ holds uniformly even at negative times. 
\end{remark}

\subsection{Speed of convergence}\label{ssec:speed}
For now let us focus on convergent flip processes. For each such corresponding rule $\rul$, we can define for each $\delta>0$ \emph{lower convergence time}
\begin{equation*}
	\tau_{\rul}^-(\delta) = \sup \{\inf\{T:\cutnd(\traj{T}{W},\dest_{\rul}(W))<\delta \}:W\in \Gra\}\;,\\
\end{equation*}
and the \emph{upper convergence time},
\begin{equation*}
	\tau_{\rul}^+(\delta) = \sup \{\sup\{T:\cutnd(\traj{T}{W},\dest_{\rul}(W))\ge \delta \}:W\in \Gra\}\;.
\end{equation*}

\begin{remark}\label{rem:exactconstantsinspeed}
	The definition of the cut norm as 
	\[
	\cutn{W}=\sup_{A,B\subseteq\Omega}\left|\int_{A\times B}W \D\pi^2\right|
	\]
	is not the only natural one, as sometimes it makes more sense to work only with 
	\[
	\cutn{W}=\sup_{A\subseteq\Omega}\left|\int_{A\times A}W \D\pi^2\right| \quad\text{or}\quad 
	\cutn{W}=\sup_{A,B\subseteq\Omega: A \cap B = \emptyset}\left|\int_{A\times B}W \D\pi^2\right|.
	\]
	As is well known (see, e.g., \cite{BCLSV}), these alternative definitions are off from the original one by a factor of at most~4. Thus, when pinning down the speed of convergence precisely, it matters which of these definitions is used.
\end{remark}
\begin{remark}
	In view of Conjecture~\ref{conj:sinkconvergeL1}, we might also define the speed using the $L^1$-distance.
\end{remark}
Below, we show that $\tau_{\rul}^-(\delta)$ is finite for flip processes with a single destination.
\begin{prop}
	If $\rul$ is a convergent flip process with a single destination, then for every $\delta\in (0,1]$ we have $\tau_{\rul}^-(\delta)<\infty$.
\end{prop}
\begin{proof}
	Let $U$ be the unique destination (and thus the unique fixed point). Hence by Proposition~\ref{prop:atleastonefixed} graphon $U$ is constant.
	
	Suppose for contradiction that the statement does not hold. Then for each $\ell\in\N$ we have a graphon $W_\ell$ such that $\inf\{T:\cutnd(\traj{T}{W_\ell},U)<\delta \}>\ell$. By the Lov\'asz--Szegedy compactness theorem, there exists a graphon $Z$ which is an accumulation point of $W_1,W_2,\ldots$ in the cut distance. Without loss of generality, let us assume that $Z$ is actually a cut distance limit of $W_1,W_2,\ldots$. Since the trajectory of $Z$ converges to $U$, there exists a time $T$ such that $\cutnd(\traj{T}{Z},U)<\delta/2$. Combining this with Theorem~\ref{thm:genome}, for any $W$ sufficiently close to $Z$ in the cut norm, we have $\cutnd(\traj{T}{W},U)<\delta/2$. Since $U$ is a constant graphon, by Section~\ref{ssec:repermuting}, the same conclusion holds for any $W$ sufficiently close to $Z$ in the cut distance. Taking $W=W_\ell$ for sufficiently large $\ell$ gives a contradiction.
\end{proof}
The same question is open for the quantity $\tau_{\rul}^+(\delta)$.

Sometimes the convergence time is infinite for trivial reasons. For example if, given a rule $\rul$, there is a graphon $W$ with $\age(W) = \infty$ such that $\lim_{t \to -\infty} \traj{t}W$ exists and $\vel W \neq 0$, then one can make the convergence time arbitrarily long by choosing the initial graphon as $\traj{-t}{W}$ with $t$ large enough.
\begin{qu}
	Suppose that $\rul$ is a convergent flip process in which every graphon is either a fixed point or has finite age. Is it then true that for each $\delta\in (0,1]$ we have $\tau_{\rul}^-(\delta)<\infty$? Is the same conclusion true if, in addition, $\rul$ has only finitely many destinations?
\end{qu}

We can go back to our motivating examples. Indeed,~\eqref{eq:solutionErdosRenyi} tells us that for the Erd\H{o}s--R\'enyi flip process (whose only destination is the constant-1 graphon) we have $\tau^+(\delta)=\ln(-\delta)/2$. The triangle removal flip process (generalised in~Section~\ref{ssec:removal}) is much more interesting. In particular, in~\cite{Flip2journal} a uniform bound on convergence time is derived using the Szemer\'edi's regularity lemma, thus giving extremely poor bounds.

It can be shown that for any flip process $\rul$ with a single destination we have a lower bound $\tau^+(\delta)=\Omega(\ln(-\delta))$. Let us discuss the possibility of any general upper bounds. Note that there is no `slowest process'. Indeed, we can slow down any rule $\rul$ by any factor $C>1$. That is, we define a new rule $\rul'$ which, for any drawn graph $H$, replaces it with $H$ with probability $1-1/C$ and replaces it with a graph sampled according to $\rul$ with probability $1/C$. Still, this slowing down does not change the order of magnitude of the function $\tau_{\rul}^-(\delta)$ as $\delta\rightarrow 0$. Hence, we could ask what is the slowest order of magnitude of convergence. 

\subsubsection{Variants} Many variants of these questions exist. For example, one might consider local versions of convergence times, such as the following. Suppose that we have a convergent flip process with a stable destination $U$ of radius $\gamma^*$. Then define
\[
\tau_{\rul}^-(\delta,\gamma)=\sup \{\inf\{T:\cutnd(\traj{T}{W},U)<\delta \}:W\in \Gra, \cutnd(W,U)<\gamma\}\;,
\]
for $\gamma^*>\gamma>\delta>0$.

\subsection{Uniqueness of the rule}\label{ssec:nonuniqueness}
 Hng~\cite{KeatUniquenessArxiv} provides a complete characterisation of equivalence classes of flip process rules (of the same order) with the same trajectories. As an application, Hng demonstrates that the symmetric deterministic rules\footnote{A rule $\rul$ of order $k$ is \emph{symmetric} if we have for each permutation $\psi$ on $[k]$ (which naturally acts on $\lgr{k}$) and each $H,J\in\lgr{k}$ that $\rul_{H,J}=\rul_{\psi(H),\psi(J)}$. A rule $\rul$ is \emph{deterministic} if for each $H\in\lgr{k}$ we have that $\rul_{H,J}=0$ for all except one $J\in\lgr{k}.$} and the rules of order~$2$ are precisely the rules which are unique in their equivalence classes. The class of symmetric deterministic rules includes several examples from Section~\ref{sec:examples}, namely the complementing rules, the component completion rules, the extremist rules, and the clique removal rules.

\subsection{Large deviations}\label{ssec:largedeviations}
Our Theorem~\ref{thm:conc_proc} provides a description of the typical evolution for flip processes. It would be interesting to work out a large deviation theory for flip processes as well. For flip processes of order~2 such a large deviation theory was worked out in~\cite{BdHM:LargeDeviations_article}, but in a variant with updates along a Poisson point process (cf.~Section~\ref{ssec:Poisson}).\footnote{As we explain in Section~\ref{ssec:Poisson} there is essentially no difference between updates at integer times and along a Poisson point process for the typical behaviour. However, there is usually a difference between these two regimes from the large deviations point of view.}

\begin{appendix}
	
	\section{Probability}
	\label{ssec:MCDiarmid}
	We recall the Hoeffding--Azuma martingale concentration inequality including necessary preliminaries. To simplify some definitions we assume that the sample space $\Omega$ of the underlying probability space $(\Omega, \cf, \Prob)$ is finite. Consider a sequence of sigma-algebras $\{\emptyset,\Omega\}= \cf_0\subset \cf_1\subset \ldots \subset \cf_t$, such that $\cf_i \subseteq \cf$, which we call a \emph{filtration}.
	We say that a sequence of random variables $Y_1, \dots, Y_t$ is a \emph{martingale difference sequence} if $Y_i$ is $\mathcal{F}_i$-measurable and $\Ec{Y_i}{\mathcal{F}_{i-1}} = 0$ for $i=1,\ldots,t$. 
	
	The following lemma is a special case of Theorem~3.10 in~\cite{McDiarmid1998}.
	\begin{lem}
		\label{lem:MBoundDiff}
		Let $c$ be a constant. If $Y_1, \dots, Y_t$ is a martingale difference sequence with $|Y_i| \le c$ for each $i$, then for any $x \ge 0$,
		\[
		\prob{\left|Y_1 + \dots + Y_t\right| \ge x } \le \exp\left(-\frac{x^{2}}{2 t c^2}\right)\;.
		\]
	\end{lem}

	\section{Graphons}
	\label{app:graphons}
	
	Here we include some more facts about graphons and the proofs for Section~\ref{ssec:graphons}.
	
	The next two facts give alternative ways of bounding the cut norm. Variants appear in the literature, but we provide proofs for completeness.
	\begin{fact}[Inequality (7.1) in \cite{BCLSV}]
		\label{fact:cutnorm1}For each $W\in \Kernel$, 
		we have $
		\cutn{W}\le 2\sup_{A\subseteq\Omega}\left|\int_{A\times A}W \D\pi^2\right|$.
	\end{fact}
	\begin{proof}
		To bound $\cutn{W}$, let $A$ and $B$ be arbitrary as in~\eqref{eq:defcutn}. Using symmetry of $W$,
		\begin{align*}
			\left|\int_{A \times B} W\right| &= \frac{1}{2}\left|\int_{A \times B} W + \int_{B \times A} W\right| = \frac{1}{2}\left| \int_{(A \cup B)^2} W + \int_{(A \cap B)^2} W - \int_{(A \setminus B)^2} W - \int_{(B \setminus A)^2} W \right| \\ 
			&\le \frac12 \cdot 4 \sup_{C\subseteq \Omega} \left|\int_{C \times C} W\right|\;.
		\end{align*}
	\end{proof}
	
	\begin{fact}\label{fact:cutnorm2}Suppose that $\mathcal{A}$ is a partition of $\Omega$ into finitely many sets. For each $W\in \Kernel$ which is constant on each cell of $\mathcal{A}\times\mathcal{A}$ we have
		\begin{equation}
			\label{eq:step_graphon_cutnorm}
			\cutn{W} \le  \sup_{S,T\subseteq\mathcal{A}}\left|\int_{\bigcup S\times \bigcup T}W \D\pi^2\right|\;.
		\end{equation}
	\end{fact}
	\begin{proof}
		To bound $\cutn{W}$, let $A$ and $B$ be arbitrary as in~\eqref{eq:defcutn}. We shall prove that for an arbitrary $D\in\mathcal{A}$ we have 
		\begin{equation}\label{eq:oneoftwo}
			\left|\int_{A\times B}W\right| \le \max \left\{ \left|\int_{(A\cup D)\times B}W\right|, \left|\int_{(A\setminus D)\times B}W\right| \right\}\;.
		\end{equation}
		A similar `rounding' statement could be obtained for $B$ as well. Once this is done,~\eqref{eq:step_graphon_cutnorm} follows by rounding each of $A$ and $B$ on each cell of $\mathcal{A}$.
		
		Let $\alpha:=\pi(A\cap D) / \pi(D)$. Since $x \mapsto \int_{y \in B} W(x, y)$ is constant on $D$, we have
		\[
		\int_{A\times B}W  
		=
		\alpha\cdot \int_{(A\cup D)\times B}W \;+\;(1-\alpha)\cdot\int_{(A\setminus D)\times B}W\;.
		\]
		Since $\alpha \in [0,1]$, inequality \eqref{eq:oneoftwo} follows from the fact that $x \mapsto |x|$ is a convex function.
	\end{proof}
	
	Here, we give a self-contained proof of Lemma~\ref{lem:complete}. This proof was suggested to us by Jan Greb\'ik and Guus Regts.
	\begin{proof}[Proof of Lemma~\ref{lem:complete}]
		The proof shall use the concept of weak* convergence of graphons (which we treat as elements of $L^\infty(\Omega^2)$ with a predual space $L^1(\Omega^2)$), which we quickly recall. We refer the reader to~\cite{MR4192828} for  details. A sequence of graphons $W_1, W_2,\ldots$ \emph{converges weak*} to a graphon $W$ if for every $S,T\subset \Omega$ and every $\eps>0$ there exists a number $N_{S,T,\eps}\in\N$ such that 
		\begin{equation}\label{eq:cauchy}
			\left|\int_{S\times T} W_n- \int_{S\times T}W \right|<\eps \quad \text{for every $n\ge N_{S,T,\eps}$}\;.
		\end{equation}
		So, the cut norm topology can be regarded as a certain uniformisation of the weak* topology, and indeed convergence in the former topology implies convergence in the latter. 
		It follows from the Banach--Alaoglu theorem that the closed unit ball in $L^\infty(\Omega^2)$ is compact with respect to the weak* topology. Since the set of graphons $\Gra$ is contained in this unit ball and is closed with respect to the $L^\infty$ norm (and therefore closed in the weak* topology), $\Gra$ is compact in the weak* topology, too.
		
		We are now ready to prove the lemma. Suppose that $W_1,W_2,\ldots$ is a Cauchy sequence in $\cutnd$. That is, for every $\eps>0$ there exists an $M_\eps \in \N$ such that 
		\begin{equation}\label{eq:cauchy2}
			\sup_{S,T\subset \Omega}\left|\int_{S\times T}W_n-\int_{S\times T}W_m\right|<\eps \quad \text{for every $n,m\ge M_\eps$}\;.
		\end{equation}
		From compactness of $\Gra$ and \eqref{eq:cauchy2} it easily follows that $(W_n)$ converges weak* to some graphon $W$. We claim that $W$ is the sought cut norm limit. Suppose that $\eps>0$ is arbitrary. We then have that for $n \ge M_\eps$
		\begin{align*}
			\cutnd(W_n, W) &=\sup_{S,T\subset \Omega} \left|\int_{S\times T}W_n-\int_{S\times T}W\right|
			\\
			&\le
			\sup_{S,T\subset \Omega} \bigg\{\left|\int_{S\times T}W_n-\int_{S\times T}W_{\max\{M_\eps,N_{S,T,\eps}\}}\right| \\
			&+\left|\int_{S\times T}W_{\max\{M_\eps,N_{S,T,\eps}\}}-\int_{S\times T}W\right|\bigg\}
			\;,
		\end{align*}
		where the terms $N_{S,T,\eps}$ are defined as in~\eqref{eq:cauchy}. From \eqref{eq:cauchy} and \eqref{eq:cauchy2} it follows that $\cutnd(W_n, W) \le \eps+\eps$, which implies $\lim_{n \to \infty} \cutnd(W_n, W) = 0$, as was needed.
	\end{proof}
	
	Next, we turn to proving Lemma~\ref{lem:contdensLinfty}. We will make use of the following fact.
	\begin{fact}\label{fact:approximatingproduct}
		Given a real number $M > 0$, any numbers $a_1,\ldots,a_\ell,b_1,\ldots,b_\ell$ of modulus at most~$M$ satisfy the inequality
		\[
		\left| \prod_i a_i -\prod_i b_i \right| \le M^{\ell - 1} \sum_i |a_i-b_i|\;.
		\]
	\end{fact}
	\begin{proof}
		To apply induction on $\ell$, write $A = \prod_{i \le \ell - 1} a_i$, $B = \prod_{i \le \ell - 1} b_i$, so that the left-hand-side is $|Aa_\ell - Bb_\ell|$. We have $Aa_\ell - Bb_\ell=  A(a_\ell - b_\ell) + (A - B)b_\ell$. Hence, by the triangle inequality and the induction hypothesis
		\begin{align*}
			|Aa_\ell - Bb_\ell| 
			&\le |A||a_\ell - b_\ell| + |A-B||b_\ell| \\
			&\le M^{\ell-1}|a_\ell - b_\ell| + |A - B|M \\
			&\le M^{\ell - 1}\sum_{i = 1}^\ell |a_i - b_i|,
		\end{align*}
		as was needed.
	\end{proof}
	
	\begin{proof}[Proof of Lemma~\ref{lem:contdensLinfty}]
		This follows by fixing $(x,y)$, writing $(\Troot{F}{a,b}U)(x,y)$ and $(\Troot{F}{a,b}W)(x,y)$ as integrals (see \eqref{eq:tindr}) with the same variables of integration, applying the inequality $|\int f | \le \int |f|$ and then applying Fact~\ref{fact:approximatingproduct} to the integrand. (Note that $W(x_i,x_j)$ and $1 - W(x_i,x_j)$ have modulus at most $1 + \eps$ because of $W \in \Gra[\eps]$.)
	\end{proof}
	
	We want to move on to prove Lemma~\ref{lem:counting_kernels}. For this we need to cite a counting lemma from~\cite[Lemma~10.24]{Lovasz2012}.
	\begin{lem}[Counting lemma for vectors of graphons]\label{lem:counting}
		Given a simple graph $F = (V,E)$ and two tuples $\mathbf{U} = (U_e : e \in E)$ and $\mathbf{W} = (W_e : e \in E)$ of graphons, we have
		\[
		|t(F,\mathbf{U}) - t(F,\mathbf{W})| \le \sum_{e \in E} \cutn{U_{e} - W_{e}}.
		\]
	\end{lem}
	\begin{proof}[Proof of Lemma~\ref{lem:counting_kernels}]
		Every $W \in \Gra[\eps]$ can be written as $W = -\eps + (1 + 2\eps)\hat W = (1 + \eps) \hat W - \eps (1 - \hat W)$ for some graphon $\hat W \in \Gra$. We have, always assuming that the arguments of $W_{ij}$ and $\hat W_{ij}$ are $(x_i, x_j)$,
		\begin{align*}
			t(F,\mathbf{W})
			&= \int \prod_{e \in E} W_e \D \pi^{V} \\
			&= \int \prod_{e \in E} \left( (1 + \eps)\hat W_e - \eps(1 - \hat W_e) \right) \D \pi^V \\
			&= \sum_{S \subseteq E} (1 + \eps)^{|S|} (-\eps)^{|E| - |S|} \int \prod_{e \in S} \hat W_e \prod_{e \in E \sm S} (1 - \hat W_e) \D \pi^V \\
			&= \sum_{S \subseteq E} (1 + \eps)^{|S|} (-\eps)^{|E| - |S|} t(F, \mathbf{\hat W}_S)
			\;,
		\end{align*}
		where we defined $\mathbf {\hat W}_S=(\hat W_{S,e}:e\in E)$ by setting $\hat W_{S,e}=\hat W_{e}$ if $e \in S$ and $\hat W_{S,e}=1 - \hat W_e$ otherwise.  
		
		Analogously defining $\mathbf {\hat U}_S$ and expressing $t(F,\mathbf U)$, it follows that 
		\begin{align*}
			|t(F, \mathbf{U}) - t(F, \mathbf{W})| 
			&\le \sum_{S \subseteq E} (1 + \eps)^{|S|}\eps^{|E| - |S|} \left| t(F,  \mathbf {\hat U}_S) - t(F,  \mathbf {\hat W}_S) \right| \\
			\justify{Lemma~\ref{lem:counting}} &\le \sum_{S \subseteq E} (1 + \eps)^{|S|}\eps^{|E| - |S|} \sum_{e \in E} \cutn{\hat U_{S,e} - \hat W_{S,e}} \\
			\justify{$\hat U_{S,e} - \hat W_{S,e} = (U_e - W_e)/(1+ 2\eps)$} &= \sum_{S \subseteq E} (1 + \eps)^{|S|}\eps^{|E| - |S|} \cdot \frac{\sum_{e \in E} \cutn{U_e - W_e}}{1+ 2\eps}\\
			\justify{binomial theorem} &= \left( 1 + 2\eps \right)^{|E|-1} \sum_{e \in E} \cutn{U_e - W_e}.
		\end{align*}
	\end{proof}
	
	Lastly, we add the proof of Lemma~\ref{lem:twinsconverge}.
	
	\begin{proof}[Proof of Lemma~\ref{lem:twinsconverge}]
		Fix a measurable $T \subset \Omega$ and let $f_n(x) = \int_T U_n(x,\cdot) \D \pi$ and $f(x) = \int_T U(x,\cdot) \D \pi$. Given a set of positive measure $B$ and $g \in L^\infty(\Omega)$, let us write $\E_B g = \frac{1}{\pi(B)}\int_B g \D \pi$.
		Since $A$ is a twin-set for $U_n$, for every $B \subset A$ of positive measure we have
		\[
		\E_B f_n = \frac{1}{\pi(B)} \int_{B \times T} U_n = \frac{1}{\pi(A)} \int_{A \times T} U_n = \E_A f_n.
		\]
		Since $U_n \to U$ in the cut norm, this implies that 
		\begin{equation}
			\label{eq:averages}
			\E_B f = \frac{1}{\pi(B)} \int_{B \times T} U = \frac{1}{\pi(A)} \int_{A \times T} U = \E_A f.
		\end{equation}
		We claim that $f$ is constant a.e.\ on $A$. Let $\eps>0$ be arbitrary. Let $B = \left\{ x \in A\!:\!f(x) \ge \eps + \E_A f \right\}$. If $\pi(B) > 0$, then $ \E_B f \ge \eps + \E_A f$, contradicting~\eqref{eq:averages}. Since $\eps$ was arbitrary, $f \le \E_A f$ a.e., and by a similar argument $f \ge \E_A f$.
		
		Recalling that the sigma-algebra of $\Omega$ is separable, let  $\mathcal C$ be a countable set of measurable sets which approximate any measurable $S \subset \Omega$. Writing $f_T(x) = \int_{T} U(x, \cdot) \D \pi$, by what we showed above, for every $T \in \mathcal C$  we have $f_T = c_T := \E_A f_T$ a.e. Hence the set $A' := \left\{ x \in A : \forall T \in \mathcal{C}\ f_T(x) = c_T\right\}$ is a conull subset of $A$.
		
		Fix $x, y \in A'$ and $\eps > 0$ and let $B_\eps = \left\{ z \in \Omega : U(x,z) - U(y, z) \ge \eps \right\}$. For any $T \in \mathcal{C}$ we have 
		\[
		\eps \pi(B_\eps) \le \int_{B_\eps} \left( U(x,\cdot) - U(y, \cdot) \right) \D \pi \le \pi(B_\eps \sdiff T) \cdot 2\Linf{U} + \int_{T} \left( U(x,\cdot) - U(y, \cdot) \right) \D \pi \;.
		\]
		Since 
		$  \int_{T} \left( U(x,\cdot) - U(y, \cdot) \right) = f_n(x) - f_n(y) = 0$
		and $\pi(B_\eps \sdiff T)$ can be made arbitrarily small, we get that $\pi(B_\eps) = 0$, that is, $U(x, \cdot) \le U(y, \cdot) + \eps$ a.e. Since $x,y, \eps$ were arbitrary, this implies that $A$ is a twin-set for $U$.
	\end{proof}

	\section{Calculus of normed-space-valued functions}
	\label{app:Banach_calc}
	
	Let $\mathcal{E}$ be a Banach space with norm $\norm{\cdot}$.
	
	\begin{defi}
		\label{def:deriv_normed}
		Suppose that a function $f : I \to \mathcal{E}$ is defined on some subset $I \subset \R$. If there is an element $z \in \mathcal{E}$ such that
		\[
		\lim_{I \ni \eps \to 0} \norm {\frac{f(t +\eps) - f(t)}{\eps} - z } = 0,
		\]
		then we say that \emph{$f$ is differentiable at $t$} and call $\frac{\D f}{\D t}(t) = f'(t) = z$ the \emph{derivative of $f$ at $t$}.
	\end{defi}
	
	For any $a, b \in \R$, let $\cs([a,b], \mathcal{E})$ be the set of \emph{step functions} $f : [a,b] \to \mathcal{E}$, that is, functions such that $f$ is fixed on each set of some finite partition of $[a,b]$ into intervals $(P_1, \dots, P_n)$. For a step function $f$, the integral $\int_a^b f$ is an element of $\mathcal{E}$ defined by
	\[
	\int_a^b f = \sum_i |P_i| f_i, \quad \text{where } f(t) = f_i \text{ on } P_i.
	\]
	Considering the space of bounded functions from $[a,b]$ to $\mathcal{E}$ with uniform norm $\Linf{f} := \sup_{t \in [a,b]} \norm{f(t)}$, 
	it is easy to show that any continuous function $f : [a,b] \to \mathcal{E}$ can be approximated by step functions (see \cite[p. 54]{AbrMar88}), that is, there is a sequence $f_n \in \cs([a,b],\mathcal{E})$ such that $\lim_{n \to \infty} \Linf{f - f_n} = 0$. Moreover, for any such approximating sequence the elements $\int_a^b f_n$ converge to the same element of $\mathcal{E}$ which we denote by $\int_a^b f$ and call the \emph{Cauchy--Bochner integral of $f$ over $[a,b]$}. This integral satisfies the following properties
	\begin{equation}
		\label{eq:norm_of_int}
		\left\lVert \int_a^b f \right\rVert \le \int_a^b \|f\| \le |b-a| \Linf{f},
	\end{equation}
	\begin{equation}
		\label{eq:Boch_int_addit}
		\int_a^b f = \int_a^c f + \int_c^b f \quad 
		\text{with a convention that} \quad \int_a^b f := - \int_b^a f.
	\end{equation}
	
	For proofs below, whenever we approximate a continuous function $f$ by a step function $f^{(n)}$, $n \in \N$, we let $[a,b] = P_1 \cup \dots \cup P_n$ be a partition into intervals of equal lengths and assume that $f^{(n)}$ on $P_i$ takes some (constant) value which $f$ takes at some point in $P_i$. We may occasionally choose this value not arbitrarily; say if $f(t)$ has some property for almost every $t$, then we may want $f^{(n)}$ to inherit this property and to achieve this we will choose the value in $P_i$ avoiding the exceptional measure-zero set for $t$. 
	
	Since $f$ is continuous on a compact interval, it is also uniformly continuous. Hence the sequence $(f^{(n)})_n$ converges to $f$ in the uniform norm implying that $\big(\int f^{(n)}\big)_n$ converges to $\int f$.
	
	We now state the fundamental theorem of calculus for normed space-valued functions. 
	\begin{thm}[Theorem 2.3.11 in \cite{AbrMar88}]
		\label{thm:fund_thm_calc}~
			If $f : [a,b] \to \mathcal{E}$ is continuous, differentiable on $(a,b)$ and if $f'$ extends to a continuous function on $[a,b]$, then  
			\[
			f(b) - f(a) = \int_a^b f' \;.
			\]
		\end{thm}
		Further, given $f, g \in L^\infty(\Omega^2)$ we write $f \ge g$ meaning $f \ge g$ a.e.
		\begin{prop}
			\label{prop:Chebyshev}
			Given a continuous function $f : [a,b] \to L^\infty(\Omega^2)$ such that for every $t$ we have $f(t) \ge 0$, we have that 
			\[
			\int_a^b f(t) \D t \ge 0 \;.
			\]
		\end{prop}
		\begin{proof}
			Let $f_n$ be the explicit step functions defined above so that 
			\[
			I_n := \int_a^b f_n(t) \D t \toLinf \int_a^b f(t) \D t =: I.
			\]
			Since $f_n$ takes values in the range of $f$, for every $t$ we have that $f_n(t) \ge 0$ and therefore $I_n \ge 0$. Assume for contradiction that $I \le -\eps$ on a set $A$ of measure $\eps > 0$. Then
			\[
			-\eps^2 \ge \int_A I = \int_A I_n + \int_A (I_n - I) \ge 0 - \Lone{I_n - I} \ge - \Linf{I_n - I} \to 0,
			\]
			giving a contradiction.
		\end{proof}
		
		The following fact says that for $L^\infty$-valued integrals one can calculate the integral pointwise.
		\begin{prop}
			\label{prop:int_pointw}
			Let $\mathcal{E}$ be the normed space $L^\infty(\Omega^2)$. Consider a continuous function $[a,b] \ni t \mapsto f_t \in \mathcal{E}$. 
			One can modify, for every $t \in [a,b]$, $f_t$ on a set of measure zero (that possibly depends on $t$) so that for every $(x,y) \in \Omega^2$ function $t \mapsto f_t(x,y)$ is continuous (and thus integrable) and the function $J \in \mathcal{E}$ defined by $J(x,y) = \int_a^b [f_t(x,y) ] \D t$ represents the same element of $\mathcal{E}$ as $I := \int_a^b f_t \D t$. 
			
			Furthermore, if we have continuous functions $t \mapsto f_t$ and $t \mapsto g_t$ such that $f_t \ge g_t$ for every $t \in [a,b]$, then each of them can be modified to satisfy the conclusion above and also $f_t(x,y) \ge g_t(x,y)$ for every $(x,y) \in \Omega^2$ and $t \in [a,b]$.
		\end{prop}
		
		\begin{proof}
			Since $[a,b]$ is compact, $f_\cdot$ is uniformly continuous, and so there is a function $\omega$ such that $\omega(x) \to 0$, $x \to 0$, for which
			\begin{equation}
				\label{eq:unicont_norm}
				\Linf{f_{t_1} - f_{t_2}} \le \omega(|t_1 - t_2|) \quad \text{for every } t_1, t_2 \in [a,b].
			\end{equation} 
			Let us choose $D_1 \subset \Omega^2$ of full measure such that
			\begin{equation}
				\label{eq:unicont_pwse}
				\sup_{(x,y) \in D_1} |f_{t_1}(x,y) - f_{t_2}(x,y)| \le \omega(|t_1 - t_2|)  \text{ for every rational } t_1, t_2 \in [a,b].
			\end{equation}
			For $(x,y) \in D_1$, the function $t \mapsto f_t(x,y)$ is continuous on $[a,b]\cap \Q$, 
			so we can extend it to a
			continuous function $t \mapsto \tilde f_t(x,y)$ defined on the whole interval $[a,b]$. For $(x,y) \in \Omega^2 \sm D_1$ and $t \in [a,b]$, we define $\tilde f_t(x,y) = 0$, say. We claim that $\Linf{f_t - \tilde f_t} = 0$ for every $t \in [a,b]$. By the triangle inequality, for any rational $\tau \in [a,b]$,
			\begin{align*}
				\Linf{f_t - \tilde f_t} &\le \Linf{f_t - f_\tau} + \Linf{f_\tau - \tilde f_\tau} + \Linf{\tilde f_\tau - \tilde f_t} \\
				&\le |t - \tau| + 0 + |\tau - t| = 2|t - \tau|.
			\end{align*}
			By sending $\tau$ to $t$ the claim follows, implying that for every $t \in [a,b]$, $f_t = \tilde f_t$ almost everywhere. 
			
			We have $\int f_t \D t = \int \tilde f_t \D t$, since $f_t$ and $\tilde f_t$ are the same element of $\mathcal{E}$. We declare $\tilde f_t$ to be the claimed modification of~$f_t$, for every $t$. For convenience let us further drop the tilde and write $f_t$ instead of $\tilde f_t$. By construction, we have that 
			\begin{equation}
				\label{eq:unicont_everywhere}
				\forall x,y \in \Omega \; \forall t_1,t_2 \in [a,b] \quad |f_{t_1}(x,y) - f_{t_2}(x,y)| \le \omega(|t_1 - t_2|).
			\end{equation}
			
			In particular, $t \mapsto f_t(x,y)$ is continuous, so $J(x,y) := \int_a^b f_t(x,y) \D t$ is well defined for \emph{every} $(x,y) \in \Omega^2$ (on $\Omega^2 \setminus D_1$ we have $J(x,y) = 0$). It remains to show that $\Linf{I-J} = 0$.
			
			Let $I$ be some bounded measurable function $I : \Omega^2 \to \R$ that represents $\int_a^b f_t \D t$.
			Let $f^{(n)}_t$ be a $n$-step function that approximates $f_t$ (recalling the convention we made above), let $I_n(x,y) := (\int f^{(n)}_t \D t) (x,y) = \int f^{(n)}_t(x,y) \D t$, where the equality follows from the definition of the integral for the step function $f^{(n)}_t$. 
			
			By \eqref{eq:norm_of_int} and \eqref{eq:unicont_norm} we get
			\begin{align}
				\label{eq:step_appr}
				\begin{split}
					\Linf{I_n - I} &\le \sum_{i \in [n]} \Linf{\int_{P_i} (f^{(n)}_t - f_t) \D t} \\
					&\le \sum_{i \in [n]} |P_i| \sup_{t \in P_i} \Linf{f^{(n)}_t - f_t} \\
					&\le \sum_{i \in [n]} |P_i| \omega\left( |P_i| \right) = (b-a) \omega\left( \frac{b-a}{n} \right) =: d_n.
				\end{split}
			\end{align}
			We define a set
			\[
			D_2 := \left\{ (x,y) \in \Omega^2 : \forall n \in \N \  \left| I_n(x,y) - I (x,y) \right| \le d_n \right\},
			\]
			which by \eqref{eq:step_appr} has full measure.
			
			Using \eqref{eq:unicont_everywhere}, we have, for every $x, y \in \Omega$,
			\begin{align*}
				\left| I_n(x,y) - J(x,y) \right| &= \left|\int_a^b [f^{(n)}_t(x,y) - f_t(x,y)] \D t\right| \\
				&\le \sum_{i \in [n]} |P_i| \sup_{t \in P_i} |f^{(n)}_t(x,y) - f_t(x,y)| \\
				&\le (b-a) \omega\left( \frac{b-a}{n} \right) = d_n.
			\end{align*}
			We get that for $(x,y) \in D_2$
			\begin{align*}
				&|I(x,y) - J(x,y)| \\
				&\le |I(x,y) - I_n(x,y)| + |I_n(x,y) - J(x,y)| \\
				&\le  2d_n.
			\end{align*}
			Since $d_n \to 0$ as $n \to \infty$, we obtain $I = J$ on $D_2$ and thus $I = J$ as elements of $\mathcal{E}$.
			
			To prove the ``furthermore'' part, we modify the set $D_1$ we chose at the beginning of the proof (still maintaining that it has full measure) so that \eqref{eq:unicont_pwse} is satisfied for both $f_\cdot$ and $g_\cdot$, while, in addition, for every rational $t$ we have $f_t(x,y) \ge g_t(x,y)$ for $(x,y) \in D_1$. Then we notice that the continuous extensions of $t \mapsto f_t(x,y)$ and $t \mapsto g_t(x,y)$ satisfy $\tilde f_t(x,y) \ge \tilde g_t(x,y)$ for $(x,y) \in D_1$ and \emph{every} $t$. Since we defined $\tilde f_t = \tilde{g}_t = 0$ outside of $D_1$, the inequality everywhere follows.
		\end{proof}
	\end{appendix}
	
 \section{Acknowledgements}
		We thank Hanka \v Rada, Joel Spencer and Lutz Warnke for their comments. We also thank an anonymous referee for their feedback.
	\bibliographystyle{plain} 
	\bibliography{DSG}       

\end{document}